\documentclass[12pt,a4paper,reqno]{amsart}
\usepackage{amsmath}
\usepackage{amsfonts}
\usepackage{amssymb}
\numberwithin{equation}{section}

\theoremstyle{plain}

\newtheorem{theorem}[subsection]{Theorem}
\newtheorem{proposition}[subsection]{Proposition}
\newtheorem{lemma}[subsection]{Lemma}
\newtheorem{corollary}[subsection]{Corollary}
\newtheorem{conjecture}[subsection]{Conjecture}

\newtheorem*{quote-thm-1}{Proposition 8.1 of \cite{green-tao-longprimeaps}}
\newtheorem*{mainthm-restate}{Theorem \ref{mainthm-dual}}
\newtheorem*{convex-lemma-restate}{Lemma \ref{convex-props}}
\newtheorem*{decomp-prop-restate}{Proposition \ref{decomp-prop}}
\newtheorem*{main-theorem}{Main Theorem}
\newtheorem*{gvn-restate}{Proposition \ref{gvn}}
\newtheorem*{gvn-restate-again}{Proposition $\mbox{\ref{gvn}}^{\prime}$}
\newtheorem*{gvn-restate-yetagain}{Proposition $\mbox{\ref{gvn}}^{\prime \prime}$}
\newtheorem*{pseudodom-repeat}{Proposition \ref{pseudodom}}

\theoremstyle{definition}

\newtheorem{definition}[subsection]{Definition}

\theoremstyle{remark}

\newtheorem{example}{Example}
\newtheorem{examples}{Examples}
\newtheorem*{remark}{Remark}
\newtheorem*{remarks}{Remarks}

\renewcommand{\leq}{\leqslant}
\renewcommand{\geq}{\geqslant}

\newsavebox{\proofbox}
\savebox{\proofbox}{\begin{picture}(7,7)%
  \put(0,0){\framebox(7,7){}}\end{picture}}

\newcommand{\md}[1]{\ensuremath{(\operatorname{mod}\, #1)}}
\newcommand{\mdsub}[1]{\ensuremath{(\operatorname{mod}\, #1)}}

\newcommand\vol{\operatorname{vol}}
\newcommand\lcm{\operatorname{lcm}}
\newcommand\codim{\operatorname{codim}}
\newcommand\dist{\operatorname{dist}}
\newcommand\GI{\operatorname{GI}}
\newcommand\MN{\operatorname{MN}}
\newcommand\vertical{\operatorname{vertical}}
\def\endproof{\hfill{\usebox{\proofbox}}\vspace{11pt}}
\def\E{\mathbb{E}}
\def\N{\mathbb{N}}
\def\Z{\mathbb{Z}}
\def\R{\mathbb{R}}

\def\C{\mathbb{C}}

\def\Q{\mathbb{Q}}
\def\b{\mathbf{b}}

\def\g{\mathbf{g}}

\def\B{\mathcal{B}}
\newcommand\HK{\operatorname{HK}}
\newcommand\id{\operatorname{id}}

\def\Lip{\operatorname{Lip}}

\renewcommand\th{\operatorname{th}}
\newcommand\eps{\varepsilon}
\newcommand\mobname{M\"obius and nilsequences\ }

\begin{document}

\title[Linear equations in primes]{Linear equations in primes}

\author{Ben Green}
\address{Centre for Mathematical Sciences, Wilberforce Road, Cambridge CB3 0WA, England}
\email{b.j.green@dpmms.cam.ac.uk}

\author{Terence Tao}
\address{Department of Mathematics, UCLA, Los Angeles CA 90095-1555, USA.
}
\email{tao@math.ucla.edu}

\thanks{While this work was carried out the first author was a Clay Research Fellow, and is pleased to acknowledge the support of the Clay Mathematics Institute. Some of this work was carried out while he was on a long-term visit to MIT. The second author was supported by a grant from the Packard Foundation.}

\begin{abstract}  Consider a system $\Psi$ of non-constant affine-linear forms $\psi_1,\ldots,\psi_t: \Z^d \to \Z$, no two of which are linearly dependent.  Let $N$ be a large integer, and let $K \subseteq  [-N,N]^d$ be convex.  A generalisation of a famous and difficult open conjecture of Hardy and Littlewood predicts an asymptotic, as $N \to \infty$, for the number of 
integer points $ n \in \Z^d \cap K$ for which the integers $\psi_1( n),\ldots,\psi_t( n)$ are simultaneously prime.  This implies many other well-known conjectures, such as the twin prime conjecture and the (weak) Goldbach conjecture.  It also allows one to count the number of solutions in a convex range to any simultaneous linear system of equations, in which all unknowns are required to be prime.

In this paper we (conditionally) verify this asymptotic under the assumption that no two of the affine-linear forms $\psi_1,\ldots,\psi_t$ are affinely related; this excludes the important ``binary'' cases such as the twin prime or Goldbach conjectures, but does allow one to count ``non-degenerate'' configurations such as arithmetic progressions.  Our result assumes two families of conjectures, which we term the \emph{inverse Gowers-norm conjecture} ($\GI(s)$) and the \emph{\mobname conjecture} ($\MN(s)$), where $s \in \{1,2,\dots\}$ is the \emph{complexity} of the system and measures the extent to which the forms $\psi_i$ depend on each other.  The case $s=0$ is somewhat degenerate, and follows from the prime number theorem in APs.

Roughly speaking, the inverse Gowers-norm conjecture $\GI(s)$ asserts the Gowers $U^{s+1}$-norm of a function $f : [N] \rightarrow [-1,1]$ is large if and only if $f$ correlates with an $s$-step nilsequence, while the \mobname conjecture $\MN(s)$ asserts that the M\"obius function $\mu$ is strongly asymptotically orthogonal to $s$-step nilsequences of a fixed complexity.  These conjectures have long been known to be true for $s=1$ (essentially by work of Hardy-Littlewood and Vinogradov), and were established for $s=2$ in two papers of the authors.  Thus our results in the case of complexity $s \leq 2$ are unconditional.

In particular we can obtain the expected asymptotics for the number of $4$-term progressions $p_1 < p_2 < p_3 < p_4 \leq N$ of primes, and more generally for any (non-degenerate) problem involving two linear equations in four prime unknowns.
\end{abstract}

\setcounter{tocdepth}{1}

\maketitle

\tableofcontents

\section{Introduction}\label{intro-sec}

\textsc{A Generalised Hardy-Littlewood Conjecture.} Let $P := \{2,3,5,\ldots\} \subset \Z$ denote the prime numbers.  We refer to the lattice points $(p_1,\ldots,p_t) \in P^t$ as \emph{prime points} in $\Z^t$.  A basic problem in additive number theory is to count the number of prime points on a given affine sublattice of $\Z^t$ in a given range.  For instance, the twin prime conjecture asserts that the number of prime points in
$\{ (n,n+2): n \in \Z \} \subset \Z^2$ is infinite.  When the affine lattice is formed by intersecting $\Z^t$ with an affine subspace, this problem is equivalent to finding solutions to simultaneous linear equations in which all unknowns are prime.
To formalise these types of problems more concretely, it is convenient to parameterise this lattice by $d$ affine-linear forms, as follows.

\begin{definition}[Affine-linear forms] Let $d,t \geq 1$ be integers.  An \emph{affine-linear form} on $\Z^d$ is a function $\psi: \Z^d \to \Z$ which is the sum $\psi = \dot \psi + \psi(0)$ of a linear form $\dot \psi: \Z^d \to \Z$ and a constant $\psi(0) \in \Z$.  A \emph{system of affine-linear forms} on $\Z^d$ is a collection $\Psi = (\psi_1,\ldots,\psi_t)$ of affine-linear forms on $\Z^d$.  To avoid trivial degeneracies we shall require that all the affine-linear forms are non-constant and no two forms are rational multiples of each other.  The entire system $\Psi$ can be thought of as an affine-linear map from $\Z^d$ to $\Z^t$, which is the sum $\Psi = \dot \Psi + \Psi(0)$ of a linear map $\dot \Psi: \Z^d \to \Z^t$ and a constant $\Psi(0) \in \Z^t$; we refer to the range $\Psi(\Z^d)$ of this map as an \emph{affine sublattice} of $\Z^t$.  We extend $\Psi$ (and $\dot \Psi$) in the obvious manner to an affine-linear map from $\R^d$ to $\R^t$.
If $N > 0$, we define the \emph{size} $\|\Psi\|_N$ of $\Psi$ relative to the scale $N$ to be the quantity
\begin{equation}\label{veci}
\| \Psi\|_N := \sum_{i=1}^t \sum_{j=1}^d |\dot \psi_i(e_j)| + \sum_{i=1}^t \left|\frac{\psi_i(0)}{N}\right|
\end{equation}
where $e_1,\ldots,e_d$ is the standard basis for $\Z^d$.
\end{definition}

\begin{example} The line $\{ (n,n+2): n \in \Z \}$ is the affine lattice associated to the system $\Psi: n \mapsto (n, n+2)$ with $d=1$ and $t=2$.  This example has bounded size for any $N \geq 1$.  The system $\Psi: n \mapsto (n,N-n)$ counts pairs of primes which sum to $N$, and has bounded size at scale $N$.
\end{example}

In order to count the number of prime points on an affine lattice, it is convenient to use the \emph{von Mangoldt function} $\Lambda: \Z \to \R^+$, defined by setting $\Lambda(n) := \log p$ when $n > 1$ is a power of a prime $p$, and $\Lambda(n) = 0$ otherwise (in particular, $\Lambda(n) = 0$ whenever $n \leq 0$).  We are then interested in estimating the sum
\begin{equation}\label{uk}
\sum_{ n \in K \cap \Z^d} \prod_{i \in [t]} \Lambda( \psi_i( n) )
\end{equation}
where $K$ is a convex subset of $\R^d$ and $[t] := \{1,\ldots,t\}$.

\begin{remark} We do not necessarily assume that $\Psi$ is injective, that is to say we allow the sum in \eqref{uk} to
count a single prime point repeatedly.  This freedom will be convenient for us at a later stage of the argument when we increase the number $d$ of parameters in order to place $\Psi$ in a certain normal form.  However, in most applications of interest it will indeed be the case that $\Psi$ is injective, and so the prime points are counted without multiplicity.
\end{remark}

The prime number theorem asserts that the average value of $\Lambda(n)$ is $1$ for positive $n$ and $0$ for negative $n$, so it is first natural (cf. Cramer's model for the primes) to consider the much simpler sum
$$ \sum_{ n \in K \cap \Z^d} \prod_{i \in [t]} 1_{\R^+}(\psi_i( n)) $$
where we use $1_E$ to denote the indicator of a set $E$ (thus $1_E(x)=1$ when $x \in E$ and $1_E(x)=0$ otherwise).  Let us assume that the convex body $K$ is contained in the box $[-N,N]^d$ for some large integer $N$, and let us also assume the size bounds $\|\Psi\|_N \leq L$ for some $L > 0$.  Then a simple volume packing argument (see Appendix \ref{convexgeom}) yields the asymptotic
\begin{equation}\label{beta-gauss}
\sum_{ n \in K \cap \Z^d} \prod_{i \in [t]} 1_{\R^+}(\psi_i( n)) = 
\beta_\infty + O_{d,t,L}(N^{d-1}) = \beta_\infty + o_{d,t,L}(N^d)
\end{equation}
where the \emph{archimedean factor} $\beta_\infty$ is defined by
\begin{equation}\label{beta-inf} \beta_{\infty} := \vol_d\big( K \cap \Psi^{-1}((\R^+)^t) \big)\end{equation}
(see \S \ref{notation-sec} for our conventions concerning asymptotic notation).
Note that the main term $\beta_\infty$ is typically of size $N^d$ or so.  
One can be much more precise about the nature of the error term, but we will not be concerned with quantitative decay rates here. Indeed the rates provided by our later arguments will be poor and often ineffective, and will dominate whatever gains one could extract from the error term in \eqref{beta-gauss}.

In view of \eqref{beta-gauss} and the prime number theorem, one might na\"{\i}vely conjecture that the expression \eqref{uk}
also enjoys the asymptotic $\beta_\infty + o_{d,t,L}(N^d)$.  However this is not the case due to local obstructions at small moduli.
For instance, we have
\begin{equation}\label{siegel}
\sum_{n=1}^N \Lambda( qn + b ) = \Lambda_{\Z_q}(b) N + o_{q}(N) 
\end{equation}
whenever $q \geq 1$ and $|b| \leq q$, where $\Lambda_{\Z_q}: \Z \to \R^+$ is the \emph{local von Mangoldt function}, that is the $q$-periodic function defined by setting $\Lambda_{\Z_q}(b) := \frac{q}{\phi(q)}$ when $b$ is coprime to $q$ and $\Lambda_{\Z_q}(b) = 0$ otherwise.
  Here $\Z_q := \Z/q\Z$ is the cyclic group of order $q$ and $\phi(q) := |\Z_q^\times|$ is the Euler totient function.  We shall refer to \eqref{siegel} as the \emph{prime number theorem in APs}. A well-known quantitative version of this result is the Siegel-Walfisz theorem, which establishes the asymptotic \eqref{siegel} uniformly in the range $q \leq \log^A N$ for any fixed $A$. In this range, the $o$-term is ineffective, and if one wishes for an effective error term it is necessary to restrict to $q \leq \log^{1-\delta} N$ for some $\delta > 0$. See \cite[p. 123]{davenport-book} for details.
  
More generally, given a system $\Psi = (\psi_1,\ldots,\psi_t)$ of affine-linear forms, one can define the \emph{local factor} $\beta_q$ for any integer $q \geq 1$ by the formula
\begin{equation}\label{beta-p} \beta_q := \E_{{n} \in \Z_q^d} \prod_{i \in [t]} \Lambda_{\Z_q}(\psi_i({n})).\end{equation}
The symbol $\E$ denotes expectation or averaging; see \S \ref{notation-sec} for more details.
From the Chinese remainder theorem we see that this factor is multiplicative, indeed we have $\beta_q = \prod_{p|q} \beta_p$, where the product is over all primes\footnote{More generally, we adopt the convention that whenever a product ranges over $p$, that $p$ is understood to be restricted to the primes.} $p$ dividing $q$.  We then have

\begin{conjecture}[Generalised Hardy-Littlewood conjecture]\label{dickson}  Let $N, d, t, L$ be positive integers, and let $\Psi = (\psi_1,\ldots,\psi_t)$ be a system of affine-linear forms with size $\| \Psi\|_N \leq L$.  
Let $K \subset [-N,N]^d$ be a convex body.  Then we have
\begin{equation}\label{dickson-eq}
\sum_{ n \in K \cap \Z^d} \prod_{i \in [t]} \Lambda( \psi_i( n) ) = \beta_\infty \prod_p \beta_p + o_{t,d,L}(N^d)
\end{equation}
where the archimedean factor $\beta_\infty$ and the local factors $\beta_p$ for each prime $p$ were defined in \eqref{beta-inf}, \eqref{beta-p}.  
\end{conjecture}

Roughly speaking, this conjecture asserts that $\Lambda$ ``behaves like'' the independent product of $1_{\R^+}$ and $\Lambda_{\Z_p}$, as $p$ ranges over primes.
In typical applications, the quantities $\beta_\infty$ and $\beta_p$ are quite easy to compute explicitly: see Examples \ref{example1}-\ref{example5} below.
We shall refer to the quantity $\prod_p \beta_p$ as the \emph{singular product}.  The local factors $\beta_p$ can be easily estimated:

\begin{lemma}[Local factor bounds]\label{sing} With the hypotheses of Conjecture \ref{dickson}, we have $\beta_p = 1 + O_{t,d,L}( p^{-1} )$.  If furthermore no two of the forms $\psi_1,\ldots,\psi_t$ are affinely related \textup{(}i.e. no two of the forms $\dot \psi_1, \ldots, \dot \psi_t$ are parallel\textup{)}, or if $p > C(t,d,L) N$ for some sufficiently large constant $C(t,d,L)$, then we have $\beta_p = 1 + O_{t,d,L}( p^{-2} )$.
\end{lemma}

\begin{proof}  Without loss of generality we may assume $p$ to be large compared to $t,d,L$, as the claim is trivial otherwise.
Let $ n$ be selected uniformly at random from $\Z_p^d$.  Since the $\psi_i$ are non-constant, we easily see
that $\Lambda_{\Z_p}(\psi_i({n}))$ will equal $\frac{p}{p-1}$ with probability $1 - \frac{1}{p}$, and $0$ otherwise.
In particular the product in \eqref{beta-p} is equal to $(\frac{p}{p-1})^t = 1 + O_t(\frac{1}{p})$ with probability $1 - O_t(\frac{1}{p})$ and zero otherwise, which gives the first bound on $\beta_p$.  Now suppose that either no two of $\psi_1,\ldots,\psi_t$ are affinely related, or that $p > C(t,d,L) N$ for some sufficiently large $C(t,d,L)$.  Then for any $1 \leq i < j \leq t$, we see from elementary linear algebra that $\psi_i({n})$ and $\psi_j({n})$ will simultaneously be divisible by $p$ with probability $O(\frac{1}{p^2})$; the point is that the hypotheses imply that\footnote{One could view this as a (very simple) manifestation of the Lefschetz principle.} $\psi_i$ and $\psi_j$ cannot be linear multiples of each other modulo $p$.  The desired bound on $\beta_p$ then follows from a simple application of the Bonferroni inequalities (that is, the fact that truncations of the inclusion-exclusion formula give upper and lower bounds alternately).
\end{proof}

In particular we see that the singular series $\prod_p \beta_p$ is always convergent (though it could vanish, thanks to the presence of the small primes $p = O_{t,d,L}(1)$).

A straightforward argument shows that Conjecture \ref{dickson} implies a conjecture which counts primes more explicitly:

\begin{conjecture}[Generalised Hardy-Littlewood conjecture, again]\label{dickson-again}  Let $N, d, t, L, \Psi, K$ be
as in Conjecture \ref{dickson}.  Then
\begin{equation}\label{dickson-eq-again}
\begin{split}
|K \cap \Z^d \cap \Psi^{-1}(P^t)| &= \#\{  n \in K \cap \Z^d: \psi_1( n), \ldots, \psi_t( n) \hbox{ prime} \} \\
&=(1 + o_{t,d,L}(1)) \frac{\beta_\infty}{\log^t N} \prod_p \beta_p + o_{t,d,L}\left(\frac{N^d}{\log^t N}\right).
\end{split}
\end{equation}
\end{conjecture}

\begin{remarks} It would be slightly more accurate to replace $\frac{\beta_\infty}{\log^t N}$ with the more precise expression \[ \int_K \prod_{j \in [t]} \frac{1_{\psi_j( x) > 2}}{\log \psi_j( x)}\ d x,\] but the difference between these two expressions can be absorbed into the qualitative $o_{t,d,L}()$ error terms.  In most (though not quite all) cases, the singular series $\prod_p \beta_p$ is bounded by $O_{t,d,L}(1)$, which
allows one to absorb the first error term into the second.  Informally speaking, this conjecture asserts that the probability that a randomly selected point in $\Psi(\Z^d) \cap \Z_+^t$ of magnitude $N$ is a prime point is asymptotically $\frac{1}{\log^t N} \prod_p \beta_p$.
\end{remarks}

\begin{proof}[Sketch proof of Conjecture \ref{dickson-again} assuming Conjecture \ref{dickson}] Let $0 < \eps < 1$ be a small quantity (depending on $N,d,t,L$) to be chosen later.  The
contribution to \eqref{dickson-eq-again} where $\min_{1 \leq i \leq t} |\psi_i( n)| \leq N^{1-\eps}$ can easily be shown
to be $o_{t,d,L,\eps}(N^{d-\eps/2})$ by crude estimates; the analogous contribution to \eqref{dickson-eq} can
similarly be shown to be $o_{t,d,L,\eps}(N^{d})$.  The contribution to \eqref{dickson-eq} where at least one of the $\psi_i( n)$ is a power of a prime $p^2, p^3, \ldots$ can similarly be shown to be $o_{t,d,L}(N^{d})$.  Finally,
for the remaining non-zero contributions to \eqref{dickson-eq}, the quantity
$\prod_{i \in [t]} \Lambda( \psi_i( n) )$ is equal to $(1 + O(t \eps)) \log^t N$.  Putting all this 
together, we see that the left-hand side of \eqref{dickson-eq-again} is
$$
(1 + O(t\eps)) \frac{\beta_\infty}{\log^t N} \prod_p \beta_p + o_{t,d,L,\eps}(\frac{N^d}{\log^t N}).
$$
Setting $\eps$ to be a sufficiently slowly decaying function of $N$ (for fixed $t,d,L$) we obtain the claim.
\end{proof}

Note that the case $d=t=1$ of the generalised Hardy-Littlewood conjecture is essentially the prime number theorem in APs \eqref{siegel}.  We have been referring to the \emph{generalised} Hardy-Littlewood conjecture because Hardy and Littlewood \cite{hardy-littlewood} in fact only conjectured an asympotic for the number of $n \leq N$ for which the forms $n + b_1, \dots, n+ b_t$ are all prime. If this were generalised to deal with the case of forms $a_1n + b_1,\dots, a_t n + b_t$ -- the case $d = 1$ of Conjecture \ref{dickson} -- then a $d$-parameter version along the lines we have been discussing would follow easily by holding $d-1$ of the variables fixed and summing in the remaining one. One has the impression that, had they thought to ask the question, Hardy and Littlewood would easily have produced a conjecture for the asymptotic formula. The name of Dickson is sometimes associated to this circle of ideas. In the 1904 paper \cite{dickson}, he noted the obvious necessary condition on the $a_i,b_i$ in order that the forms $a_1n + b_1,\dots, a_tn + b_t$ might all be prime infinitely often and suggested that this condition might also be sufficient. 

Dickson also suggested that the ``experts in the new Dirichlet theory'' try their hand at establishing this. His hope has yet to be realised, however, since the $d=1$, $t > 1$ case of Conjecture \ref{dickson} seems to be extremely difficult. The twin prime, Sophie Germain, and weak\footnote{That is, the conjecture that every \emph{sufficiently large} even number is the sum of two primes.} even Goldbach conjectures, for instance, follow easily from the $d=1$, $t=2$ case of the conjecture. These cases are probably well beyond the reach of current technology, although we remark that if one replaces the von Mangoldt function $\Lambda$ with substantially simpler weight functions arising from the Selberg $\Lambda^2$ sieve then such asymptotics can be obtained by standard sieve theory methods (see Theorem \ref{fundlemma}).  This
in turn leads to \emph{upper} bounds on \eqref{uk} which differ from \eqref{dickson-eq} only by a multiplicative constant 
depending only on $d,t,L$. 

Note also that it is possible to establish the case $d = 1$, $t > 1$ of the Hardy-Littlewood conjecture \emph{on average} over the choice of forms $\psi_1,\dots,\psi_t$ in a certain sense: see \cite{Balog}.  This essentially amounts to increasing $d$, which can place one back in the ``finite complexity'' regime discussed below.

\textsc{Complexity.} We will not make any progress on the $d=1$, $t > 1$ case here, but instead focus on the substantially simpler cases when $d > 1$ and
the system is ``finite complexity'' in the following sense.

\begin{definition}[Complexity]  Let $\Psi = (\psi_1,\ldots,\psi_t)$ be a system of affine-linear forms.  If $1 \leq i \leq t$ and $s \geq 0$, we say that $\Psi$ has \emph{$i$-complexity at most $s$} if one can cover the $t-1$ forms $\{ \psi_j: j \in [t] \backslash \{i\} \}$ by $s+1$ classes, such that $\psi_i$ does not lie in the affine-linear span of any of these classes.  The \emph{complexity} of the $\Psi$ is defined to be the least $s$ for which the system has $i$-complexity at most $s$ for all $1 \leq i \leq t$, or $\infty$ if no such $s$ exists.
\end{definition}

\begin{remark} It is easy to see that one can replace ``cover $\ldots$ by'' by ``partition $\ldots$ into'' in the above definition without affecting the definition of $i$-complexity or complexity.  While partitions are slightly more natural here than covers, we prefer to use covers as it makes it a little easier to compute the complexity in some cases.
\end{remark}

\begin{examples}\label{simple} The system $\Psi(n_1,\ldots,n_d) := (n_1,\ldots,n_d)$, which counts $d$-tuples of independent primes, has complexity $0$, because no form $u_i$ lies in the affine span of all the other forms.  For any $k \geq 2$, the system $\Psi(n_1,n_2) := (n_1, n_1 + n_2, \ldots, n_1 + (k-1) n_2)$, which counts arithmetic progressions of primes of length $k$, has complexity $k-2$, because each form does not lie in the affine span of any other \emph{individual} form, though it is in the affine span of any two other forms.  The system $\Psi(n_1,n_2) := (n_1,n_2,N-n_1-n_2)$, which counts triples of primes that sum to a fixed number $N$, has complexity $1$.  The system $\Psi(n_1,n_2) := (n_1, n_2, n_1+n_2-1, n_1+2n_2-2)$, which counts progressions of primes of length three, whose difference $n_2-1$ is one less than a prime, has complexity $2$.
The system $\Psi(n_1) := (n_1,n_1+2)$, which counts twin primes, has infinite complexity. So too does the system $\Psi(n_1) := (n_1, N-n_1)$, which counts pairs of primes which sum to a fixed number $N$, as well as $\Psi(n_1) = (n_1,2n_1+1)$, which counts Sophie Germain primes.  More generally, any system with $d=1$ and $t > 1$ has infinite complexity.
\end{examples}

\begin{example}[Cubes]\label{cube-ex}  Let $d \geq 2$ and $t := 2^{d-1}$.  Then the system
\[ \Psi(n_1,\ldots,n_d) := \big( n_1 + \sum_{j \in A} n_j \big)_{A \subseteq  \{2,\ldots,d\}}, \]
(which counts $(d-1)$-dimensional cubes whose vertices are all prime)
has a very large value of $t$, but has complexity at most $d-2$.  For instance, if one considers the form $n_1$, then one can cover the other $t-1$ forms by $d-1$ classes, with the $i^{\th}$ class consisting of those forms which involve $n_{i+1}$, then $n_1$ is not in the affine span of any of these classes because the $i^{\th}$ class always assigns the same coefficient to both $n_1$ and $n_{i+1}$.  The other forms can be treated similarly after ``reflecting'' the cube appropriately.
\end{example}

\begin{example}[$\mbox{IP}_0$ cubes]\label{ipex} Let $d \geq 1$ and $t := 2^d-1$.  Then the system
\[ \Psi(n_1,\ldots,n_d) := \big( 1 + \sum_{j \in A} n_j \big)_{A \subseteq  [d]; A \neq \emptyset},\]
which counts $d$-dimensional cubes pinned at the origin whose remaining vertices are one less than a prime,
also has a large value of $t$ but has complexity at most $d-1$, for reasons similar to the previous example.
\end{example}

In fact in Example \ref{cube-ex} the complexity is \emph{exactly} $d-2$, whilst in Example \ref{ipex} it is \emph{exactly} $d-1$. We leave the proofs to the reader.

\begin{example}[Balog's example]\label{balog-ex} Let $d \geq 2$ and $t := \frac{d(d+1)}{2}$.  Then the system \[ \Psi(n_1,\ldots,n_d) := (n_i+n_j+1)_{1 \leq i \leq j \leq d},\] which counts $d$-tuples of odd primes $p_1,\ldots,p_d$, all of whose midpoints $\frac{p_i+p_j}{2}$ are also prime, has complexity $1$, even though $t$ is quite large.  Indeed, if
one considers the form $n_i+n_j+1$ with $i < j$, one can partition the other $t-1$ forms into two classes, those which do not involve $n_i$, and those which do involve $n_i$ (and hence do not involve $n_j$), and $n_i+n_j+1$ is an affine-linear combination of neither of these two classes.  If instead one considers the form $n_i+n_i+1=2n_i+1$, one can partition the other $t-1$ forms into two classes, those which involve $n_i$ (and one other $n_j$), and those which do not involve $n_i$ at all, and again $2n_i+1$ is an affine-linear combination of neither of these two classes.
\end{example}

The complexity is a little difficult to compute directly, but the following lemma gives some easy bounds on this quantity.

\begin{lemma}[Complexity bounded by codimension]\label{complex}   Let $\Psi = (\psi_1,\ldots,\psi_t)$ be a system of affine-linear forms.  Then this system has finite complexity if and only if no two of the $\psi_i$ are affinely dependent.  Furthermore, in this case the complexity of the system is less than or equal to $t - \dim( \dot \Psi )$.
\end{lemma}

\begin{proof} If two of the forms $\psi_i$ and $\psi_j$ are affinely related, then it is not possible for the $i$-complexity to be finite, as $\psi_i$ will lie in the affine span of any collection of forms which contain $\psi_j$.  Conversely, if no two of the $\psi_i$ are affinely related, then the $i$-complexity is at most $t-2$, as we can partition the $t-1$ forms $\{ \psi_j: j \in [t] \backslash \{i\}\}$ into singletons.  This gives the first claim of the lemma.

Now suppose that no two of the $\psi_i$ are affinely dependent. Write $r := \dim(\dot \Psi)$. Choose any homogeneous form, say $\dot \psi_1$; this will be nonzero. Relabelling if necessary, we may suppose that $\{\dot \psi_1,\dots, \dot \psi_r\}$ is a basis for $\dot \Psi$. Consider the set $\{ \psi_2,\dots, \psi_r\}$ along with the singleton sets $\{ \psi_{r+1}\},\dots,\{ \psi_t\}$. Clearly $\psi_1$ is not in the affine-linear span of any such set, and so the system has $1$-complexity at most $t-r$. Since this is true with any $\psi_i$ in place of $\psi_1$, the claim follows.
\end{proof}

\begin{remark}
This lemma is sharp in all the cases treated in Examples \ref{simple}, but is very far from sharp in Examples \ref{cube-ex}-\ref{balog-ex}.  It asserts that the infinite complexity systems are precisely those which encode a ``binary'' problem such as the twin prime, Goldbach, Sophie Germain, or prime tuples conjectures.  Observe from Lemma \ref{complex} and Lemma \ref{sing} that if the system has finite
complexity, then $\beta_p = 1 + O_{t,d,L}(\frac{1}{p^2})$ and so the singular series $\prod_p \beta_p$ is either zero, or is bounded above and below by constants depending only on $t,d,L$.  In particular we can eliminate the first error term in \eqref{dickson-eq-again} in this setting.
\end{remark}

For systems of complexity $0$, The generalised Hardy-Littlewood conjecture follows easily from the prime number theorem in APs \eqref{siegel}.  For systems of complexity $1$, the conjecture  can be treated by the Hardy-Littlewood circle method (see e.g. \cite{Balog,balog-cite}).  Systems of complexity $2$ or higher, on the other hand, are largely out of reach of the circle method and the conjecture has remained open in these cases.

We mention two directions in which a partial approach to high complexity cases of the generalised Hardy-Littlewood conjecture has been made. The first is that a version of the conjecture remains true if one is willing to enlarge sufficiently many of the $\Lambda$ factors, replacing primes with some notion of an \emph{almost prime}, and adjust the singular series appropriately; see for instance Theorem \ref{fundlemma} for a simplified version of this result.  One consequence of this is that \emph{upper bounds} in \eqref{dickson-eq} (or \eqref{dickson-eq-again}) are known which are only off by a multiplicative constant of $O_{t,d,L}(1)$.  

For certain special systems a \emph{lower bound} of the correct order of magnitude is available.
For some systems such as the cube systems in Example \ref{cube-ex} this is rather simple, involving nothing more than a few applications of the Cauchy-Schwarz inequality, despite the fact that such systems can have arbitrarily high complexity.  However, the task of obtaining asymptotics here is just as difficult as obtaining asymptotics for other systems; see \cite{host-kra} for some related discussion of this phenomenon.

There is also the system $\Psi(n_1,n_2) := (n_1, n_1 + n_2, \ldots, n_1 + (k-1) n_2)$ of arithmetic progressions of length $k$, for which the powerful tool of \emph{Szemer\'edi's theorem} \cite{szemeredi} was available.  Despite the fact that these systems can have arbitrarily high complexity, a \emph{lower bound} for \eqref{dickson-eq} and \eqref{dickson-eq-again} was established which was again only off by a multiplicative constant. In particular this implied that the primes contain arbitrarily long arithmetic progressions; see \cite{green-tao-longprimeaps}.  

Our arguments in this paper borrow many ideas and results from \cite{green-tao-longprimeaps}, in particular drawing heavily on the \emph{transference principle} developed in that paper. However we shall not use Szemer\'edi's theorem in this paper, as it does not apply to the general systems of affine-linear forms studied here. Roughly speaking, one only expects Szemer\'edi-type theorems for systems which are \emph{homogeneous} (so $\Psi(0)=0$) and \emph{translation invariant}, that is the lattice $\dot \Psi(\Z^d)$ contains the diagonal generator $(1,\ldots,1)$. In any case Szemer\'edi's theorem only provides lower bounds and not asymptotics.
  
\textsc{Main result.} Our main result settles the generalised Hardy-Littlewood conjecture for any system of affine-linear forms of finite complexity, 
conditional on two simpler, partially resolved, conjectures.

\begin{main-theorem}[Generalised Hardy-Littlewood conjecture, finite complexity case]\label{main}  Suppose that the inverse Gowers-norm conjecture $\GI(s)$ and the \mobname conjecture $\MN(s)$ are true for some finite $s \geq 1$. Both of these conjectures will be stated formally in \S \ref{gowers-sec}. Then the generalised Hardy-Littlewood conjecture is true for all systems of affine-linear forms of complexity at most $s$.
\end{main-theorem}

We have deferred the precise statement of the conjectures $\GI(s)$ and $\MN(s)$ to \S \ref{gowers-sec} on account of the fact that both of them are somewhat technical to state formally. The impatient reader may wish to jump to that section to view these conjectures, but for now we settle for informal one-line statements of them.

The inverse Gowers-norm conjecture $\GI(s)$ gives an explicit criterion as to when a bounded sequence of complex numbers is ``Gowers uniform of order $s$'', this being a measure of pseudorandomness of the sequence; namely, this Gowers uniformity holds whenever the sequence fails to be correlated with any $s$-step nilsequence.  

The \mobname conjecture $\MN(s)$ asserts that the M\"obius function $\mu(n)$ (which is of course closely related to $\Lambda(n)$) does indeed have negligible correlation with all $s$-step nilsequences.  

Neither of these two conjectures are fully resolved at present.  However, the case $s = 1$ is classical and was essentially already present in the work of Hardy-Littlewood and Vinogradov, though not in this language.  The conjecture $\GI(2)$ was settled more recently in \cite{green-tao-u3inverse}, while the conjecture $\MN(2)$ was settled in \cite{green-tao-u3mobius}.  Because of this, we have the following unconditional result:

\begin{corollary}\label{maincor}  The generalised Hardy-Littlewood conjecture is true for all systems of affine-linear forms of complexity at most $2$.
In particular, thanks to Lemma \ref{complex}, the generalised Hardy-Littlewood conjecture is true for any system $\Psi = (\psi_1,\ldots,\psi_t)$ in which no two $\psi_i,\psi_j$ are affinely dependent, and such that $\codim(\dot \Psi(\R^d)) \leq 2$.
\end{corollary}

We expect both $\GI(s)$ and $\MN(s)$ to be settled shortly for general $s$, and hope to report on progress on both of these conjectures in the not-too-distant future\footnote{Note added in April 2008: in a recent preprint, the authors have fully resolved the $\MN(s)$ conjecture for every $s$.}.  We therefore expect to settle the generalised Hardy-Littlewood conjecture entirely in the finite complexity case, or in other words we should be able to remove the last hypothesis in Corollary \ref{maincor}.  The only unresolved case of the generalised Hardy-Littlewood conjecture would then be the presumably very hard ``binary'' or ``infinite complexity'' case in which two or more of the forms are affinely related.

Let us now state some particular new consequences of our results.  The first three are unconditional, while the last two require further progress on the inverse Gowers-norm and \mobname conjectures.

\begin{example}[APs of length $4$]\label{example1} The number of $4$-tuples of primes $p_1 < p_2 < p_3 < p_4 \leq N$ which lie in arithmetic progression is $(1 + o(1)) \mathfrak{S}_1 \frac{N^2}{\log^{4} N}$, where 
\[ \mathfrak{S}_1 := \frac{3}{4} \prod_{p \geq 5} \big( 1 - \frac{3p-1}{(p-1)^3}\big) \approx 0.4764.\]
This follows from Corollary \ref{maincor} with the system
$\Psi(n_1,n_2) := (n_1,n_1+n_2,n_1+2n_2,n_1+3n_2)$, with $K$ being the convex region $\{ (n_1,n_2): 1 \leq n_1 \leq n_1 + 3 n_2 \leq N \}$; one has
$\beta_\infty = N^2/6$, $\beta_2 = 4$, $\beta_3 = 9/8$, and $\beta_p = 1 - \frac{3p-1}{(p-1)^3}$ for $p \geq 5$.
Note that the results in \cite{green-tao-longprimeaps} do not give this asymptotic, instead yielding a \emph{lower} bound of $(c + o(1)) \frac{N^2}{\log^4 N}$
for some explicitly computable but rather small constant $c > 0$.
\end{example}

\begin{example}[APs of length $3$ with common difference $p \pm 1$]\label{example2} The number of triples of primes $p_1 < p_2 < p_3 \leq N$ in arithmetic progression, in which the common difference $p_2 - p_1$ is equal to a prime plus 1, is  $(1 + o(1)) \mathfrak{S}_2 N^2\log^{-4} N$, where
\[ \mathfrak{S}_2 := \prod_{p \geq 3} \big( 1 -  \frac{p^2 - 4p + 1}{(p-1)^4}  \big) \approx 1.0481.\]
The same asymptotic holds for progressions in which $p_2 - p_1$ is a prime minus 1.  This follows from a similar application of Corollary \ref{maincor} as in Example \ref{example1}.  
\end{example}

\begin{example}[Vinogradov $3$-primes theorem with a constraint]\label{example3} Let $N$ be a large odd integer. Then the number of distinct representations of $N$ as $p_1 + p_2 + p_3$ in which $p_1 - p_2$ is equal to a prime minus 1 is 
$(\mathfrak{S}_3(N) + o(1)) \frac{N^2}{\log^{4} N}$,  where 
\[ \mathfrak{S}_3(N) := \frac{1}{3} \prod_{\substack{p \geq 3\\ p | N^3 - N}} \big( 1 - \frac{p^2 - 4p + 1}{(p-1)^4}\big)\prod_{\substack{p \geq 3\\ p \nmid N^3 - N}} \big( 1 + \frac{4p - 1}{(p-1)^4}\big).\] 
Thanks to Lemma \ref{sing}, we see that $\mathfrak{S}_3(N)$ is bounded above and below by absolute positive constants independently of $N$.  Again, this result follows from a specific application of Corollary \ref{maincor}.
\end{example}

\begin{example}[APs of length $k$]\label{example4} Let $k \geq 2$ be a fixed integer. Assume the $\GI(k-2)$ conjecture and the $\MN(k-2)$ conjecture. Then the number of $k$-tuples of primes $p_1 < p_2 < \dots < p_k \leq N$ which lie in arithmetic progression is 
$$\left(\frac{1}{2(k-1)} \prod_p \beta_p + o_{k}(1)\right) \frac{N^2}{\log^{k} N}$$
where
\[ \beta_p := \left\{ \begin{array}{ll} \frac{1}{p} \left(\frac{p}{p-1}\right)^{k-1} & \mbox{if $p \leq k$} \\ \left(1 - \frac{k-1}{p}\right) \left( \frac{p}{p-1}\right)^{k-1} & \mbox{if $p \geq k$}.\end{array}\right.\]
The $k=4$ case of this is Example \ref{example1}; the $k=3$ case is due to van der Corput \cite{van-der-corput}; and the $k=1,2$ cases are equivalent to the prime number theorem.
For comparison, the arguments in \cite{green-tao-longprimeaps} give an unconditional lower bound of $(c_k + o(1)) \frac{N^2}{\log^k N}$ for some $c_k > 0$.
\end{example}

\begin{example}[$P-1$ and $P + 1$ are $\mbox{IP}_0$-sets]\label{example5} Assume $s \geq 0$ is such that the $\GI(s)$ and $\MN(s)$ conjectures are true. Then (thanks to Example \ref{ipex}) there exist infinitely many $s+1$-tuples $(n_1,\ldots,n_{s+1})$ of distinct positive integers such that all of the sums $\{ \sum_{i \in A} n_i : A \subseteq  [s+1], A \neq \emptyset\}$, are equal to a prime minus $1$. Similarly for the primes plus $1$.  In particular, we unconditionally have the new result that there are infinitely many distinct $n_1,n_2,n_3$ such that $n_1,n_2,n_3,n_1+n_2,n_1+n_3,n_2+n_3,n_1+n_2+n_3$ are all one less than a prime.
\end{example}

Another consequence of the Main Theorem concerns counting the number of solutions in a given range to a system of linear equations, in which all unknowns are required to be prime:

\begin{theorem}[Linear equations in primes]\label{main2} Assume the $\GI(s)$ and $\MN(s)$ conjectures. Let $A = (a_{ij})$ be an $s \times t$ matrix of integers, where $s \leq t$. Assume the \emph{non-degeneracy conditions} that $A$ has full rank $s$, and that the only element of the row-space of $A$ over $\mathbb{Q}$ with two or fewer non-zero entries is the zero vector.   Let $N > 1$, let $ b = (b_1,\ldots,b_s) \in \Z^s$ be a vector in $A \Z^t = \{ A {x}: {x} \in \Z^t\}$, and suppose that the coefficients $|a_{ij}|$ and the quantities $|b_i/N|$ are uniformly bounded by some constant $L$. Let $K \subseteq  [-N,N]^t$ be convex. Then we have
\begin{equation}\label{ziggy}
 \sum_{\substack{{x} \in K \cap \Z^t \\ A {x} = {b}}} \prod_{i \in [t]} \Lambda(x_i) = \alpha_{\infty} \prod_p \alpha_p + o_{t,L,s}(N^{t-s}),
\end{equation}
where the \emph{local densities} $\alpha_p$ are given by
\begin{equation}\label{alpha-p} \alpha_p := \lim_{M \to \infty} \E_{x \in [-M,M]^t, A {x} = {b}} \prod_{i \in [t]} \Lambda_{\Z_p}(x_i)
\end{equation} 
and the \emph{global factor} $\alpha_{\infty}$ is given by
\begin{equation}\label{alpha-inf} \alpha_{\infty} := \#  \{ {x} \in \Z^t : {x} \in K, A {x} = {b} , x_i \geq 0\} .\end{equation}
\end{theorem}

Theorem \ref{main2} follows easily from the Main Theorem and some elementary linear algebra: the details may be found in \S \ref{linalg}.  The quantities $\alpha_p$ and $\alpha_\infty$ can be easily computed in practice.
One can also formulate an analogue of Theorem \ref{main2} which counts prime solutions to $A {x} = {b}$, just as Conjecture \ref{dickson-again} could be deduced from Conjecture \ref{dickson}. We leave the details to the reader.
Theorem \ref{main2} is not the most general consequence of the Main Theorem, but it is rather representative. For instance, it already implies Examples \ref{example1}--\ref{example4} (and also implies Example \ref{example5} if $\GI(s)$ and $\MN(s)$ are known for all $s$).

Another simple ``qualitative'' consequence of the Main Theorem is the following.

\begin{corollary}[Qualitative generalised H-L conjecture for finite complexity systems] Suppose that $\GI(s)$ and $\MN(s)$ are true for some $s \geq 1$.  Let $\Psi = (\psi_1,\ldots,\psi_t): \Z^d \to \Z^t$ be a system of complexity at most $s$, and let $K \subset \R^d$ be an open convex cone, that is to say an open convex set which is closed under dilations. Suppose that we have the following two local solvability conditions:
\begin{itemize}
\item \textup{(Solvability at $p$)} For each prime $p$, there exists $ n \in \Z^d$ such that the forms $\psi_1( n),\ldots,\psi_t( n)$ are all coprime to $p$.
\item \textup{(Solvability at $\infty$)} There exists $ n \in K \cap \Z^d$ such that $\dot \psi_1( n), \ldots, \dot \psi_t( n) > 0$.
\end{itemize}
Then there exist infinitely many $ n \in K \cap \Z^d$ such that $\psi_1( n), \ldots, \psi_t( n)$ are all prime.
\end{corollary}

\begin{remark} This significantly generalises the main theorem in \cite{green-tao-longprimeaps} that the primes contain infinitely many progressions 
of length $k$, though for progressions of length $k > 4$ the argument here is conditional on the conjectures $\GI(k-2)$ and $\MN(k-2)$.\end{remark}

\begin{proof} If we truncate $K$ to $[-N,N]^d$, then the hypotheses ensure that $\beta_\infty \gg_{K,d} N^d$ and $\beta_p \neq 0$ for all $p$.  From Lemma \ref{sing} we conclude that $\beta_\infty \prod_p \beta_p \gg_{K,\Psi,d,t} N^d$, and the claim now follows by letting $N \to \infty$.
\end{proof}

\textsc{Acknowledgement.} The authors would like to thank the two referees, who both produced extremely careful and helpful reports which have improved the presentation of this paper.

\section{Overview of the paper}

This section is a kind of roadmap for the rest of the paper, and is somewhat informal in nature.  Also, it employs some
terminology which will only be rigorously defined in later sections.

The bulk of the paper will be concerned with the proof of the Main Theorem.  
A substantial portion of our argument consists of reprising the transference principle machinery from \cite{green-tao-longprimeaps}. This allows us to model certain unbounded functions, such as $\Lambda$, by bounded ones. Another large component of this paper consists of some facts on nilmanifolds which are essentially contained in papers in the ergodic literature, particularly that of Host and Kra \cite{host-kra}. Unfortunately, as our situation here is slightly different from that in \cite{green-tao-longprimeaps} we cannot simply cite the results we need directly from that paper, and for similar reasons we cannot cite the nilmanifold material directly.  Thus we have placed a large number of appendices in this paper in which we slightly modify the arguments from these sources to suit our present needs.  

In \S \ref{linalg} we use linear algebra to deduce Theorem \ref{main2} from the Main Theorem, and also to reduce the Main Theorem to a simplified form, Theorem \ref{main-normal}, in which the archimedean factor $\beta_\infty$ is not present and the system $\Psi$ is in a certain ``normal form''.  Then we use the ``$W$-trick'' from \cite{green-tao-longprimeaps} to eliminate the local factors $\beta_p$ and reduce matters to establishing a discorrelation estimate, Theorem \ref{main-normal-w-again}, for certain variants $\Lambda'_{b,W}-1$
of the von Mangoldt function.

In \S \ref{envelope-sec}, we recall one of the main ingredients of \cite{green-tao-longprimeaps}. This is the idea that the von Mangoldt function $\Lambda$, or more precisely the variants $\Lambda'_{b,W}-1$, are dominated by a certain \emph{enveloping sieve} $\nu$ which obeys some good pseudorandomness properties.  The verification of these properties is essentially given in \cite[Ch. 9,10]{green-tao-longprimeaps}. We take the opportunity, in Appendix \ref{gy-sec}, to give a simpler variant along the lines of unpublished notes of the second author \cite{tao-gy-notes}.

In \S \ref{gowersdef-sec} we recall the \emph{generalised von Neumann theorem} from \cite{green-tao-longprimeaps}, which allows us to use the pseudorandom enveloping sieve $\nu$ to deduce the desired discorrelation estimate, Theorem \ref{main-normal-w-again}, from a Gowers uniformity estimate on
$\Lambda'_{b,W}-1$. This latter estimate is the content of Theorem \ref{gowers-norm}. We in fact provide a more general type of generalised von Neumann theorem: the one in \cite{green-tao-longprimeaps} was specific to the case of arithmetic progressions, and did not allow one to count points inside an arbitrary convex body $K$. The basic theory of Gowers uniformity norms is reviewed in Appendix \ref{gowersnorm-sec}, whilst the generalised von Neumann theorem itself is proved in Appendix \ref{gvn-app}, following some preliminaries on convex geometry in Appendix \ref{convexgeom}.

To prove the Gowers uniformity estimate, we begin by stating in \S \ref{gowers-sec} the two conjectures we need, namely the inverse Gowers-norm conjecture $\GI(s)$ and the \mobname conjecture $\MN(s)$.  At this point we pause to present some easy consequences of these conjectures, deducing in \S \ref{mob-lio-sec} some results concerning the behaviour of the M\"obius and Liouville functions along systems of linear forms. These functions have an advantage over $\Lambda$, in that they are bounded by $1$.

In \S \ref{sec5} we apply the transference principle technology from \cite{green-tao-longprimeaps} to extend the inverse Gowers-norm conjecture $\GI(s)$ to cover functions which are bounded only by a pseudorandom measure. This result, Proposition \ref{gow-ps}, is in a sense the conceptual heart of the paper.  Once this is done the matter is reduced to the task of showing that $\Lambda'_{b,W}-1$ is asymptotically orthogonal to nilsequences. The precise statement of such a result is Proposition \ref{manortho}.

At this point we need a technical reduction, replacing a nilsequence by a slightly better behaved \emph{averaged nilsequence}.  This
reduction is carried out in \S \ref{average-sec}, and uses some basic structural facts about nilmanifolds and the cubes within them. These facts are somewhat difficult to extract from the literature, so we give them in Appendix \ref{nil-app}. In preparing this appendix we benefitted much from conversations with Sasha Leibman.

Finally, to show that $\Lambda'_{b,W}-1$ is asymptotically orthogonal to an averaged nilsequence, we split $\Lambda$ into a ``smooth'' part $\Lambda^\sharp$ and a ``rough'' part $\Lambda^\flat$. This is a fairly standard construction in analytic number theory which we learnt from \cite{iwaniec-kowalski}. The contribution of the smooth part $\Lambda^{\sharp}$ can be handled by the Gowers-Cauchy-Schwarz inequality \eqref{gcz-u}, combined with correlation estimates for truncated divisor sums. The latter type of estimates are given in Appendix \ref{gy-sec} -- the technology is that we used to build the enveloping sieve.  The rough part $\Lambda^{\flat}$ can be handled by the \mobname conjecture $\MN(s)$, thus concluding the proof.

In \S \ref{remarks-sec} we gather some concluding remarks concerning possible extensions of our results, as well as possibilities for making our estimates effective. We also indicate a proof of (say) the asymptotic in Example \ref{example1} which is somewhat shorter than the one given here, but is harder to motivate from the conceptual point of view.

In \S \ref{bounds-sec} we gather some remarks concerning bounds for the error terms in our main results. The most interesting part of this discussion focusses on what can be said assuming GRH, since unconditionally all error terms are at present completely ineffective.

The remainder of the paper consists of appendices which supply proofs for various results that we need, but which require techniques which are either standard or somewhat outside the line of the main portion of the paper.

\section{General notation}\label{notation-sec}

Our conventions for asymptotic notation are as follows.
We use $O_{a_1,\ldots,a_k}(X)$ to denote a quantity which is bounded in magnitude by $C_{a_1,\ldots,a_k} X$ for some finite positive quantity $C_{a_1,\ldots,a_k}$ depending only on $a_1,\ldots,a_k$; we also write $Y \ll_{a_1,\ldots,a_k} X$ or $X \gg_{a_1,\ldots,a_k} Y$ for 
the estimate $|Y| \leq O_{a_1,\ldots,a_k}(X)$.  

In this paper we always think of the parameter $N$ as ``large'' or ``tending to infinity''. Thus we use $o_{a_1,\ldots,a_k}(X)$ to denote a quantity bounded by $c_{a_1,\ldots,a_k}(N) X$, where $c_{a_1,\ldots,a_k}(N)$ is a quantity which goes to zero as $N \to \infty$ for each fixed $a_1,\ldots,a_k$. We do not assume that the convergence is uniform in these parameters $a_1,\ldots,a_k$.

We do not require the implied constants $C_{a_1,\ldots,a_k}$, $c_{a_1,\ldots,a_k}(N)$ to be effective.  While the arguments presented in this paper are entirely effective, the bounds that arise in the \mobname conjecture $\MN(s)$, Conjecture \ref{mnconj}, inevitably involve Siegel zeroes and are thus ineffective with current technology. They are, however, effective if the GRH is assumed.

The $o$-notation being reserved for functions which become small as $N \rightarrow \infty$, we introduce a further notation, the $\kappa$-notation, for functions which tend to zero as their parameters become \emph{small}. Thus $\kappa(\delta)$ denotes a quantity which tends to $0$ as $\delta \rightarrow 0$. Once again the $\kappa$ may be subscripted by other parameters, indicating a rate of decay which depends on those parameters.

We will frequently take advantage of the fact that two errors involving different parameters can often be concatenated by choosing one of the parameters properly.  To give a typical example, suppose we have a quantity $Q(N)$ for which we have established the bound
\begin{equation}\label{qnform}
Q(N) \leq o_{\epsilon}(1) + \kappa(\epsilon)
\end{equation}
where $\epsilon \in (0,1)$ is a parameter at our disposal and $Q(N)$ does not depend on $\epsilon$.  Then we can concatenate the two error terms by optimising in $\epsilon$ and conclude that
\begin{equation}\label{qsmall}
 Q(N) = o(1).
\end{equation}
Indeed for fixed $\epsilon$ one may choose $N$ so large that the $o_{\epsilon}(1)$ term in \eqref{qnform} is at most $\epsilon$. This means that $Q(N) = \epsilon + \kappa(\epsilon)$, still a function of the form $\kappa(\epsilon)$. Since $\epsilon$ can be as small as one likes, one obtains $Q(N) = o(1)$. Note that this kind of trick was already used to deduce Conjecture \ref{dickson-again} from Conjecture \ref{dickson}.

If $A$ is a finite non-empty set and $f: A \to \C$ is a function, we write $|A|$ for the cardinality of $A$ and
$\E_{x \in A} f(x) := \frac{1}{|A|} \sum_{x \in A} f(x)$ for the average of $f$ on $A$.  We extend this notation to functions of
several variables in the obvious manner, thus for instance $\E_{x \in A, y \in B} f(x,y) := \frac{1}{|A| |B|} \sum_{x \in A} \sum_{y \in B} f(x,y)$.

For any integer $N \geq 1$, we use $[N]$ to denote the discrete interval $[N] := \{1,\ldots,N\}$, while $\Z_N$ denotes the cyclic group $\Z_N := \Z/N\Z$.  At some places in the argument it will be convenient to pass from intervals $[N]$ to cyclic groups $\Z_N$, possibly after modifying $N$ by a constant multiplicative factor.

The letter $i$ is too important for use only as the square-root of minus one. Occasionally it will be used in this capacity and as an index in the same formula. This ought not to cause any confusion; an earlier attempt to write $\sqrt{-1}$ throughout made several of our formulae rather difficult to read.

In an earlier version of the paper we used vector notation such as $\vec{x}$ to indicate that certain elements lay in product spaces such as $\Z^d$. It was discovered that consistent use of this notation rendered certain of our expressions rather difficult to read, and so we have abandoned this practice. The reader may, at certain times, need to carefully remind herself of the spaces in which certain variables take values.

\textsc{Important convention.} For the rest of the paper, the parameters $t, d, s, L$ (which control the size and complexity of our system $\Psi = (\psi_i)_{i \in [t]}$ of linear forms). All implied constants in the $\ll$, $O(\;)$, or $o(\;)$ notation are understood to be dependent on these parameters $t,d,s,L$, even if we do not subscript them explicitly.  In particular, any quantity depending just on $t,d,s,L$ is automatically $O(1)$.  Note however that we do allow our system $\Psi$ to vary (for instance, in order to encompass Vinogradov's three-primes theorem, $\Psi$ must depend on $N$), and our estimates will be uniform in the choice of $\Psi$ so long as the parameters $t,d,s,L$ remain fixed.

\section{Linear algebra reductions}\label{linalg}

In this section we show how the Main Theorem implies Theorem \ref{main2}, and also reduce the Main Theorem to the case in which the system $\Psi$ is placed in a suitable ``normal form''. More precisely, in this section we reduce both the Main Theorem and Theorem \ref{main2} to the simpler Theorem \ref{main-normal}.
Our methods here use only elementary linear algebra.  In particular we do not require precise knowledge of exactly what the conjectures $\GI(s)$, $\MN(s)$ are at this point.  We will however restrict to the case $s \geq 1$, because the case $s=0$ follows from the $s=1$ case (note that the conjectures $\GI(1)$, $\MN(1)$ are known to be true) and in any event the $s=0$ case can be easily deduced from \eqref{siegel}.  This allows us to avoid some degeneracies later on.

\textsc{Derivation of Theorem \ref{main2} from the Main Theorem.} Suppose that we are in the situation of the Main Theorem.  Because $A$ has full rank, and ${b}$ lies in the set $A \Z^t$, the set $\Gamma := \{ {x} \in \Z^t: A {x} = {b} \}$ is a non-empty affine sublattice of $\Z^t$ of rank $d := t-s$.  Since ${b} = O(N)$ and $A$ have bounded integer coordinates, it is not hard to see that $\Gamma$ must contain at least one point of magnitude $O(N)$. For instance, one could apply any standard linear algebra algorithm to produce an element of $\Gamma$, which will then necessarily have magnitude $O(N)$ from inspection of the algorithm.  Furthermore, the generators of this lattice can also be chosen to have magnitude $O(1)$, again by applying standard linear algebra algorithms.  Thus we have a multiplicity-free parameterisation
$\Gamma = \Psi(\Z^{t-s})$  for some system of affine-linear forms $\Psi = (\psi_1,\ldots,\psi_t)$ with
$\|\Psi\|_N = O(1)$.

The full rank of $A$ ensures that the codimension of $\Psi(\Z^d)$ is the minimal value, namely $s$.  We can then write the left-hand side of \eqref{ziggy} as
$$
\sum_{ n \in K' \cap \Z^{t-s}} \prod_{i \in [t]} \Lambda( \psi_i( n) )$$
where $K' \subset \R^{t-s}$ is the convex body
$$ K' := \{  y \in \R^{t-s}: \Psi( y) \in K \}.$$ Note that $K'$ is contained in the box $[-N',N']^{t-s}$ for some $N' = O(N)$.

If two of the $\psi_i$ were affinely dependent then two of the coordinates of lattice points in $\Gamma$ would obey an affine-linear constraint. This is equivalent to the row space of $A$ containing a non-trivial vector with at most two non-zero entries, which is contrary to assumption.  From Lemma \ref{complex} we conclude that $\Psi$ has complexity at most $s$.  We now invoke the Main Theorem.  Comparing \eqref{dickson-eq} with
\eqref{ziggy} we see that we will be done as soon as we show that 
$\alpha_\infty \prod_p \alpha_p = \beta_\infty \prod_p \beta_p + o(N^d)$.
For any fixed prime $p$,
the set $\{  n \in \Z^{t-s}: \Psi( n) \in [-M,M]^t \}$ is asymptotically uniformly distributed in
residue classes in $\Z_p^{t-s}$ in the limit $M \to \infty$ and hence $\alpha_p = \beta_p$.  Since the product $\prod_p \beta_p$ is either zero or comparable to $1$, it thus suffices to show that
$\alpha_\infty = \beta_\infty + o(N^d)$.  But this follows from \eqref{beta-gauss}.  \endproof

\textsc{Elimination of the archimedean factor.} We now return to the task of proving the Main Theorem, using some simple linear algebra to obtain some reductions.  

First of all, we can use the following easy trick to hide the ``archimedean factor'' $\beta_\infty$ from view.  Clearly
we may intersect $K$ with the convex set $\Psi^{-1}( (\R^+)^t )$ and reduce to the case where $\psi_i > 0$ on $K$; in this case $\beta_\infty$ is simply the
volume of $K$.  In light of \eqref{beta-gauss} and the boundedness of the product $\prod_p \beta_p$, 
we can then rewrite \eqref{dickson-eq} as
\begin{equation}\label{dick-2}
\sum_{ n \in K \cap \Z^d} \big(\prod_{i \in [t]} \Lambda( \psi_i( n) )  - \prod_p \beta_p\big) = 
o(N^d).
\end{equation}

\begin{remark} One can easily verify the ``local'' version of this formula,
$$
\sum_{ n \in K \cap \Z^d} \big(\prod_{i \in [t]} \Lambda_{\Z_p}( \psi_i( n) ) - \beta_p\big) = 
o_{p}(N^d);$$
indeed this is a variant of the identity $\alpha_p = \beta_p$ discussed previously.
\end{remark}

It turns out to be convenient to strengthen the condition $\psi_i > 0$ slightly, say to $\psi_i > N^{9/10}$. The exact power of $N$ is not important so long as it lies between $0$ and $1$. One can easily verify, by estimating $\Lambda$ crudely by $\log N$, that for each $i$ the contribution of the case $0 \leq \psi_i( n) \leq N^{9/10}$ to \eqref{dick-2} is 
$o(N^d)$.
We have thus reduced to showing

\begin{theorem}[Finite complexity generalised H-L conjecture, again]\label{main-alt}  Let $s \geq 1$, and let $\Psi: \Z^d \to \Z^t$ be a  system of affine-linear forms of complexity $s$.  Suppose that the inverse Gowers-norm conjecture $\GI(s)$ and the \mobname conjecture $\MN(s)$ are true.  Let $N > 1$ and suppose that $\|\Psi\|_N = O(1)$.  Let $K \subset [-N,N]^t$ be a convex body such that $\psi_1,\ldots,\psi_t > N^{9/10}$ on $K$.  Then \eqref{dick-2} holds.
\end{theorem}

\textsc{Normal form reduction of the Main Theorem.} We now reduce Theorem \ref{main-alt} further by placing the system $\Psi$ in a convenient ``normal form''.  We denote
the standard basis of $\Z^d$ by $e_1,\ldots,e_d$.

\begin{definition}[Normal form]\label{indep-def}  Let $\Psi = (\psi_1,\ldots,\psi_t)$ be a system of affine-linear forms on $\Z^d$, and let $s \geq 0$.  We say that $\Psi$ is in \emph{$s$-normal form} if for every $i \in [t]$, there exists a collection $J_i \subseteq  \{e_1,\ldots,e_d\}$ of basis vectors of cardinality $|J_i| \leq s+1$ such that $\prod_{e \in J_i} \dot \psi_{i'}(e)$ is non-zero for $i' = i$ and vanishes otherwise.
\end{definition}

If a system is in $s$-normal form, then we can explicitly see that for each $i \in [t]$ the $i$-complexity of the system is at most $s$.  Indeed, we can cover the $t-1$ forms $\{ \psi_j: j \in [t] \backslash \{i\} \}$ by $|J_i|$ classes, where the class associated to a basis vector $e \in J_i$ is simply the collection of all the forms $\psi_{i'}$ for which $\dot \psi_{i'}(e) = 0$;
since $\dot \psi_i(e) \neq 0$, we see that $\psi_i$ cannot lie in the affine span of such a class. It is, therefore, necessary that a system be of a finite complexity $s$ before admitting an $s$-normal form.  We now investigate the converse relationship, beginning with some illustrative examples.

\begin{example}\label{ex1}  The system of affine-linear forms $\Psi(n_1,n_2) := (n_1, n_1 + n_2, n_1 + 2n_2, n_1 + 3n_2)$, which counts progressions of length four, has complexity $2$ but is not in $s$-normal form for any $s$. However the system of affine-linear forms
$$\Psi'(n_1,n_2,n_3,n_4) := (n_2 + 2n_3 + 3n_4, -n_1 + n_3 + 2n_4, -2n_1 - n_2 + n_4, -3n_1 - 2n_2 - n_3),$$
which also counts progressions of length four, is also of complexity $2$ and is now in $2$-normal form.
\end{example}

\begin{example}\label{ex2} The system in Example \ref{cube-ex}, which counts $(d-1)$-dimensional cubes,
has complexity $d-2$ but is not in $s$-normal form for any $s$.  However the system
$$ \Psi'( n_1,\ldots,n_{d-1}, n'_1,\ldots,n'_{d-1} ) = \big( \sum_{i \in A} n_i + \sum_{i \in [d-1] \backslash A} n'_i \big)_{A \subset [d-1]},$$
which also counts $(d-1)$-dimensional cubes, is also of complexity at most $d-2$ and is now in $(d-2)$-normal form.
\end{example}

\begin{example}\label{balogex}  Let $t := \frac{d(d+1)}{2}$, and consider the system of affine-linear forms \[ \Psi(n_1,\ldots,n_d) := (n_i+n_j+1)_{1 \leq i \leq j \leq d}\] from Example \ref{balog-ex}.  This system has complexity $1$ but is not in $s$-normal form for any $s$.  However, if we increase the number of parameters from $d$ to $2d$, and consider the system
$$ \Psi'(n_1,\ldots,n_d,n_{d+1},\ldots,n_{2d}) := \big( n_i+n_j+1 + n_{d+i}+n_{d+j} - \sum_{k=d+1}^{2d} n_k \big)_{1 \leq i \leq j \leq d},$$
which count the same type of pattern, then this system still has complexity $1$ and is now in $1$-normal form. Indeed for the off-diagonal forms $i < j$ we may use the basis vectors $e_i,e_j$, while for the diagonal forms $i=j$ we may use the basis vectors $e_i, e_{d+i}$.    
\end{example}

\begin{remark} Informally speaking, if $(\psi_1,\ldots,\psi_t)$ is in $s$-normal form, then for each form $\psi_i$ there exist a set of at most $s+1$ variables $(n_j)_{j \in J_i}$, such that $\psi_i$ is the only form which truly utilises all the variables at once.  As we shall see later, this property will be convenient for establishing a ``generalised von Neumann theorem'' (Proposition \ref{gvn}), which roughly speaking controls averages such as \eqref{dick-2} in terms of Gowers uniformity norms, which we shall recall in Appendix \ref{gowersnorm-sec}. \end{remark}

Now we investigate the converse question, namely whether every system of complexity $s$ has a normal form representation.  
To formalise this we first need the concept of \emph{extending} a system of affine-linear forms by adding some ``dummy'' parameters:

\begin{definition}[Extensions]  Let $\Psi: \Z^d \to \Z^t$ be a system of affine-linear forms.  An \emph{extension} of this system is a system $\Psi': \Z^{d'} \to \Z^t$ with $d' \geq d$, such that
\begin{equation}\label{lattice-eq}
\Psi'(\Z^{d'}) = \Psi(\Z^d)
\end{equation}
and furthermore if we identify $\Z^d$ with the subset $\Z^d \times \{0\}^{d'-d}$ of $\Z^{d'}$ in the obvious manner, then $\Psi$ is the restriction of $\Psi'$ to $\Z^d$.
\end{definition}

We note that if $\Psi$ is in $s$-normal form at $i$, and if $\Psi'$ is an extension of $\Psi$, then $\Psi'$ is also in $s$-normal form at $i$. By the same token, we note also that if $\Psi = (\psi_i)_{i = 1}^d$ is in $s$-normal form, then so is any subsystem $(\psi_i)_{i \in I}$, $I \subset \{1,\dots,d\}$.

\begin{example} In Example \ref{balog-ex}/Example \ref{balogex}, $\Psi'$ is an extension of $\Psi$.  This is not quite the case in Examples \ref{ex1}, \ref{ex2}, because $\Psi$ is not a restriction of $\Psi'$.  However in these two examples, the direct sum $\Psi \oplus \Psi'$ of the two systems is both an extension of $\Psi$ \emph{and} in normal form; for instance, in Example \ref{ex1} the system
\begin{align*}
\Psi \oplus \Psi'(n_1,n_2,n'_1,n'_2,n'_3,n'_4) &:= (n_1 + n'_2 + 2n'_3 + 3n'_4, n_1+n_2-n'_1 + n'_3 + 2n'_4, \\
&\quad n_1+2n_2-2n'_1 - n'_2 + n'_4, n_1+3n_2-3n'_1 - 2n'_2 - n'_3)
\end{align*}
is an extension of $\Psi$ which is in $2$-normal form.
\end{example}

\begin{lemma}[Existence of normal forms]\label{normal-lemma} Let $\Psi: \Z^d \to \Z^t$ be a system of affine-linear forms of some finite complexity $s$.  Then there exists an extension $\Psi': \Z^{d'} \to \Z^t$ of $\Psi$ which is in $s$-normal form, where
$d' = O(1)$. Furthermore if the original system $\Psi$ had size $\|\Psi\|_N = O(1)$, then the same is true of the extended system $\Psi'$.
\end{lemma}

\begin{proof} Let us fix $i \in [t]$.  We shall obtain an extension $\Psi': \Z^{d'} \to \Z^t$ of $\Psi$ which in $s$-normal form at $i$, by which we mean that there is a collection $J_i \subseteq \{e_1,\ldots,e_{d'}\}$ of basis vectors of cardinality $|J_i| \leq s+1$ such that $\prod_{e \in J_i} \dot \psi'_{i'}(e)$ is non-zero for $i' = i$ and vanishes otherwise.  Applying this extension procedure once for each value of $i$ we shall obtain the result.

By hypothesis, $\Psi$ has $i$-complexity at most $s$, and so we can cover $[t] \backslash \{i\}$ by $s+1$ classes $A_1,\ldots,A_{s+1}$, such that $\psi_i$ is not in the affine-linear span of $\{ \psi_j: j \in A_k \}$ for $k \in [s+1]$.
In particular, this implies that one can find vectors $f_1,\ldots,f_{s+1} \in \Q^d$ which ``witness this fact'', that is to say such that $\dot \psi_j(f_k) = 0$ and $\dot \psi_i(f_k) \neq 0$
all $k \in [s+1]$ and $j \in A_k$.  By clearing denominators we can take $f_1,\ldots,f_{s+1} \in \Z^d$.  Since
$\dot \Psi$ has bounded integer coefficients we also see that $f_1,\ldots,f_{s+1} = O(1)$.
If we now let $d' := d+s+1$ and let $\Psi': \Z^{d'} \to \Z^t$ be the system
$$ \Psi'(  n, m_1, \ldots, m_{s+1} ) := \Psi(  n + m_1 f_1 + \ldots + m_{s+1} f_{s+1} )$$
for all $ n \in \Z^d$ and $m_1,\ldots,m_{s+1} \in \Z$, we easily verify that $\Psi'$ satisfies the desired $s$-normal form property at $i$, as well as the size bounds on $\Psi'$. By repeating this procedure once for each $i$ we obtain the claim.
\end{proof}

Using this lemma it is not hard to show that, in order to prove the Main Theorem, it suffices to prove the following result for $s$-independent systems.  

\begin{theorem}[Primes in affine lattices in normal form]\label{main-normal}  Let $s \geq 1$, and let $\Psi: \Z^d \to \Z^t$ be a  system of affine-linear forms of complexity $s$ in $s$-normal form.  Suppose that the inverse Gowers-norm conjecture $\GI(s)$ and the \mobname conjecture $\MN(s)$ are true.  Let $N > 1$ and suppose that $\|\Psi\|_N = O(1)$.  Let $K \subseteq [-N,N]^t$ be a convex body such that $\psi_1,\ldots,\psi_t > N^{8/10}$ on $K$. Then \eqref{dick-2} holds, that is to say \[ \sum_{ n \in K \cap \Z^d} \big(\prod_{i \in [t]} \Lambda( \psi_i( n) )  - \prod_p \beta_p\big) = 
o(N^d). \] 
\end{theorem}

\begin{proof}[Proof of the Main Theorem assuming Theorem \ref{main-normal}]
By our earlier reduction it suffices to show that Theorem \ref{main-alt} holds.
Let $\Psi$, $K$, $N$ be as in
Theorem \ref{main-alt}.  We may assume $N$ large as the claim is trivial for $N$ small.

Let $\Psi': \Z^{d'} \to \Z^t$ be the $s$-normal form extension given by Lemma \ref{normal-lemma}.  
An inspection of the proof of that lemma allows us to find vectors $f_{d+1},\ldots,f_{d'} \in \Z^d$ of magnitude $O(1)$ such that
$$ \Psi'(  n, m_{d+1}, \ldots, m_{d'} ) := \Psi(  n + m_{d+1} f_{d+1} + \ldots + m_{d'} f_{d'} ).$$
(One can also deduce the existence of these vectors directly from the conclusions of Lemma \ref{normal-lemma}.)  We observe
that the local factors $\beta'_p$ associated to the system $\Psi'$ are precisely the same as the local factors $\beta_p$ associated to $\Psi$; this is ultimately due to the translation-invariance of $\Z_p$.  Now let $K' \subseteq \R^{d'}$ be the convex body
$$ K' := \{ (  n, m_{d+1}, \ldots, m_{d'} ) \in \R^d \times [-N,N]^{d'-d}:  n + m_{d+1} f_{d+1} + \ldots + m_{d'} f_{d'} \in K \}.$$
This is contained in $[-N',N']^{d'}$ for some $N' = O(N)$.  Applying Theorem \ref{main-normal} we conclude
$$
\sum_{( n, m) \in K' \cap \Z^{d'}} \big(\prod_{i \in [t]} \Lambda( \psi'_i( n, m) ) - \prod_p \beta_p \big) = o(N^{d'}).$$
Making the change of variables $ r :=  n + m_{d+1} f_{d+1} + \ldots + m_{d'} f_{d'}$, the left-hand side can be simplified to 
$$ \left|[-N,N]^{d'-d} \cap \Z^{d'-d}\right| \sum_{ r \in K \cap \Z^d} \big(\prod_{i \in [t]} \Lambda( \psi_i( r) ) - \prod_p \beta_p\big)$$
and \eqref{dick-2} follows upon dividing out by $(2N+1)^{d'-d}$.
\end{proof}

This completes our linear algebra manipulations. It now remains to prove Theorem \ref{main-normal}, a task which will occupy the remainder of the paper.

\section{The $W$-trick}

In the preceding section we were able to eliminate the archimedean factor $\beta_\infty$ by assuming that $\psi_1,\ldots,\psi_t$
were non-negative on $K$, and using the formulation \eqref{dick-2}.  Now we use a somewhat similar trick, which we term the ``$W$-trick''. This was a vital trick in \cite{green-roth-primes,green-tao-longprimeaps,green-tao-selberg}, where it was used in similar fashion to eliminate the local factors $\beta_p$.
Once again, the reductions here will not actually require any knowledge of the two conjectures $\GI(s)$ and $\MN(s)$,
which we shall finally introduce in \S \ref{gowers-sec}.

\textsc{Important convention.} From now on in the paper, fix some slowly growing function $w = w(N)$. Any function such that $w(N) \leq \frac{1}{2}\log\log N$ and $\lim_{N \to \infty} w(N) = \infty$ would suffice; for sake of definiteness we shall conservatively set $w := \log\log\log N$. The exact choice of $w$ is only relevant for determining the decay rate of the $o()$ terms, but as our final decay bounds are ineffective we will not attempt to optimise in $w$.

We define the quantity $W = W(w)$ by
$$ W := \prod_{p \leq w} p;$$
since $w \leq \frac{1}{2} \log \log N$ we have $W = O( \log^{1/2} N )$.
For each $b \in [W]$ with $\gcd(b,W)=1$, let $\Lambda_{b,W}: \Z^+ \to \R^+$ be the function
\begin{equation}\label{lbn}
 \Lambda_{b,W}(n) := \frac{\phi(W)}{W} \Lambda(Wn+b)
\end{equation}
where we recall that $\phi(W) = \#\{ b \in [W]: \gcd(b,W) = 1 \}$ is the Euler totient function of $W$.  Thus for instance the prime number theorem in APs \eqref{siegel} asserts\footnote{In order to obtain this statement for $w$ as large as $\frac{1}{2} \log \log N$, one needs a more quantitative version of \eqref{siegel} such as the Siegel-Walfisz theorem.} that $\Lambda_{b,W}(n)$ has average value $1$ as $n \to \infty$.  Actually it will be slightly more convenient to work with the variant
$$ \Lambda'_{b,W}(n) := \frac{\phi(W)}{W} \Lambda'(Wn+b)$$
where $\Lambda'$ is the restriction of $\Lambda$ to the primes, i.e. $\Lambda'(p) = \log p$ for all primes $p$ and $\Lambda'(n) = 0$ for non-prime $p$.   Thus $\Lambda'$ only differs from $\Lambda$ on the (negligible) set of prime powers $p^2, p^3, \ldots$. 

Recall that we reduced the task of proving the Main Theorem to that of proving Theorem \ref{main-normal}. We now make a further reduction, showing that it suffices to prove the following.

\begin{theorem}[W-tricked primes in affine lattices]\label{main-normal-w}  Let $s \geq 1$, and suppose that $\Psi = (\psi_1,\ldots,\psi_t): \Z^d \to \Z^t$ is a system of affine-linear forms in $s$-normal form and with $\Vert \Psi \Vert_N = O(1)$.  Suppose that the inverse Gowers-norm conjecture $\GI(s)$ and the \mobname conjecture $\MN(s)$ are true.  Let $K \subseteq [-N,N]^t$ be any convex body on which $\psi_1,\ldots,\psi_t > N^{7/10}$.  Then for any $b_1,\ldots,b_t \in [W]$ which are coprime to $W$, we have
$$
\sum_{ n \in K \cap \Z^d} \big(  \prod_{i \in [t]} \Lambda'_{b_i,W}( \psi_i( n) )  - 1\big) = 
o(N^d).
$$
\end{theorem}

\begin{remark} Note that the bounds on the right do not depend on $b_1,\ldots,b_t$.  The philosophy here is that the functions $\Lambda'_{b,W}$ should behave ``pseudorandomly'' with average value one; this is in contrast with $\Lambda$, which has many local irregularities with respect to small moduli which necessitate the introduction of the local factors $\beta_p$.  This philosophy of passing from $\Lambda$ to the more uniformly distributed $\Lambda'_{b,W}$ underlies the arguments in \cite{green-tao-longprimeaps}.  In \S \ref{wtrick-sec} we will have to invert the $W$-trick and deduce some correlation estimates on $\Lambda'_{b,W}$ from that on $\Lambda$.
\end{remark}

\begin{proof}[Proof of the Main Theorem assuming Theorem \ref{main-normal-w}] By previous reductions, it suffices to establish Theorem \ref{main-normal}. Let $\Psi, K$ be as in Theorem \ref{main-normal}.  We may then replace $\Lambda$ by $\Lambda'$ as the contribution of the prime powers is easily seen to be negligible.
To prove \eqref{dick-2}, it then suffices by \eqref{beta-gauss} to show that
\begin{equation}\label{beat}
\sum_{ n \in K \cap \Z^d} \prod_{i \in [t]} \Lambda'( \psi_i( n) ) = \vol_d(K) \prod_p \beta_p + o(N^d).
\end{equation}
We may take $N$ to be large, since the claim is trivial otherwise.  

Now the upper bound on $w$ ensures that $W \leq \log N$.  From Lemma \ref{sing} followed by the multiplicativity of the local factors $\beta$ 
we have 
$$\prod_p \beta_p = \prod_{p \leq w} \beta_p + o(1) = \beta_W + o(1);$$
since $\vol_d(K) = O(N^d)$, we conclude that
$$ \vol_d(K) \prod_p \beta_p = \vol_d(K) \beta_W + o(N^d).$$
Now let $A$ be the set
$$A := \{  a \in [W]^d: \gcd(\psi_i( a),W)=1 \hbox{ for all } i \in [t] \}.$$
Then from \eqref{beta-p} we have $\beta_W = \left(\frac{W}{\phi(W)}\right)^t |A| / W^d$, which implies that
\begin{equation}\label{volbeta}
\vol_d(K) \prod_p \beta_p = \sum_{ a \in A} \bigg(\frac{W}{\phi(W)}\bigg)^t W^{-d} \vol_d(K) + o(N^d).
\end{equation}
Also, from Lemma \ref{sing} we know that $\beta_W$ is comparable to $1$, and so
\begin{equation}\label{a-card}
|A| \ll \left(\frac{\phi(W)}{W}\right)^t W^d.
\end{equation}
Next, note that by a simple expansion we have
\begin{equation}\label{vuk}
 \sum_{ n \in K \cap \Z^d} \prod_{i \in [t]} \Lambda'( \psi_i( n) )
= \sum_{ a \in [W]^d} 
\sum_{\substack{ n \in \Z^d \\ W  n +  a \in K}} \prod_{i \in [t]} \Lambda'( \psi_i(W  n +  a) ).
\end{equation}
If $ a$ does not lie in $A$, then $\psi_i(W  n +  a)$ will not be coprime to $W$ for some $i \in [t]$.
Since $\psi(W  n +  a) > N^{7/10}$ by hypothesis, and $W$ is so small compared to $N$, we see that $\Lambda'(\psi_i(W  n +  a)) = 0$.  Thus we may restrict $ a$ to $A$.  Now for each $ a \in A$ and $i \in [t]$, we can write
$$\psi_i( W  n +  a ) = W \tilde \psi_{i, a}( n) + b_i( a)$$
where $b_i( a)$ lies in $[W]$ and is coprime to $W$, while $\tilde \psi_{i, a}$ is a translate of $\psi_i$ whose constant term $\tilde \psi_{i, a}(0)$ is $O(N/W)$.  Indeed $b_i( a)$ is simply the remainder formed when dividing $\psi_i( a)$ by $W$.  We then have
$$ \Lambda'( \psi_i( W  n +  a ) ) = \frac{W}{\phi(W)} \Lambda'_{b_i( a), W}( \tilde \psi_{i,a}(  n ) ).$$
It follows from \eqref{vuk} that
\begin{equation}\label{3.5a} \sum_{{n} \in K \cap \Z^d} \prod_{i \in [t]} \Lambda'(\psi_i(n)) = \sum_{ a \in A} \left(\frac{W}{\phi(W)}\right)^t \sum_{\substack{ n \in \Z^d\\ W  n + a \in K}} 
\prod_{i \in [t]} \Lambda'_{b_i( a),W}( \tilde \psi_{i,a}(n) ).\end{equation}
However from Theorem \ref{main-normal-w} (with $N$ replaced by $\tilde N = O(N/W)$ and $\tilde K := (K - {a})/W$: note that $\Vert \tilde \Psi \Vert_{\tilde N} = O(1)$) we have
$$ \sum_{\substack{ n \in \Z^d\\ W  n + a \in K}} \big(\prod_{i \in [t]} \Lambda'_{b_i(a),W}( \tilde \psi_{i,a}(n) ) - 1\big) =  o\left(\frac{N}{W}\right)^d.$$
Recalling \eqref{a-card}, this together with \eqref{3.5a} implies that
\begin{equation}\label{eq3.5b} \sum_{{n} \in K \cap \Z^d} \prod_{i \in [t]} \Lambda'(\psi_i ({n})) = \sum_{a \in A} \left(\frac{W}{\phi(W)}\right)^t \sum_{\substack{ n \in \Z^d \\ W  n + a \in K}} \!\!\!1 \;\;\;
+  o(N^d).\end{equation}
On the other hand a simple volume-packing argument (cf. Appendix \ref{convexgeom}) yields
$$ \sum_{\substack{ n \in \Z^d\\ W  n + a \in K}} \!\!\!1 = W^{-d} \vol_d(K) + o\left(\frac{N}{W}\right)^d$$
and so, using \eqref{a-card} once more together with \eqref{eq3.5b}, we see that 
$$\sum_{{n} \in K \cap \Z^d} \prod_{i \in [t]} \Lambda'(\psi_i ({n})) = \sum_{a \in A} \left(\frac{W}{\phi(W)}\right)^t W^{-d} \vol_d(K)
+ o(N^d).$$
Subtracting this against \eqref{volbeta} we see that the left-hand side of \eqref{beat} is  $o(N^d)$.  This proves the claim.
\end{proof}

Theorem \ref{main-normal-w}, as we have just seen, implies the Main Theorem. Before moving on to the more substantial arguments in this paper, we give one further simple reduction, deducing Theorem \ref{main-normal-w} from the following variant.

\begin{theorem}[Final technical reduction]\label{main-normal-w-again}  Let $s \geq 1$, and let $\Psi = (\psi_1,\ldots,\psi_t): \Z^d \to \Z^t$ be a system of affine-linear forms in $s$-normal form.  Suppose that the inverse Gowers-norm conjecture $\GI(s)$ and the \mobname conjecture $\MN(s)$ are true.  Let $K \subseteq [-N,N]^t$ be any convex body on which $\psi_1,\ldots,\psi_t > N^{7/10}$.  Then for any $b_1,\ldots,b_t \in [W]$ which are coprime to $W$, we have
$$
\sum_{ n \in K \cap \Z^d} \prod_{i \in [t]} (\Lambda'_{b_i,W}( \psi_i( n) ) - 1)  = 
o(N^d).
$$
\end{theorem}

Indeed, Theorem \ref{main-normal-w} follows immediately from Theorem \ref{main-normal-w-again} by splitting each
$\Lambda'_{b_i,W}$ as $(\Lambda'_{b_i,W}-1) + 1$, expanding out the product in Theorem \ref{main-normal-w}, and using
Theorem \ref{main-normal-w-again} repeatedly, noting that any subsystem of $\Psi$ will still be in $s$-normal form.\endproof

The remainder of the paper shall be devoted to establishing Theorem \ref{main-normal-w-again}.

\section{The enveloping sieve}\label{envelope-sec}

In previous sections we have reduced matters to establishing a certain discorrelation estimate, Theorem \ref{main-normal-w-again}, for the functions
$\Lambda'_{b_i,W}-1$.  A major difficulty in the analysis here is that these functions are not bounded uniformly in $N$.  However, as
in \cite{green-tao-longprimeaps,green-tao-selberg}, we shall be able to import tools from sieve theory. In particular, we use the principle of the ``enveloping sieve''. This is a well-behaved function $\nu$, some constant multiple of which provides a pointwise bound for the functions $\Lambda'_{b_i,W}-1$. Of course, the function $\nu$ will not be bounded as $N \rightarrow \infty$; however it does obey a number of very good correlation or \emph{pseudorandomness} estimates which assert, roughly speaking, that $\nu$ ``effectively behaves like'' the bounded function $1$. 

To define the notion of pseudorandomness properly we recall the \emph{linear forms condition} and \emph{correlation condition} from \cite{green-tao-longprimeaps}, modified slightly for the application at hand.
In the following three definitions we assume that $N$ is a large positive integer, and that
$N' = N'(N)$ is a prime number of size $N < N' \leq O_{s,t,d,L}(N)$.

\begin{definition}[Measures]  A \emph{measure} on $\Z_{N'}$ is a function $\nu : \Z_{N'} \rightarrow \R^+$ (depending of course on $N'$ and hence on $N$) with 
\begin{equation}\label{norm-condition} \E_{n \in \Z_{N'}} \nu(n) = 1 + o(1).
\end{equation}
\end{definition}

\begin{definition}[Linear forms condition]\label{linear-forms-condition}
Let $\nu$ be a measure on $\Z_{N'}$, and let $m_0,d_0$ and $L_0$ be positive integer parameters. Then we say that $\nu$ satisfies the $(m_0,d_0,L_0)$-linear forms condition if the following holds: given $1 \leq d \leq d_0$, $1 \leq t \leq m_0$, and any finite complexity system $\Psi = (\psi_1,\ldots,\psi_t)$ of affine-linear forms on $\Z^d$ with all coefficients of $\dot \Psi$ bounded in magnitude by $L_0$, we have
\begin{equation}\label{lfc}
 \E_{ n \in \Z_{N'}^d} \prod_{i \in [t]} \nu( \psi_i(  n) ) = 1 
 + o_{m_0,d_0,L_0}(1).
\end{equation}
In this expression we induce the affine-linear forms $\psi_j: \Z_{N'}^d \to \Z_{N'}$ from their global counterparts $\psi_j: \Z^d \to \Z$ in the obvious manner. 
\end{definition}

\begin{remarks} Note that \eqref{lfc} includes \eqref{norm-condition} as a special case. Strictly speaking, it would be more accurate to call measures ``probability densities'', and the linear forms condition is really an ``affine-linear forms condition'', but we will keep the notation as above for brevity and compatibility with \cite{green-tao-longprimeaps}.  In \cite{green-tao-longprimeaps} the coefficients of the affine-linear forms were allowed to be rational with bounded numerator and denominator. Since $N'$ is a large prime, it is always possible in practice to clear denominators and deal only with forms having integer coefficients.  Note that Theorem \ref{main-normal-w} is a (conditional) assertion that the $\Lambda_{b,W}$ essentially obey the linear-forms condition.  Thus trying to establishing the linear forms condition for $\Lambda_{b,W}$ would essentially be as hard as trying to prove the Main Theorem. The point of the definition, however, is that it will suffice to achieve the much simpler task of \emph{majorising} $\Lambda_{b,W}$ by constant multiples of measures $\nu$ which obey this condition.  Finally, we note that the error term in \eqref{lfc} is uniform over all choices of constant term $\Psi(0)$.
\end{remarks}

\begin{definition}[Correlation condition]\label{correlation-condition}
Let $\nu : \Z_{N'} \rightarrow \mathbb{R}^{+}$ be a measure, and let $m_0$ be a positive integer parameter. We say that $\nu$ satisfies the \emph{$m_0$-correlation condition} if for every $1 < m \leq m_0$
there exists a weight 
function $\tau = \tau_{m}: \Z_{N'} \to \R^+$ which obeys the moment conditions
\begin{equation}\label{eq3.1}
 \E_{n \in \Z_{N'}}\tau^q(n) \ll_{m,q} 1
\end{equation}
for all $1 \leq q < \infty$
and such that
\begin{equation}\label{eq3.2}
\E_{n \in \Z_{N'}} \prod_{i \in [m]} \nu(n+h_i)
\leq \sum_{1 \leq i < j \leq m} \tau(h_i-h_j)
\end{equation}
for all $h_1, \ldots, h_m \in \Z_{N'}$, not necessarily distinct.
\end{definition}

\begin{remarks}  Because we are only seeking upper bounds here rather than asymptotics, this condition would follow from a standard upper bound sieve such as Selberg's sieve.  One should compare this condition with the much more difficult prime tuples conjecture, which is part of the ``infinite complexity'' case $d=1$, $t > 1$ of the generalised Hardy-Littlewood conjecture.
The correlation condition will only be used implicitly in this paper, as it is needed in the proof of
\cite[Proposition 8.1]{green-tao-longprimeaps}, which is in turn used in the proof of Proposition \ref{kvn-decomp}.
\end{remarks}

Let $D$ be a positive integer. We call a measure \emph{$D$-pseudorandom} if it obeys the $(D,D,D)$-linear forms and $D$-correlation conditions.
In practice, we shall work with measures which are $D$-pseudorandom where $D$ is a sufficiently large function of $s,d,t,L$. The exact value will not be terribly important for our arguments and, whilst it could be specified explicitly, we shall not do so.

Our next task is to show that the functions $\Lambda'_{b_1,W},\ldots,\Lambda'_{b_t,W}$ can be dominated by a $D$-pseudorandom measure for any fixed $D$ that we choose, providing we are willing to concede multiplicative constants that depend on $D$. 

\begin{proposition}[Domination by a pseudorandom measure]\label{pseudodom}  
Let $D > 1$ be arbitrary. Then there is a constant $C_0 := C_0(D)$ such that the following is true. Let $C \geq C_0$, and suppose that $N' \in [CN,2CN]$.  Let
$b_1,\ldots,b_t \in \{0,1,\ldots,W-1\}$ be coprime to $W := \prod_{p \leq w} p$.  Then there exists a $D$-pseudorandom 
measure $\nu: \Z_{N'} \to \R^+$ which obeys the pointwise bounds
$$ 1 + \Lambda'_{b_1,W}(n) + \ldots + \Lambda'_{b_t,W}(n) \ll_{D,C} \nu(n) $$
for all $n \in [N^{3/5}, N]$, where we identify $n$ with an element of $\Z_{N'}$ in the obvious manner.
\end{proposition}

The proof of this proposition is a minor variant of that in \cite{green-tao-longprimeaps}.  For the sake of completeness we present
a proof in Appendix \ref{gy-sec}.  The constant $C$ is a technicality needed to avoid certain ``wraparound'' issues when passing from $[N]$ to $\Z_{N'}$ and can be largely ignored.

The philosophy of the \emph{transference principle} developed in \cite{green-tao-longprimeaps} is that functions which are dominated by pseudorandom measures behave almost as if they were bounded, for the purposes of computing correlations and other multilinear averages. We shall see examples of this in later sections.  For now, we turn to the first significant step in the paper, namely
the reduction of matters to establishing a Gowers uniformity norm estimate for $\Lambda'_{b,W}-1$.

\section{Reduction to a Gowers norm estimate}\label{gowersdef-sec}

We shall informally refer to a function $f: [N] \to \C$ as being \emph{Gowers uniform of order $s$} if its Gowers uniformity norm
$\|f\|_{U^{s+1}[N]}$ is small; see Appendix \ref{gowersnorm-sec} for definitions and basic properties of this norm.  A basic principle is that Gowers uniform functions of order $s$ have a negligible impact on multilinear averages of complexity $s$ or less. An example of this is \cite[Proposition 5.3]{green-tao-longprimeaps}, but we will prove a much more general result of this type here. We refer to such statements as \emph{generalised von Neumann theorems}. The name originally came from results in ergodic theory such as \cite[Theorem 11.1]{host-kra}, but it has been convenient to use the name to describe a large number of contexts in additive combinatorics in which some kind of expression is bounded using Gowers norms\footnote{Another example of this is the \emph{Koopman von Neumann theorem,} which we will introduce  in \S \ref{sec5}.}.

A crucial observation in \cite{green-tao-longprimeaps} is that this type of principle also applies to \emph{unbounded} functions, so long as these unbounded functions are in turn dominated pointwise by a suitably pseudorandom measure. 

\begin{proposition}[Generalised von Neumann theorem]\label{gvn} Let $s,t,d,L$ be positive integer parameters as usual. Then there are constants $C_1$ and $D$, depending on $s,t,d$ and $L$, such that the following is true.  Let $C_1 \leq C \leq O_{s,t,d,L}(1)$ be arbitrary and suppose that $N' \in [CN,2CN]$ is a prime. Let $\nu:\Z_{N'} \to \R^+$ be a $D$-pseudorandom measure, and suppose that $f_1,\dots, f_t : [N] \rightarrow \R$ are functions with $|f_i(x)| \leq \nu(x)$ for all $i \in [t]$ and $x \in [N]$. Suppose that $\Psi = (\psi_1,\ldots,\psi_t)$ is a system of affine-linear forms in $s$-normal form with $\|\Psi\|_N \leq L$.  Let $K \subseteq  [-N,N]^d$ be a convex body such that $\Psi(K) \subseteq  [N]^t$.  Suppose also that
$$ \min_{1 \leq j \leq t}  \Vert f_j \Vert_{U^{s+1}[N]} \leq \delta$$
for some $\delta > 0$. Then we have
\begin{equation}\label{lin-w} 
\sum_{{n} \in K} \prod_{i \in [t]} f_i(\psi_i({n})) =  \; 
o_{\delta,C}(N^d) + \kappa_C(\delta)N^{d}.\end{equation}
\end{proposition}

\begin{remarks}
For an explanation of the $\kappa$-notation, we refer the reader to \S \ref{notation-sec}.
One could specify explicit values for $C_1,D$, but we have not done so. In applications to the primes we will always take $C \geq C_0(D)$, where $C_0$ is the function defined in Proposition \ref{pseudodom}.
\end{remarks}

This proposition is a variant of \cite[Proposition 5.3]{green-tao-longprimeaps}. It is somewhat more elaborate than that result in that it applies to a general system of affine linear forms, and one has the flexibility of summing over an arbitrary convex body. Once the convex body is handled by standard techniques, however, the only real tool that is needed is several applications of the Cauchy-Schwarz inequality. This is a common feature of generalised von Neumann theorems. We give a proof of Proposition \ref{gvn} in Appendix \ref{gvn-app}, which uses some preliminaries in Appendices \ref{convexgeom}, \ref{gowersnorm-sec} but is otherwise self-contained.  Using Propositions \ref{pseudodom} and \ref{gvn} we reduce Theorem \ref{main-normal-w-again}, and hence the Main Theorem, to the following Gowers uniformity estimate.

\begin{theorem}[Gowers uniformity estimate]\label{gowers-norm} 
Let $N, w > 1$, and let $b \in [W]$ be coprime to $W = \prod_{p \leq w} p$.  Suppose that the inverse Gowers-norm conjecture $\GI(s)$ and the \mobname conjecture $\MN(s)$ are true for some $s \geq 1$.  Then we have
$$ \| \Lambda'_{b,W} - 1 \|_{U^{s+1}[N]} = o(1).$$
\end{theorem}

\begin{remark} Observe (cf. Examples \ref{cube-ex} and \ref{ex2}) that this theorem is a special case of Theorem \ref{main-normal-w-again}.  Thus the generalised von Neumann theorem, Proposition \ref{gvn}, can be viewed as an assertion that the $U^{s+1}$ average is ``universal'' or ``characteristic'' among all multilinear averages of complexity $s$, even when dealing with functions that are bounded only by a pseudorandom measure.  
\end{remark}

\begin{proof}[Proof of Main Theorem assuming Theorem \ref{gowers-norm}] By previous reductions, is suffices to prove Theorem \ref{main-normal-w-again}. Let the notation and assumptions
be as in that theorem.  By enlarging $N$ by a multiplicative factor of $O(1)$ if necessary we may assume that $\Psi(K) \subseteq [N]^t$.
Let $D = D_{s,t,d,L}$ be the constant in Proposition \ref{gvn}, and set $C := \max(C_1,C_0(D))$, where $C_0$ is the function appearing in Proposition \ref{pseudodom} and $C_1$ is the one appearing in Proposition \ref{gvn}.
Applying Bertrand's postulate, we may select a prime $N'$ such that $CN \leq N' \leq 2CN$. Let $\nu$ be the $D$-pseudorandom measure given by \eqref{pseudodom}. Then the functions
$f_i(n) := c \cdot (\Lambda'_{b_i,W} - 1)$ will be pointwise dominated in magnitude by $\nu$ for some suitably small constant $c = c_{s,t,d,L} > 0$.  Applying Theorem \ref{gowers-norm} and Proposition \ref{gvn}, we obtain the desired estimate after dividing out the factors of $c$.
\end{proof}

We have now completed yet another reduction, and it remains to prove Theorem \ref{gowers-norm}.   Note that we have eliminated the system $\Psi$ of affine-linear forms, as well as the
convex body $K$, replacing them both with the Gowers norm $U^{s+1}[N]$; the parameters $d,t$ have also disappeared.  In order to proceed further, we need to exploit some deeper facts and conjectures concerning the Gowers norm.  In particular we shall shortly need the \emph{inverse Gowers-norm conjecture} $\GI(s)$, to which we now turn.

\section{The inverse Gowers-norm and \mobname conjectures}\label{gowers-sec}

\textsc{Nilsequences.} The purpose of this section is to state the two conjectures $\GI(s)$ and $\MN(s)$ which have appeared in many of the above theorems, most recently in Theorem \ref{gowers-norm}.  Both conjectures revolve around the concept of a \emph{nilsequence}, which we now pause to recall.

\begin{definition}[Nilmanifolds and nilsequences]\label{nil-def}
Let $G$ be a connected, simply connected, Lie group. We define the \emph{central series} $G_0 \supseteq G_1 \supseteq G_2 \supseteq \dots$ by defining $G_0 = G_1 = G$, and $G_{i+1} = [G, G_i]$ for $i \geq 2$, where the commutator group $[G,G_i]$ is the group generated by $\{ g h g^{-1}h^{-1}: g \in G, h \in G_i \}$. We say that $G$ is \emph{$s$-step nilpotent} if $G_{s+1} = {1}$. Let $\Gamma \subseteq  G$ be a discrete, cocompact subgroup. Then the quotient $G/\Gamma$ is called an $s$-step nilmanifold. If $g \in G$ then $g$ acts on $G/\Gamma$ by left multiplication, $x \mapsto g\cdot x$. By a \emph{an $s$-step nilsequence}, we mean a sequence of the form $(F(g^n  x))_{n \in \N}$, where $x \in G/\Gamma$ is a point and $F : G/\Gamma \rightarrow \R$ is a continuous function.  We say that the nilsequence is \emph{1-bounded} if $F$ takes values in $[-1,1]$.
\end{definition}

\begin{remark} For a full technical treatment of nilsequences, see \cite{corwin-greenleaf}. The reader might consult \cite{bergelson-host-kra,host-kra,kra-icm-lecture} for the ergodic theory perspective, or other papers of the authors \cite{green-icm,green-tao-u3inverse,green-tao-u3mobius} for various discussions more-or-less in the spirit of additive combinatorics.

As remarked above, the exact definition of a nilsequence will not be terribly important to our arguments here. In the $s=2$ case, representative examples of nilsequences are those associated to the \emph{Heisenberg nilmanifold}, which is discussed in detail in \cite{bergelson-host-kra,green-icm,green-tao-u3inverse,green-tao-u3mobius}.  See also the proof of Proposition \ref{gi2-prop}.
\end{remark}

\begin{remark} Note that we are requiring our nilpotent groups to be connected and simply connected.  The latter hypothesis is not overly restrictive, since if $G$ \emph{is} connected, then it may be assumed to be simply connected by passing to a universal cover. 
The connectedness assumption however is more substantial; the nilpotent groups constructed in the ergodic theory literature (e.g. in \cite{host-kra}) are not always shown to be connected.  However, Sasha Leibman \cite{sasha-personal} has indicated to us that it suffices, in the context of the $\GI(s)$ conjecture, to deal with connected $G$. We will elaborate on this point in a future paper if necessary, but the issue does not need to be addressed here. This is because the arguments used in proving the cases $s \leq 2$, which are the only cases of the conjectures established so far, give connectedness as a byproduct.
\end{remark}

As we shall need to be rather quantitative regarding these nilmanifolds, we shall arbitrarily endow\footnote{Strictly speaking, we are abusing notation here; a nilmanifold should not be represented solely by the quotient space $G/\Gamma$, but rather as a quadruplet $(G, \Gamma, G/\Gamma, d_{G/\Gamma})$ (and the Lie group $G$ should in turn be expanded to explicitly mention the group operations, coordinate charts, etc.).  Similarly, the nilsequence should not be represented solely as $F(g^n  x)$, but should really be the octuplet $(G, \Gamma, G/\Gamma, d_{G/\Gamma}, g, x, F, (n \mapsto F(g^n  x)))$.  However we shall continue to abuse notation in order to simplify the exposition.}
each nilmanifold $G/\Gamma$ with a smooth Riemannian metric $d_{G/\Gamma}$.  We then define the \emph{Lipschitz constant} of a nilsequence $F(g^n  x)$ to be the Lipschitz constant of $F$.

\begin{remark} Note that the Lipschitz constant of a nilsequence depends on the choice of metric $d_{G/\Gamma}$ one places on the nilmanifold; there is no obvious canonical metric to assign to any given nilmanifold, and so the Lipschitz constant is a somewhat arbitrary quantity.  However if one replaces the metric with another smooth Riemannian metric then from the compactness of $G/\Gamma$ we see that the Lipschitz constant is only affected by at most a multiplicative constant.  One could replace the Lipschitz constant here by other quantitative measures of regularity, such as H\"older continuity norms or $C^k$ norms, but this will not significantly affect the statements of the conjectures here, basically because a function which is controlled in one of these norms can be approximated in a quantitative manner as the uniform limit of functions controlled in any other of these norms.
\end{remark}

\begin{remark}\label{algebra-remark} The Lipschitz nilsequences form an algebra in the following sense: if $f(n)$ is an $s$-step nilsequence on $G/\Gamma$ with Lipschitz constant $M$, and $\widetilde{f}(n)$ is an $s$-step nilsequence on $\widetilde{G}/\widetilde{\Gamma}$ with Lipschitz constant $\widetilde{M}$, and both nilsequences are bounded by $O(1)$, then $f(n)\pm \widetilde{f}(n)$ or $f(n) \widetilde{f}(n)$ is an $s$-step nilsequence on the product nilmanifold $(G/\Gamma) \times (\widetilde{G}/\widetilde{\Gamma})$ with Lipschitz constant $O_{M,\widetilde{M}}(1)$.
However, nilsequences as we have defined them are not closed under uniform limits. This leads to a slight conflict between the nomenclature of the present paper and that of (for example) \cite{bergelson-host-kra}. In that paper the objects we have called nilsequences are referred to as \emph{basic} nilsequences; a nilsequence is then a uniform limit of basic nilsequences. Since our analysis is essentially finitary in nature we will not make any further mention of this distinction.
\end{remark}

\textsc{The inverse Gowers-norm  conjecture.} An important feature of $s$-step nilmanifolds is that they have significant ``constraints'' connecting arithmetic progressions of length $s+2$, or cubes of dimension $s+1$.  Roughly speaking, given the first $s+1$ elements $x, g \cdot x, g^2\cdot x, \ldots, g^s \cdot x$
of a progression in an $s$-step nilmanifold $G/\Gamma$, the next element $g^{s+1}\cdot x$ of the progression and all further elements are essentially completely determined as continuous functions of these first $s+1$ elements. For a precise formulation of
this assertion see \cite[Lemma 12.7]{green-tao-u3inverse}.  Similarly, when considering an $s$-dimensional ``cube''
$\{ g_1^{\omega_1} \ldots g_s^{\omega_s} \cdot x: (\omega_1,\ldots,\omega_s) \in \{0,1\}^s \}$ in $G/\Gamma$, the final
vertex $g_1 \ldots g_s \cdot x$ of this cube is essentially a continuous function of the other $2^{s-1}$ elements of this cube.
See Appendix \ref{nil-app} for more precise formulations of this statement, which we will make heavy use of in this paper. As a consequence of either of these facts,
we can relate nilsequences to the $U^{s+1}$ norm. The next result is in this direction, but it is not sufficiently general for our later applications. We state it now to introduce the concept of nilsequences obstructing uniformity, and because it can be proved using earlier results.

\begin{proposition}[Nilsequences obstruct uniformity]\label{nil-gow} Let $s \geq 1$ be an integer and let $\delta \in (0,1)$ be real. Let $G/\Gamma = (G/\Gamma, d_{G/\Gamma})$ be an $s$-step nilmanifold with some fixed smooth metric $d_{G/\Gamma}$, and let $(F(g^n  x))_{n \in \N}$ be a bounded $s$-step nilsequence with Lipschitz constant at most $M$.  Let $f : [N] \rightarrow [-1,1]$ be a function for which
\[ \E_{n \in [N]} f(n) F(g^n\cdot  x) \geq \delta.\]
Then we have
\[ \Vert f \Vert_{U^{s+1}[N]} \gg_{s,\delta,M,G/\Gamma} 1.\]
\end{proposition}

\begin{proof}
See \cite[Prop. 12.6]{green-tao-u3inverse}. The lower bound arising in that proposition was stated to depend on the continuous function $F : G/\Gamma \rightarrow \C$, and not just on $\Vert F \Vert_{\Lip}$. However, an examination of the proof reveals that the argument can be made uniform in $F$, for a given value of $\Vert F \Vert_{\Lip}$.
\end{proof}

\begin{remark} It turns out that one can relax the assumption that $f$ be uniformly bounded, requiring only that $f$ be bounded in $L^1$ norm; see Corollary \ref{nil-gow-cor}.
\end{remark}

The inverse Gowers-norm conjecture is an assertion in the converse direction, that nilsequences are the \emph{only} obstruction
to uniformity.  More precisely, we have for each $s \geq 1$ the following conjecture:

\begin{conjecture}[$\GI(s)$ conjecture]\label{gow-inv-conj}
Suppose that $0 < \delta \leq 1$.  Then there exists a finite collection ${\mathcal M}_{s,\delta}$ of $s$-step nilmanifolds
$G/\Gamma = (G/\Gamma,d_{G/\Gamma})$ with the following property. Given any $N$ and any
$f : [N] \rightarrow [-1,1]$ such that
\[ \Vert f \Vert_{U^{s+1}[N]} \geq \delta,\]
there is a nilmanifold $G/\Gamma \in {\mathcal M}_{s,\delta}$ and a $1$-bounded $s$-step 
nilsequence $(F(g^n  x))_{n \in \N}$ on it with Lipschitz constant $O_{s,\delta}(1)$, such that 
\[ |\E_{n \in [N]} f(n) F(g^n x)| \gg_{s,\delta} 1.\]
\end{conjecture}

This conjecture in this form is due to the authors. It was hinted at in \cite[\S 13]{green-tao-u3inverse} and is being stated formally for the first time here. The evidence in favour of it is strong. First of all we know that the cases $s = 1,2$ are true. The case $s = 1$ is an exercise in harmonic analysis. Indeed in this case one can take $G/\Gamma$ to just be the standard unit circle $\R/\Z$, so that ${\mathcal M}_{1,\delta}$ is a singleton set independent of $\delta$. The case $s = 2$ was established, with some effort, in \cite{green-tao-u3inverse} and is stated in Proposition \ref{gi2-prop} below. Note that things are not so simple when $s > 1$, and it is known that as $\delta$ decreases to zero, the collection ${\mathcal M}_{s,\delta}$ of nilmanifolds $G/\Gamma$ that one must employ must have cardinality going to infinity\footnote{This seems to be related to the fact, known to the ergodic theorists, that the inverse limit of $1$-step nilsystems is a $1$-step nilsystem, but the same is not true for $s$-step nilsystems, $s \geq 2$. }. 

\begin{proposition}[The $\GI(2)$ conjecture, {\cite{green-tao-u3inverse}}]\label{gi2-prop} 
The $\GI(2)$ conjecture holds in the form stated above. In fact the group $G$ may be taken to be a product of $O(\delta^{-O(1)})$ Heisenberg groups $\left(\begin{smallmatrix} 1 & \R & \R \\ 0 & 1 & \R \\ 0 & 0 & 1\end{smallmatrix}\right)$, and the discrete cocompact subgroup $\Gamma$ may be taken to be a product of copies of $\left(\begin{smallmatrix} 1 & \Z & \Z \\ 0 & 1 & \Z \\ 0 & 0 & 1\end{smallmatrix}\right)$.
\end{proposition}

\begin{proof} This is almost \cite[Thm. 12.8]{green-tao-u3inverse}. In that theorem, a nilsequence was constructed in a somewhat \emph{ad hoc} manner from another type of object, a generalised quadratic phase. In the argument of that paper, however, the nilpotent groups constructed were not all Heisenberg groups. Some of them were isomorphic to $\R^2 \times \Z$, which is not connected and hence, with our definition, cannot be used to construct a nilmanifold.

More precisely, in the proof of \cite[Thm. 12.8]{green-tao-u3inverse} it is shown that if $\Vert f \Vert_{U^3} \geq \delta$ then 
\[ \big|\E_{n \in [N]} f(n) F_1(g^n  x) e(n^2 \theta)\big| \gg \exp(-\delta^{-O(1)}),\]
where $F_1(g^n  x)$ is a product of nilsequences coming from $O(\delta^{-O(1)})$ Heisenberg groups (which are all connected and simply-connected), $\theta \in \R/\Z$, and
$e(x) := e^{2\pi i x}$. In \cite[Thm. 12.8]{green-tao-u3inverse} we proceeded by constructing $e(n^2 \theta)$ as a nilsequence coming from a skew torus which, being a quotient of the disconnected nilpotent Lie group 
$\left(\begin{smallmatrix} 1 & \R & \R \\ 0 & 1 & \Z \\ 0 & 0 & 1\end{smallmatrix}\right)$, is not immediately helpful in the present context. However we might just as easily have observed that
\[ \left(\begin{smallmatrix} 1 & -\theta & -\theta & \\ 0 & 1 & 2 \\ 0 & 0 & 1 \end{smallmatrix}\right)^n = \left(\begin{smallmatrix} 1 & -n\theta & -n^2 \theta \\ 0 & 1 & 2n \\ 0 & 0 & 1\end{smallmatrix}\right)\] which, upon quotienting by the right action of $\left(\begin{smallmatrix} 1 & \Z & \Z \\ 0 & 1 & \Z \\ 0 & 0 & 1\end{smallmatrix}\right)$, leads to
\[ \big[ \left(\begin{smallmatrix} 1 & -\theta & -\theta & \\ 0 & 1 & 2 \\ 0 & 0 & 1 \end{smallmatrix}\right)^n \big] = \big[ \left(\begin{smallmatrix} 1 & \{-n\theta\} & \{n^2 \theta\} \\ 0 & 1 & 0\\ 0 & 0 & 1\end{smallmatrix}\right)\big].\]
Here we have moved our matrix under the right action of $\Gamma$ so that it lies in the fundamental domain
\[ \mathcal{F} := \{\left( \begin{smallmatrix} 1 & x & z \\ 0 & 1 & y \\ 0 & 0 & 1 \end{smallmatrix}\right) : \textstyle -\frac{1}{2} < x,y,z \leq \frac{1}{2} \};\]
see \cite{green-tao-u3mobius} for further discussion. The fractional parts $\{t\}$ are chosen to lie in $(-\frac{1}{2},\frac{1}{2}]$. 

This almost exhibits $e(n^2 \theta)$ as a nilsequence coming from the Heisenberg group, but there is one small problem: the function
\[ \left( \begin{smallmatrix} 1 & x & z \\ 0 & 1 & y \\ 0 & 0 & 1 \end{smallmatrix}\right) \mapsto e(z)\] from $\mathcal{F}$ to $\C$
does not extend to a continuous function on $G/\Gamma$, since there are discontinuities on the boundary $\partial \mathcal{F}$.

To get around this one may introduce a smooth partition of unity $(\chi_{j})_{j \in J}$ on $(\R/\Z)^2$, where each function $\chi_j$ is supported on (say) a square of width $1/100$. Each function
\[ \left( \begin{smallmatrix} 1 & x & z \\ 0 & 1 & y \\ 0 & 0 & 1 \end{smallmatrix}\right) \mapsto \chi_j(x,y)e(z)\]
\emph{does} extend to a Lipschitz function on $G/\Gamma$.
This makes it clear that $e(n^2 \theta)$ may, after all, be realised as a nilsequence coming from a product of $O(1)$ Heisenberg groups.
\end{proof}

For higher values of $s$, the conjecture $\GI(s)$ remains open.  However, significant support in favour of this conjecture arises from the combinatorial and Fourier-analytic work of Gowers \cite{gowers-long-aps}, in which a ``local'' form of this conjecture was established in order to provide a new proof of Szemer\'edi's theorem. Further substantial support for the conjecture comes from the ergodic-theoretic work of Host-Kra \cite{host-kra}. 

\textsc{The \mobname conjecture.} Our main results are concerned with the von Mangoldt function $\Lambda(n)$ and with functions derived from $\Lambda$, such as $\Lambda'_{b,W}$. It turns out, however, to be convenient to rewrite
this function in terms of the closely related \emph{M\"obius function} $\mu: \Z \to \{-1,0,+1\}$, defined by setting $\mu(n) := (-1)^d$ when $n$ is the product of $d$ distinct primes, and $\mu(n) = 0$ otherwise.  The main advantage of doing so is that $\mu$ is a $1$-bounded function, whereas $\Lambda$ patently is not. As is well known, $\Lambda$ and $\mu$ are related by the identity
\begin{equation}\label{lambda-mu}
\Lambda(n) = \sum_{d|n} \mu(d) \log\frac{n}{d} = - \sum_{d|n} \mu(d) \log d 
\end{equation}
for all $n \geq 1$.  In principle this allows us to reduce the task of estimating correlations involving $\Lambda$ to that of estimating correlations involving $\mu$, although when doing so the unbounded weight $\log \frac{n}{d}$ and the summation over $d$ will introduce some dangerous factors of $O(\log N)$ which must be handled with some caution.

Suppose we formally apply Conjecture \ref{gow-inv-conj} to the task of proving Theorem \ref{gowers-norm}, ignoring for now the
significant issue that $\Lambda'_{b,W}-1$ is not uniformly bounded.  Then we expect to reduce this theorem to the assertion that $\Lambda'_{b,W}-1$ has small correlation with any $s$-step nilsequence.  In the light of \eqref{lambda-mu}, we expect this statement to
be related to the corresponding assertion for the M\"obius function $\mu$.  We formalise this latter statement as the following conjecture.

\begin{conjecture}[$\MN(s)$ conjecture]\label{mnconj} Let $G/\Gamma = (G/\Gamma,d_{G/\Gamma})$ be an $s$-step nilmanifold with smooth metric $d_{G/\Gamma}$, and let $(F(g^n  x))_{n \in [N]}$ be a bounded $s$-step nilsequence with Lipschitz constant $M$. Then we have the bound
\[ \big| \E_{n \leq N} \mu(n) F(g^n  x) \big| \ll_{A,M,G/\Gamma,s} \log^{-A} N\]
for any real number $A> 0$.
\end{conjecture}

\begin{remark} It is important to note that the implied constant is \emph{not} allowed to depend on $g$ and $x$. The case $s = 1$ can be reduced to a classical result of Davenport \cite{davenport-old}; see \cite[\S 6]{green-tao-u3mobius} for details. The case $s = 2$ was the main result of \cite{green-tao-u3mobius}.  The case $s > 2$ remains open; however, we certainly expect
$\MN(s)$ to be true in this case because of the \emph{M\"obius randomness} heuristic from analytic number theory, which states that $\mu$ exhibits a substantial degree of orthogonality to any suitably ``Lipschitz'' function. Moreover, it seems likely that the techniques we developed to prove $\MN(2)$ will eventually extend to cover $\MN(s)$, $s \geq 3$, as well. This is another ongoing area of research.  As is well known, even when $s = 1$ the current technology for establishing this conjecture yields ineffective implied constants in the $\ll_{A,M,G/\Gamma}$ due to our lack of knowledge regarding the existence of Siegel zeroes.   This ultimately makes the decay rates in the Main Theorem (and its corollaries) similarly ineffective.  If the GRH is assumed, the estimates do become effective. However they are still somewhat poor for $s \geq 2$, largely because the bounds in the $\GI(2)$ conjecture obtained in \cite{green-tao-u3inverse} are a little weaker than one might hope for.
\end{remark}

\section{Correlation estimates for M\"obius and Liouville}\label{mob-lio-sec}

Perhaps the heart of the present paper is \S \ref{sec5}, in which it is shown how, in certain circumstances, the requirement of $1$-boundedness can be dropped in the $\mbox{GI}(s)$ conjecture. This section is an aside to the main line of our argument, in which we use what we already have to obtain estimates similar to the generalised Hardy-Littlewood conjecture for the M\"obius function and the related 
\emph{Liouville function} $\lambda: \N \to \{-1,+1\}$, defined to be the unique completely multiplicative function such that $\lambda(p) = -1$ for all primes $p$.

\begin{proposition}[Correlation estimates for $\mu$ and $\lambda$]\label{mu-lam-cor}  Let $d, t, L$ be positive integers, let $N$ be a large positive integer parameter, and let $\Psi = (\psi_1,\ldots,\psi_t)$ be a system of affine-linear forms with size $\| \Psi\|_N \leq L$ and complexity at most $s$. Assume the $\GI(s)$ and $\MN(s)$ conjectures. Let $K \subset [-N,N]^d$ be a convex body.  Then we have
\begin{equation}\label{mob-eq}
\sum_{ n \in K \cap \Z^d} \prod_{i \in [t]} \mu( \psi_i( n) ) =  o_{s,t,d,L}(N^d)
\end{equation}
and 
\begin{equation}\label{lam-eq}
\sum_{ n \in K \cap \Z^d} \prod_{i \in [t]} \lambda( \psi_i( n) ) =  o_{s,t,d,L}(N^d).
\end{equation} 
\end{proposition}

\begin{remark} Note the lack of any local factors $\beta_p,\beta_{\infty}$. This makes Proposition \ref{mu-lam-cor} rather appealing from a certain point of view. It also provides an instance of the ``M\"obius randomness heuristic'' alluded to above.\end{remark}

\proof We begin by applying Proposition \ref{gvn}, the generalised von Neumann theorem. Since $\mu$ and $\lambda$ are $1$-bounded, this may be applied with the pseudorandom measure $\nu$ set equal to the constant function $1$, which is obviously $D$-pseudorandom for all $D$. We note that in this case the proof of Proposition \ref{gvn} that we give in Appendix \ref{gvn-app} is rather simpler than in the case of a more general $\nu$; specifically, one can use Corollary \ref{gczow-2} in place of Corollary \ref{gczow-gvn}, while the verification of \eqref{rgv-1}, \eqref{rgv-2} is trivial when $\nu=1$. 

The application of Proposition \ref{gvn} reduces \eqref{mob-eq} to the statement
\begin{equation}\label{eq71}
\Vert \mu  \Vert_{U^{s+1}[N]} = o_{s}(1).
\end{equation}
Applying the $\GI(s)$ conjecture, it is sufficient to establish that 
\begin{equation}\label{eq70} \E_{n \leq N} \mu(n) F(g^n  x) = o_{s,M,\delta}(1)\end{equation}
unifromly over all $G/\Gamma \in \mathcal{M}_{s,\delta}$ and all $1$-bounded $M$-Lipschitz nilsequences $(F(g^n  x))_{n \leq N}$ on $G/\Gamma$. Indeed the truth of such a statement implies, by the $\GI(s)$ conjecture, that $\Vert \mu \Vert_{U^{s+1}[N]} \leq \delta$, and one may then take $\delta$ arbitrarily small to deduce \eqref{eq71}. Recalling that $|\mathcal{M}_{s,\delta}| = O_{\delta,s}(1)$, we see that \eqref{eq70} follows immediately from (a weak form of) the $\MN(s)$ conjecture. This proves \eqref{mob-eq}.

The proof of \eqref{lam-eq} proceeds similarly. It suffices to establish the analogue of \eqref{eq70}, that is to say the bound
\begin{equation}\label{eq72}
\E_{n \leq N} \lambda(n) F(g^n  x) = o_{s,M,\delta}(1)
\end{equation}
uniformly over all $G/\Gamma \in \mathcal{M}_{s,\delta}$ and all $1$-bounded $M$-Lipschitz nilsequences $(F(g^n  x))_{n \leq N}$ on $G/\Gamma$. We begin by noting the identity
\[ \lambda(n) := \sum_{d^2 | n} \mu(\frac{n}{d^2}).\]
This implies that for any positive real $X$, any fixed $G/\Gamma \in \mathcal{M}_{s,\delta}$ and any $1$-bounded $M$-Lipschitz nilsequence $(F(g^n  x))_{n \leq N}$ on $G/\Gamma$ we have 
\begin{align}\nonumber
\E_{n \leq N} \lambda(n) F(g^n  x) & = \E_{n \leq N}\sum_{d^2 | n} \mu(\frac{n}{d^2}) F(g^n  x) \\ \nonumber & = \sum_{d \leq X} \E_{n \leq N} 1_{d^2 | n} \mu(\frac{n}{d^2}) F(g^n  x) + \sum_{d > X} \E_{n \leq N} 1_{d^2 | n} \mu(\frac{n}{d^2}) F(g^n  x) \\ & =  \sum_{d \leq X} \E_{k \leq N/d^2} \mu(k) F(g^{d^2 k}  x) + O(X^{-1}).\label{eq74}
\end{align}
By replacing $g$ by $g^{d^2}$ in the $\MN(s)$ conjecture 
we obtain the bound
\[ \E_{k \leq N/d^2} \mu(k) F(g^{d^2 k}  x) = o_{G/\Gamma,M,d}(1).\]
Substituting into \eqref{eq74} we obtain
\[ \E_{n \leq N} \lambda(n) F(g^n  x) = o_{G/\Gamma,M,X}(1) + O(X^{-1}).\]
Let $\eps > 0$ be arbitrary. Taking $X := 1/\eps$, we may make this expression smaller than a constant times $\eps$ by taking $N$ sufficiently large. This implies that
\[ \E_{n \leq N} \lambda(n) F(g^n  x) = o_{G/\Gamma,M}(1).\]
Recalling once more that $|\mathcal{M}_{s,\delta}| = O_{s,\delta}(1)$, we therefore obtain \eqref{eq72} and hence \eqref{lam-eq}.\endproof

Let us remark that, as with the Main Theorem, Proposition \ref{mu-lam-cor} is unconditional in the cases $s = 1,2$. 

We conclude with a mention of a conjecture of Chowla \cite{chowla}, which asserts that $\lambda$ is uniformly distributed on any polynomial, thus for instance
\begin{equation}\label{chowla}
\E_{y_1, y_2 \leq N} \lambda(P(y_1, y_2)) = o_P(1)
\end{equation}
for any polynomial $P: \N \times \N \to \N$ of two variables. Our results imply (for instance) the following case of Chowla's conjecture.

\begin{proposition}\label{liouville}  Let $P: \N \times \N \to \N$ be a polynomial of degree at most $4$ which is the product of homogeneous linear factors over $\Q$, and which is not a rational multiple of a perfect square.  Then we have
\[ \E_{y_1, y_2 \leq N} \lambda(P(y_1, y_2)) = o_P(1).\]
\end{proposition}
The proof is immediate from \eqref{lam-eq} and the complete multiplicativity of $\lambda$; note that we can easily eliminate any repeated factors in $P$ and so the system of linear forms associated to $P$ will be non-degenerate. We remark that this conjecture was also recently verified for all homogeneous polynomials of degree at most three in \cite{helfgott1,helfgott2}. 
Removing the homogeneity assumption looks hopeless with current technology; the case $P(y_1,y_2) = y_1(y_1+2)$ is already roughly of the same order of difficulty as the twin prime conjecture.

\section{Transferring the inverse Gowers-norm conjecture}\label{sec5}

Recall that we are trying to use the inverse Gowers-norm and M\"obius and nilsequences conjectures to prove Theorem \ref{gowers-norm}.
We cannot apply the Gowers Inverse conjecture directly to prove Theorem \ref{gowers-norm}, because $\Lambda'_{b,W}-1$ is not bounded
uniformly in $N$.  The difficulty here is similar to that encountered in \cite{green-tao-longprimeaps}, in which Szemer\'edi's theorem, which ostensibly only establishes multiple recurrence bounds for bounded functions, needed to be extended to an unbounded
function such as $\Lambda'_{1,W}$.  We will use a similar resolution to that in \cite{green-tao-longprimeaps}, namely to \emph{transfer} the inverse Gowers-norm conjecture to the situation of a function bounded by a pseudorandom measure.  More
precisely, the purpose of this section is to prove the following result.

\begin{proposition}[Relative inverse Gowers-norm conjecture]\label{gow-ps}
Assume the $\GI(s)$ conjecture. For any $0 < \delta \leq 1$ and any $C \geq 20$, there exists a finite collection ${\mathcal M}_{s,\delta,C}$ of nilmanifolds $G/\Gamma = (G/\Gamma,d_{G/\Gamma})$ with the following property.  Let $N \geq 1$, suppose that $N' \in [CN,2CN]$ is a prime, that $\nu: \Z_{N'} \to \R^+$ 
is an $(s+2)2^{s+1}$-pseudorandom measure, that
$f : [N] \rightarrow \R$ is a function with $|f(n)| \leq \nu(n)$
for all $n \in [N]$ and that $\Vert f \Vert_{U^{s+1}[N]} \geq \delta$. Then there exists $G/\Gamma \in {\mathcal M}_{s,\delta,C}$ together with a $1$-bounded $s$-step nilsequence $(F(g^n x))_{n \in \Z}$ with Lipschitz constant $O_{s,\delta,C}(1)$,
such that
\[ |\E_{n \leq N} f(n) F(g^n x)| \gg_{s,C,\delta} 1.\]
\end{proposition}

\begin{remarks}
This looks significantly more complicated than the ordinary $\GI(s)$ conjecture, but this is something of an illusion. Most of the complexity comes from the need for the additional dependence on $C$. A largeish value of $C$ might be required in order to construct an appropriate pseudorandom measure $\nu$ on $\Z_{N'}$ (cf. Proposition \ref{pseudodom}) and so we leave $C$ unspecified in this proposition. 
\end{remarks}

In view of Proposition \ref{gow-ps} and Proposition \ref{pseudodom}, it is not hard to see that Theorem \ref{gowers-norm}, and hence the Main Theorem, follows from the next proposition. All one need do is choose $C := \max(C_0((s+2)2^{s+1}),20)$, where $C_0$ is the function appearing in Proposition \ref{pseudodom}. This ensures that an appropriate pseudorandom measure $\nu$ can be constructed.

\begin{proposition}[W-tricked von Mangoldt orthogonal to nilsequences]\label{manortho}
Let $s \geq 1$, and assume the $\MN(s)$ conjecture.
Let $G/\Gamma = (G/\Gamma,d_{G/\Gamma})$ be an $s$-step nilmanifold with smooth metric $d_{G/\Gamma}$, and let $(F(g^n x))_{n \in [N]}$ be a bounded $s$-step nilsequence with Lipschitz constant $M$. Let $b \in [W]$ be coprime to $W$.  Then we have the bound
\[  \E_{n \in [N]} (\Lambda'_{b,W}(n)-1) F(g^n x)  = o_{M,G/\Gamma,s}(1).\]
\end{proposition}

\begin{remark} In principle, Proposition \ref{manortho} is substantially easier to establish than the preceding reductions of the Main Theorem, such as Theorem \ref{gowers-norm}.  This is because we are now computing the correlation of $\Lambda$ (or $\Lambda'_{b,W}-1$) with respect to a ``low complexity'' sequence $F(g^n x)$, rather than the more complicated task of computing a multilinear correlation of $\Lambda$ with \emph{itself}.  In particular one can now hope to use tools such as Vinogradov's method to establish this proposition.  Indeed, the computation of exponential sums such as $\sum_{n \in [N]} \Lambda(n) e(\alpha n)$, or more generally $\sum_{n \in [N]} \Lambda(n) e(\alpha n^k)$, are essentially model cases of Proposition \ref{manortho} and are well-known to be treatable by Vinogradov's method. However, Proposition \ref{manortho} is somewhat more general as it also (for example) asserts some control on generalised polynomial exponential sums such as $\sum_{n \in [N]} \Lambda(n) e( \alpha n \lfloor \beta n \rfloor )$, where $\lfloor \cdot \rfloor$ is the greatest integer function.  See \cite{green-tao-u3mobius} for further discussion of the link between such functions and $2$-step nilsequences.  Thus we see that the inverse Gowers-norm conjecture $\GI(s)$ is 
a powerful tool for establishing bounds on the Gowers norms $U^{s+1}$, and thence to all multilinear averages of complexity 
at most $s$.
\end{remark}

We prove Proposition \ref{manortho} in later sections.  For the remainder of this section we derive Proposition \ref{gow-ps} from the inverse Gowers-norm conjecture.

\textsc{A Koopman-von Neumann theorem.\footnote{This term has something in common with the term ``generalised von Neumann theorem'' in that it originally came from analogies with ergodic theory. We now use it in our work to describe a range of theorems whose general aim is to decompose a given function $f$ into the sum of a function $f_1$ which is somehow less complicated than $f$, together with an error $f_2$ which is small in some Gowers norm.}} The primary tool in deducing Proposition \ref{gow-ps} from the Gowers Inverse conjecture is the following structure theorem, which allows us to decompose an arbitrary function $f$ which is bounded pointwise by $\nu$ into a bounded function and a Gowers-uniform function.

\begin{proposition}[Koopman -- von Neumann theorem]\label{kvn-decomp}
Let $s \geq 1$ and let $N' \geq N \geq 1$ be an integer. Suppose that $\nu$ is an $(s+2)2^{s+1}$-pseudorandom measure on $\Z_{N'}$, and that $f : \Z_{N'} \rightarrow \R$ is a function such that $|f(n)| \leq \nu(n)$ pointwise. Then we may decompose $f = f_1 + f_2$, where
\begin{equation}\label{tf}
 \sup_{n \in \Z_{N'}} |f_1(n)| \leq 1 
 \end{equation}
and
\begin{equation}\label{tf2}
 \Vert f_2 \Vert_{U^{s+1}(\Z_{N'})} = o(1).
 \end{equation}
If furthermore $f$ is supported in $\{-N,\ldots,N\}$ for some $N < N'/10$, then we may arrange matters so that $f_1$ and $f_2$ are both supported on $\{-2N, \ldots, 2N \}$.
\end{proposition}

\begin{remark} Informally, this theorem asserts that in the $U^{s+1}$ topology, bounded functions are dense in the class of functions bounded by $\nu$.  This fact (and refinements thereof), in conjunction with generalised von Neumann theorems such as Proposition \ref{gvn}, underlie the ``transference principle'' from \cite{green-tao-longprimeaps} which allow one to convert results for multilinear averages of $1$-bounded functions to results for multilinear averages of functions bounded by a pseudorandom measure.
This principle is essential for our arguments here, as it allows us in many cases to manipulate functions such as $\Lambda_{b,W}$ as if they were uniformly bounded.
\end{remark}

\begin{proof} Let us first make the observation that we can weaken \eqref{tf} to
\begin{equation}\label{tf-0}
 \sup_{n \in \Z_{N'}} |f_1(n)| \leq 1 + o(1) 
\end{equation}
since one could simply transfer the $o(1)$ error in \eqref{tf-0} to the $f_2$ component afterwards, using
the triangle inequality on \eqref{tf2}.

We shall rely heavily on a similar result from \cite[Proposition 8.1]{green-tao-longprimeaps}.  Before we give this result we 
need some notation.

\begin{definition}[Conditional expectation] If $f: \Z_{N'} \to \R$ is a function and $1 \leq p \leq \infty$, we denote
$\|f\|_{L^p(\Z_{N'})} := (\E_{n \in \Z_{N'}} |f(n)|^p)^{1/p}$, with the usual convention that $\|f\|_{L^\infty(\Z_{N'})} := \sup_{n \in \Z_{N'}} |f(n)|$.  If $\B$ is a $\sigma$-algebra on $\Z_{N'}$, that is to say the Boolean algebra generated by the atoms of a partition of $\Z_{N'}$, we define the \emph{conditional expectation} $\E(f|\B)$ of $f$ relative to $\B$ to be the orthogonal projection in $L^2(\Z_{N'})$ from $f$ to the $\B$-measurable functions.
\end{definition}

In our current notation, Proposition 8.1 from \cite{green-tao-longprimeaps} asserts\footnote{In \cite{green-tao-longprimeaps} the result is only stated when $0 \leq f(n) \leq \nu(n)$, but exactly the same proof applies under the more general assumption that $|f(n)| \leq \nu(n)$.  In any case, in order to prove Proposition \ref{kvn-decomp} one could always decompose $f$ into non-negative and negative parts $f^+ + f^-$ and follow the proof for each part separately. The key point to note is that the function $f^+_1$ is non-negative, whilst $f^-_1 \leq 0$. Thus $f_1 = f^+_1 + f^-_1$ satisfies the requisite $L^{\infty}$ bound \eqref{tf-0}.} the following.

\begin{quote-thm-1} Suppose that $N' \geq N$ and that $\nu : \Z_{N'} \rightarrow \R_{\geq 0}$ is an $(s+2)2^{s+1}$-pseudorandom measure.  Let $f: \Z_{N'} \to \R$ be such that
 $|f(n)| \leq \nu(n)$ for all $n \in \Z_{N'}$. Let $\eps \in (0,1)$ be a small parameter, and assume $N'$ is sufficiently large depending on $\eps$. Then there exists a $\sigma$-algebra $\B$ and an exceptional set $\Omega \in \B$ such that
\begin{itemize}
\item \textup{(smallness condition)} 
\begin{equation}\label{small-exceptions-final}
 \E_{\Z_{N'}}( \nu 1_{\Omega} ) = o_{\eps}(1);
\end{equation}
\item \textup{($\nu$ is uniformly distributed outside of $\Omega$)}
\begin{equation}\label{mu-regular-final}
 \| (1 - 1_{\Omega}) \E( \nu - 1 | \B) \|_{L^\infty(\Z_{N'})} = o_{\eps}(1)
\end{equation} and
\item \textup{(Gowers uniformity estimate)}
\begin{equation}\label{unif-ok}
 \| (1 - 1_{\Omega}) (f - \E(f|\B)) \|_{U^{s+1}(\Z_{N'})} \leq \eps^{1/2^{s+2}} = \kappa_s(\epsilon).
\end{equation}
\end{itemize}
\end{quote-thm-1}

Let $\eps$ be chosen later (it will eventually be a slowly decaying function of $N$).
If $N$ is sufficiently large depending on $\eps$, we can invoke the above theorem.  
Write
\[ f = f_1 + f_2 = f_1 + f^{(1)}_2 + f_2^{(2)},\]
where
\[ f_1 := (1 - 1_{\Omega}) \E(f | \B),\]
\[ f_2^{(1)} := (1 - 1_{\Omega}) (f - \E(f | \B))\]
and 
\[ f_2^{(2)} := 1_{\Omega} f.\]
Then by \eqref{mu-regular-final} we have 
\begin{equation}\label{f1-bound} \Vert f_1 \Vert_{L^\infty(\Z_{N'})} \leq 1 + o_{\eps}(1).\end{equation}
Also, by \eqref{unif-ok} we have 
\begin{equation}\label{f2-bound} \Vert f_2^{(1)} \Vert_{U^{s+1}(\Z_{N'})} = \kappa_s(\epsilon).\end{equation}
Next, we claim that
\begin{equation}\label{to-prove-10} \Vert f_2^{(2)} \Vert_{U^{s+1}(\Z_{N'})} = o_{\eps}(1).
\end{equation} 
To see this, first note that from \eqref{small-exceptions-final} we have
\begin{equation}\label{f3-bound} \Vert f_2^{(2)} \Vert_{L^1(\Z_{N'})} = o_{\eps}(1).\end{equation}
Secondly, we prove that for functions $g$ for which $|g|$ is bounded pointwise by a pseudorandom measure $\nu$, the $L^1(\Z_{N'})$ norm controls the $U^{s+1}(\Z_{N'})$-norm.
Indeed for such a function we have
\begin{align*}
\Vert g \Vert_{U^{s+1}(\Z_{N'})}^{2^{s+1}} &= \E_{n \in \Z_{N'}, {h} \in \Z_{N'}^{s+1}}g(n) \prod_{\substack{\omega \in \{0,1\}^{s+1} \\ \omega \neq 0}} g(n + \omega \cdot {h}) \\ &\leq  \E_{n \in \Z_{N'}} |g(n)| \sup_n \E_{{h} \in \Z_{N'}^{s+1}} \prod_{\substack{\omega \in \{0,1\}^{s+1} \\ \omega \neq 0}} \nu(n + \omega \cdot {h}) \\ & = \Vert \mathcal{D} \nu \Vert_{L^\infty(\Z_{N'})} \Vert g \Vert_{L^1(\Z_{N'})} ,\end{align*}
where 
$$\mathcal{D}\nu(n) := \prod_{\substack{\omega \in \{0,1\}^{s+1} \\ \omega \neq 0}} \nu(n + \omega \cdot {h})$$
is the \emph{dual function} associated to $\nu$. However a simple application of the linear forms condition, given in detail in \cite[Lemma 6.1]{green-tao-longprimeaps}, confirms that
\[ \Vert \mathcal{D}\nu \Vert_{L^\infty(\Z_{N'})} \leq 1 + o(1).\]
This concludes the proof of \eqref{to-prove-10}.  From this, \eqref{f2-bound}, and the triangle inequality for the $U^{s+1}(\Z_{N'})$ norm we conclude that
$$ \| f_2 \|_{U^{s+1}(\Z_{N'})} \leq o_{\eps}(1) + \kappa(\epsilon).$$
Choosing $\eps$ to be a sufficiently slowly decaying function of $N$ we obtain the first part of Proposition \ref{kvn-decomp}.

It remains to deal with the situation where $f$ is supported\footnote{An alternate way to proceed at this point is to modify the proof of \cite[Proposition 8.1]{green-tao-longprimeaps}, where the $\sigma$-algebra $\B$ is initialised not at the trivial factor, but rather at the factor generated by $\{-N,\ldots,N\}$.} in $\{-N,\ldots,N\}$.  
We can write $f(n) = f(n) \psi(n)$, where $\psi: \Z_{N'} \to [0,1]$ equals $1$ on $\{-N,\ldots,N\}$, vanishes outside of $\{-2N,\ldots,2N\}$ and interpolates smoothly in the range $N \leq |n| \leq 2N$. One could, for example, take $\psi$ to be a de la Vall\'ee Poussin kernel.  If $f = f_1 + f_2$ is the previous decomposition, then upon multiplying by $\psi$ we obtain $f = \tilde f_1 + \tilde f_2$, where $\tilde f_1 := f_1 \psi$ and $\tilde f_2 := f_2 \psi$.  The function $\tilde f_1$ continues to enjoy the bound \eqref{tf-0} but now also has the desired support property.  To confirm that $\tilde f_2$ enjoys the bound \eqref{tf2}, simply use Fourier series to break $\psi$ up as a rapidly convergent linear combination of linear phases $e(n \xi/N)$, and use the triangle inequality combined with the phase invariance \eqref{character} of the $U^{s+1}$ norm. This concludes the proof of Proposition \ref{kvn-decomp}.
\end{proof}

\textsc{Proof of Proposition \ref{gow-ps}.} Suppose that $N' \in [CN,2CN]$ is prime, that $\nu : \Z_{N'} \rightarrow \R$ is an $(s+2)2^{s+1}$-pseudorandom measure, that $f : [N] \rightarrow \R$ is a function with $|f(n)| \leq \nu(n)$ for all $n \in [N]$ and that $\Vert f \Vert_{U^{s+1}[N]} \geq \delta$. Applying Proposition \ref{kvn-decomp} we may decompose
\[ f = f_1 + f_2,\]
where $\Vert f_1 \Vert_{L^\infty(\Z_{N'})} \leq 1$ and $\Vert f_2 \Vert_{U^{s+1}(\Z_{N'})} = o(1)$. Since $C > 10$, we may further assume that both $f_1$ and $f_2$ are supported in $\{-2N,\dots, 2N\}$. By Lemma \ref{normlemma} the assumption that $\Vert f \Vert_{U^{s+1}[N]} \geq \delta$ implies that $\Vert f \Vert_{U^{s+1}(\Z_{N'})} \gg_{C,s} \delta$, and hence that $\Vert f_1 \Vert_{U^{s+1}(\Z_{N'})} \gg_{C,s} \delta$. Applying Lemma \ref{normlemma} once more, we conclude that $\Vert f_1 \Vert_{U^{s+1}(\{-2N,\dots,2N\})} \gg_{C,s} \delta$. 

We now apply the inverse Gowers-norm conjecture $\GI(s)$, translating $\{-2N,\ldots,2N\}$ to the interval $[4N+1]$, to conclude that there exists an $s$-step nilmanfold $G/\Gamma = (G/\Gamma,d_{G/\Gamma})$ from
a fixed finite collection $G/\Gamma \in {\mathcal M}_{s,\delta,C}$, together with a bounded $s$-step 
nilsequence $(F(g^n x))_{n \in \N}$ generated by this nilmanifold and with Lipschitz constant $O_{s,\delta,C}(1)$, such that 
\[ |\E_{-2N \leq n \leq 2N} f_1(n) F(g^n x)| \gg_{s,\delta,C} 1.\]
On the other hand, from \eqref{tf2} and the contrapositive of Proposition \ref{nil-gow} we have
\[ |\E_{-2N \leq n \leq 2N} f_2(n) F(g^n x)| = o_{G/\Gamma,s,\delta,C}(1).\]
If $N \geq N_0(s,\delta,C)$ is large depending on $s$, $\delta$ and $C$, we conclude that
\[ |\E_{-2N \leq n \leq 2N} f(n) F(g^n x)| \gg_{s,\delta,C} 1\]
and the claim follows (since $f$ is supported on $[N]$).

If by contrast $N = O_{s,\delta,C}(1)$ then the claim is trivial, since all norms on $[N]$ are then equivalent up to factors of $O_N(1) = O_{s,\delta,C}(1)$, and 
all functions on $[N]$ can be expressed as nilsequences (say on the torus $\R/\Z$) with
Lipschitz constant $O_N(1) = O_{s,\delta,C}(1)$. \endproof

\section{Averaging the nilsequence}\label{average-sec}

To summarise so far, we have reduced the task of showing that the $\GI(s)$ conjecture implies the Main Theorem to 
the much easier task of establishing Proposition \ref{manortho}.  This, recall, is an estimate on the correlation between
the number-theoretic function $\Lambda'_{b,W}(n)-1$ and the nilsequence $F(g^n x)$.

The purpose of this section is to perform a rather technical modification to the nilsequence $F(g^n x)$, which is necessary
for the following reason.
At a later stage in the proof we would like to discard certain ``small'' components of the function $\Lambda'_{b,W}(n)-1$ from this correlation.  Some of these components will be easy to discard; for instance, any error which is small in $L^1$ norm will be easily removed since the nilsequence is bounded.  However, there will be one component of $\Lambda'_{b,W}(n)-1$ that we shall encounter
(namely, the term arising from the ``smooth'' component $\Lambda^\sharp$ of the von Mangoldt function) which will not be small in $L^1$, but is instead small in the Gowers norm $U^{s+1}[N]$.  In principle, Proposition \ref{nil-gow} or Corollary \ref{nil-gow-cor} would allow us to safely drop such terms.  Unfortunately, a problem arises because the component of $\Lambda'_{b,W}(n)-1$ that we are trying to discard is not bounded, and we have also not been able to dominate this component by a pseudorandom measure or even to establish a bound for it in $L^1$.  To get around this problem, we need to improve the ``regularity'' of the nilsequence $F(g^n x)$. In particular we must convert it to an object which we can bound in the dual norm $U^{s+1}[N]^*$, defined as usual by the formula
$$ \|F\|_{U^{s+1}[N]^*} := \sup \{ |\E_{n \in [N]} f(n) F(n)|: \|f\|_{U^{s+1}[N]} \leq 1 \}.$$
This dual norm also appeared in \cite{green-tao-longprimeaps}, and plays a similar r\^ole there as it does here.

It would be very pleasant if every $s$-step nilsequence was automatically bounded in the $U^{s+1}[N]^*$ norm.  Unfortunately, this
statement is false even in the $s=1$ case, as in that case it amounts to a certain $l^{4/3}$ summability estimate on the Fourier coefficients of Lipschitz functions on a compact abelian group. There is no such estimate if the group is of sufficiently high dimension.  Of course one can rectify this by replacing the Lipschitz functions with smooth functions.  It seems likely that a similar claim is true for higher $s$, but it also seems likely that a proof would involve a finer analysis of the structure of nilmanifolds than we need for the rest of our argument.

Fortunately, however, we can achieve an adequate substitute result by replacing the concept of a nilsequence by its convex hull. Definition \ref{def11.1} provides a precise definition. 


\begin{definition}[Averaged nilsequences]\label{def11.1}  Let $G/\Gamma = (G/\Gamma,d_{G/\Gamma})$ be an $s$-step nilmanifold, and let $M > 0$.
An $s$-step \emph{averaged nilsequence} on $G/\Gamma$ with Lipschitz constant at most $M$ is a function $F(n)$ having the form
\[ F(n) = \E_{i \in I} F_i(g_i^n x_i),\]
where $I$ is some finite index set, 
and for each $i$, $F_i(g_i^n x_i)$ is a bounded $s$-step nilsequence on $G/\Gamma$ with Lipschitz constant at most $M$.
\end{definition}

\begin{remark}
An averaged nilsequence of the type just described is a genuine nilsequence on the nilmanifold $(G/\Gamma)^I$. However the averaging set $I$ will, in applications, have size comparable to $N$ and so in our finitary world these averaged nilsequences should be thought of as a strict generalisation of the notion of a nilsequence. Were it not for the desire to avoid issues of measurablility, we might even have replaced the finite averaging operator $\E_{i \in I}$ by an integration over a suitable probability space.
\end{remark}

We now state the crucial technical lemma we need, which allows us to replace a nilsequence by an averaged nilsequence with a good
$U^{s+1}[N]^*$ bound.

\begin{proposition}[Decomposition of nilsequences]\label{decomp-prop}
Let $G/\Gamma = (G/\Gamma,d_{G/\Gamma})$ be an $s$-step nilmanifold, and let $M > 0$.
Suppose that $(F(g^n x))_{n \in \N}$ is a bounded $s$-step nilsequence on $G/\Gamma$ with Lipschitz constant at most $M$. Let $\eps \in (0,1)$ and suppose that $N \geq 1$. Then we may effect the decomposition
\begin{equation}\label{fgn}
 F(g^n x) = F_1(n) + F_2(n),
\end{equation}
where $F_1: \N \to [-1,1]$ is an averaged nilsequence on $(G/\Gamma)^{2^{s+1}-1}$ with Lipschitz constant $O_{M,\eps,G/\Gamma}(1)$ and obeying the dual norm bound
\begin{equation}\label{gow-dual-bd} \Vert F_1 \Vert_{U^{s+1}[N]^*} \ll_{M,\eps,G/\Gamma} 1,\end{equation}
while $F_2: \N \to \R$ obeys the uniform bound
\begin{equation}\label{fun}
\Vert F_2\Vert_{\infty} = O(\eps). 
\end{equation}
\end{proposition}

\begin{remark} At present, our decomposition \eqref{fgn} depends on the parameter $N$.  It is possible to modify the argument below in such a way that the decomposition is independent of $N$, but this requires generalising the notion of an averaged nilsequence by replacing the averaging over a finite set $I$ with an integral over a continuous probability measure.  As this introduces some minor technical issues such as measurability, we shall settle for the slightly weaker formulation of Proposition \ref{decomp-prop} given above, as it still suffices for our application.
\end{remark}

We shall prove Proposition \ref{decomp-prop} shortly.  Assuming it for the moment, we may make yet another reduction of the Main Theorem. This we do by reducing Proposition \ref{manortho} (which, as we have already shown, implies the Main Theorem) to the following result.

\begin{proposition}[W-tricked $\Lambda$ orthogonal to averaged nilsequences]\label{manortho-avg}
Let $s \geq 1$, and assume the $\MN(s)$ conjecture. 
Let $G/\Gamma = (G/\Gamma,d_{G/\Gamma})$ be an $s$-step nilmanifold with smooth metric $d_{G/\Gamma}$, and let $F_1(n)$ be an averaged $s$-step nilsequence with Lipschitz constant $M$. Let $b \in [W]$ be coprime to $W$.  
Suppose we also have the dual norm bound 
\begin{equation}\label{f1un}
 \Vert F_1 \Vert_{U^{s+1}[N]^*} \leq M'.
\end{equation}
Then we have the bound
\[  \E_{n \in [N]} (\Lambda'_{b,W}(n)-1) F_1(n) = o_{M,M',G/\Gamma,s}(1).\]
\end{proposition}

Indeed, to deduce Proposition \ref{manortho} from Proposition \ref{manortho-avg}, let $\eps \in (0,1)$ be arbitrary and apply Proposition \ref{decomp-prop}.  The contribution of $F_2$ will be bounded by $O(\eps) + o_{\eps}(1)$ thanks to
\eqref{siegel} and \eqref{fun}.  The contribution of $F_1$ can be controlled using Proposition \ref{manortho-avg}. Putting these estimates together leads to the bound
\[  \E_{n \in [N]} (\Lambda'_{b,W}(n)-1) F(g^n x)  = o_{M,G/\Gamma,s,\eps}(1) + o_{\eps}(1)
+ O(\eps).\]
Letting $\eps$ go to zero sufficiently slowly, we obtain the claim.

In later sections we shall prove Proposition \ref{manortho-avg}.  For now we turn to the task of
proving Proposition \ref{decomp-prop}.

\begin{proof}[Proof of Proposition \ref{decomp-prop}]
Fix $G/\Gamma$, $s$, $M$.  Observe that if we have proven the proposition for a single Lipschitz function $F$, 
then if we perturb $F$ in the $L^\infty$ norm by $\eps$ then the statement is still true for the perturbed function (with slightly worse implied constants in the $O()$ notation).  On the other hand, since $G/\Gamma$ is a compact metric
space, we know from the Arzel\`a-Ascoli theorem that the space of Lipschitz functions $F$ on $G/\Gamma$ with Lipschitz constant at most $M$ is equicontinuous and hence compact in the uniform topology.  In particular, it can be covered by finitely many balls in the uniform metric of radius $\eps/2$, say.  In view of this compactness\footnote{One could also use the compactness of $G/\Gamma$ to remove the requirement that all bounds be uniform in $x$.  However the parameter $g$ ranges over the non-compact group $G$ and cannot be eliminated so easily; the range of the parameter $n$ is similarly non-compact.  Thus we will be forced to look for constraints in the orbit $g^n x$ which are \emph{independent} of $g$ and $n$.  This helps motivate our introduction of cubes below.}, we see that it will suffice to establish the \emph{qualitative} version of the Proposition, namely given any continuous function $F$ (not necessarily Lipschitz) and any $\eps > 0$, we have a decomposition \eqref{fgn} for all $N \geq 1$, $g \in G$ and $x \in G/\Gamma$, where $F_1$ is an averaged nilsequence on $G/\Gamma$ with Lipschitz constant uniform in $g,x,N$, and with dual norm $\|F_1\|_{U^{s+1}[N]^*}$ bounded uniformly in $N,g,x$, and $F_2$ obeys the bound \eqref{fun}.

Fix $F$ and $\eps$. To proceed further we need to detect some ``constraints'' on the orbit $n \mapsto g^n x$ in $G/\Gamma$.
The most convenient framework for giving such constraints will be the $(s+1)$-dimensional parallelepipeds in $G/\Gamma$, as studied in \cite{host-kra}.

\begin{definition}[Parallelepipeds in nilmanifolds]  Let $(G/\Gamma)^{\{0,1\}^{s+1}}$ denote the space of all $2^{s+1}$-tuples $(x_\omega)_{\omega \in \{0,1\}^{s+1}}$.  
An \emph{$(s+1)$-dimensional parallelepiped} is any element of $(G/\Gamma)^{\{0,1\}^{s+1}}$ having
the form
$$ ( g^{n + \omega \cdot {h}} x )_{\omega \in \{0,1\}^{s+1}}$$
for some $g \in G$, $x \in G/\Gamma$, $n \in \Z$, and ${h} \in \Z^{s+1}$. Here, and for the remainder of the paper, we write $\omega \cdot h := \omega_1 h_1 + \dots + \omega_{s+1}h_{s+1}$ where $\omega = (\omega_1,\dots,\omega_{s+1})$ and $h = (h_1,\dots,h_{s+1})$.
\end{definition}

A fundamental property of $s$-step nilmanifolds is that the value of any one vertex of a parallelepiped (say, the zero vertex $x_{0^{s+1}}$, where $0^{s+1} := (0,\ldots,0)$) is determined ``continuously'' by all the other vertices.  In the following proposition, and for the remainder of the paper, write $\{0,1\}_*^{s+1} := \{0,1\}^{s+1} \setminus \{0^{s+1}\}$.

\begin{proposition}[Parallelepiped constraint]\label{cube} There exists a compact set \[ \Sigma \subseteq (G/\Gamma)^{\{0,1\}^{s+1}_*}\] and
a continuous function $P: \Sigma \to G/\Gamma$ such that, for any $(s+1)$-dimensional parallelepiped $(x_\omega)_{\omega \in \{0,1\}^{s+1}}$, we have $(x_\omega)_{\omega \in \{0,1\}^{s+1}_*} \in \Sigma$ and the constraint
$$ x_{0^{s+1}} = P( (x_\omega)_{\omega \in \{0,1\}^{s+1}_*} ).$$
\end{proposition}

This proposition is a topological and algebraic statement about the structure of nilmanifolds, and it was essentially proved in \cite{host-kra}. We supply a complete and self-contained proof in Appendix \ref{nil-app}, taking the opportunity to introduce the \emph{Host-Kra cube groups}. A closely related 
statement regarding arithmetic progressions in nilmanifolds appeared in \cite[Lemma 12.7]{green-tao-u3inverse}; results of this latter type seem to have been around in the ergodic theory community for some time and feature, for instance, in the papers of Furstenberg \cite{furst-von-neumann,furst-erdos-survey}.

For now, we shall simply illustrate this proposition with two model examples before continuing with the proof of Proposition \ref{decomp-prop}.

\begin{example}[Abelian shift] Take $s=1$, let $G$ be an abelian Lie group, and let $\Gamma$ be a cocompact lattice in $G$. Thus $G/\Gamma$ is a compact abelian Lie group, and any action of $g \in G$ on $G/\Gamma$ has the form of a shift $x \mapsto x + g$.  Of course, $G/\Gamma$ is a $1$-step nilmanifold.  A $2$-dimensional parallelepiped in this nilmanifold takes the form $(x+ng, x+(n+h_1)g, x+(n+h_2)g, x+(n+h_1+h_2)g)$.  The first vertex is a function of the other three. In the notation of Proposition \ref{cube} we can take $\Sigma := (\R/\Z)^3$ and $P: \Sigma \to G/\Gamma$ be the map
$P( y_{10}, y_{01}, y_{11} ) := y_{01} + y_{10} - y_{11}$ and we easily verify that $y_{00} = P( y_{10}, y_{01}, y_{11} )$ whenever $(y_{00}, y_{10}, y_{01}, y_{11})$ is a $2$-dimensional parallelepiped.
\end{example}

\begin{example}[Skew shift] For the sake of illustration, we consider a quotient $G/\Gamma$ where $G$ is $2$-step nilpotent but not connected. The way we have set things up in this paper, then, $G/\Gamma$ does not qualify as a nilmanifold; however one can modify this example so that it genuinely takes place in a nilmanifold (cf. the proof of Proposition \ref{gi2-prop}).

Set $G := \left( \begin{smallmatrix} 1 & \Z & \R \\ 0 & 1 & \R \\ 0 & 0 & 1\end{smallmatrix} \right)$ and $\Gamma := \left( \begin{smallmatrix} 1 & \Z & \Z \\ 0 & 1 & \Z \\ 0 & 0 & 1\end{smallmatrix} \right)$. Then $G$ is $2$-step nilpotent, and $G/\Gamma$ may be identified with the torus $(\R/\Z)^2$ via the map
\[ (x,y) \mapsto \left( \begin{smallmatrix} 1 & 0 & y \\ 0 & 1 & x \\ 0 & 0 & 1\end{smallmatrix} \right)\Gamma.\]
Taking $g := \left( \begin{smallmatrix} 1 & 1 & 0 \\ 0 & 1 & \alpha \\ 0 & 0 & 1\end{smallmatrix} \right)$, it is easy to check the action of $g$ on $G/\Gamma$ is given by $(x,y) \mapsto (x + \alpha, y + x)$.
The $3$-dimensional parallelepipeds of this nilflow take the form
$$ \big( x+(n+\omega \cdot  h)\alpha, y + \textstyle\frac{1}{2}\displaystyle(n + \omega \cdot  h)(n+\omega \cdot  h+1) \alpha
+ (n+\omega \cdot  h) x \big)_{\omega \in \{0,1\}^3}.$$
The key point to note here is that the first coordinate is at most linear in $n, h$, while the second coordinate is at most quadratic.  Take the set $\Sigma$ to be the set of all $7$-tuples $((x_\omega,y_\omega))_{\omega \in \{0,1\}^3_*}$ with the linear constraints
\begin{align*}
x_{000} + x_{011} &= x_{010} + x_{001} \\
x_{000} + x_{101} &= x_{001} + x_{100} \\
x_{000} + x_{110} &= x_{100} + x_{010}.
\end{align*}
The map $P: \Sigma \to (\R/\Z)^2$ is given by the alternating sum
$$ P( ((x_\omega,y_\omega))_{\omega \in \{0,1\}^3_*} ) := 
- \sum_{\omega \in \{0,1\}^3_*} (-1)^{|\omega|} (x_\omega,y_\omega).$$
This is ultimately a reflection of the fact that linear and quadratic functions have vanishing third derivative.
Note, in contrast to the previous example, that for the skew shift a vertex of a $2$-dimensional parallelepiped is \emph{not} determined continuously by the other three vertices.
\end{example}

Now we return to the task of proving Proposition \ref{decomp-prop}.  Let $P$ and $\Sigma$ be as in Proposition \ref{cube}.  
The function $ x \mapsto F(P( x))$ is continuous on the compact metric space $\Sigma$.  By the Stone-Weierstrass theorem, 
we may approximate this function to uniform accuracy $O(\eps)$ by a finite linear combination of tensor 
products of bounded Lipschitz functions on $G/\Gamma$, obtaining the uniform approximation
$$ F( P(  x ) )
= \sum_{\alpha \in A}  \prod_{\omega \in \{0,1\}^{s+1}_*} H_{\omega,\alpha}(x_\omega) + O(\eps)$$
for some finite index set $A$ and some $1$-bounded Lipschitz functions $H_{\omega,\alpha}: G/\Gamma \to [-1,1]$.  In particular, since $(g^{n + \omega \cdot h})_{\omega \in \{0,1\}^{s+1}_*}$ lies in $\Sigma$ and the image of this point under $P$ is $g^n x$, we have
$$ 
 F( g^n x ) = \sum_{\alpha \in A} \prod_{\omega \in \{0,1\}^{s+1}_*}
 H_{\omega,\alpha}( g^{n + \omega \cdot  h} x ) + O(\eps)$$
for all $g \in G$, $x \in G/\Gamma$, $n \in \Z$, and $ h \in \Z$.
  
Now we introduce the parameter $N \geq 1$ and average\footnote{One could take a limit here as $N \to \infty$, using an ergodic theorem to ensure suitable convergence; this would make the decomposition $F = F_1 + F_2$ independent of $N$, but at the cost of replacing the finite averaging in the definition of an averaged nilsequence with an infinite one.  We omit the details.}  the $ h$ parameter over the box $[N]^{s+1}$. In fact it is necessary to perform this averaging somewhat smoothly, to which end we take a smooth function cutoff $\sigma : \R \rightarrow [0,1]$ which is supported on $[-1,2]$ and equals $1$ on $[0,1]$, and then set
$$ 
 F( g^n x ) = F_1(n) + F_2(n)$$
 where
$$ F_1(n) := \sum_{\alpha \in A} \E_{ h \in [N]^{s+1}} 
 \sigma(h_1/N) \dots \sigma(h_{s+1}/N)\prod_{\omega \in \{0,1\}^{s+1}_*}
 H_{\omega,\alpha}( g^{n + \omega \cdot  h} x ) $$
and $F_2(n) = O(\eps)$.  In particular, we have $\Vert F_1 \Vert_{\infty} \leq 1 + O(\eps)$ since $F$ is bounded by $1$.  By shrinking the
Lipschitz functions $H_{\omega,\alpha}$ by a multiplicative factor of $1-O(\eps)$, and transferring the error over to $F_2$,
we may in fact ensure that $\Vert F_1 \Vert_{\infty} \leq 1$.
  
Now observe that for each fixed $\alpha,\omega$ and $ h$ the function given by
$n \mapsto H_{\omega,\alpha}( g^{n + \omega \cdot  h} x ) = H_{\omega,\alpha}( g^{n} (g^{\omega \cdot {h}} x) )$ is a Lipschitz nilsequence
on the $s$-step nilmanifold $G/\Gamma$, with Lipschitz constant independent of $N$, $g$ and $x$. We remarked, in \S \ref{gowers-sec}, that the Lipschitz nilsequences form an algebra in a certain sense. From this remark we conclude that $F_1$ is an averaged Lipschitz
$s$-step nilsequence on the product space $(G/\Gamma)^{\{0,1\}^{s+1}_*}$,
again with Lipschitz constant independent of $N,g$ and $x$.  To conclude the proof it suffices to show that $F_1$ is also
bounded in $U^{s+1}[N]^*$ uniformly in $N,g$ and $x$.  By the triangle inequality and the definition of the $U^{s+1}[N]^*$ norm,
it thus suffices to show that the absolute value of
\begin{equation}\label{eq8.11} \E_{n \in [N];  h \in [N]^{s+1}} f(n) \sigma(h_1/N) \dots \sigma(h_{s+1}/N)\prod_{\omega \in \{0,1\}^{s+1}_*}
 H_{\omega,\alpha}( g^{n + \omega \cdot  h} x )
\end{equation}
is uniformly bounded in $N,g,x$ whenever $f: [N] \to \R$ satisfies $\Vert f \Vert_{U^{s+1}[N]} \leq 1$. From this point onwards we do not care what the functions $n \mapsto H_{\omega,\alpha}(g^n x)$ actually are: it is merely important that they are $1$-bounded. For that reason we write $\b_{\omega}(n) = H_{\omega,\alpha}(g^{n} x)$, whereupon the quantity \eqref{eq8.11} that we are to show is uniformly bounded becomes 
\begin{equation}\label{eq8.12}
\E_{n \in [N];  h \in [N]^{s+1}} f(n) \sigma(h_1/N) \dots \sigma(h_{s+1}/N)\prod_{\omega \in \{0,1\}^{s+1}_*} \b_{\omega}(n + \omega \cdot  h).
\end{equation}

At this point we transfer to a group $\Z_{N'}$ where $N' = 10sN$ (say). Slightly abusing notation, the expression \eqref{eq8.12} is, up to factors of $O_s(1)$, equal to 
\begin{equation}\label{eq8.13}
\E_{n \in \Z_{N'}; {h} \in \Z_{N'}^{s+1}} f(n) \sigma(h_1/N) \dots \sigma(h_{s+1}/N)\prod_{\omega \in \{0,1\}^{s+1}_*} \b_{\omega}(n + \omega \cdot {h}).
\end{equation}
Here we have extended $f$ from $[N]$ to all of $\Z_{N'}$ by defining it to be zero outside of $[N]$. Now by taking a Fourier expansion on $\Z_{N'}^{s+1}$ we may write
\[ \sigma(h_1/N) \dots \sigma(h_{s+1}/N) = \sum_{r_1,\dots,r_{s+1}} c_{r_1,\dots, r_{s+1}} e\big( (r_1 h_1 + \dots + r_{s+1} h_{s+1})/N\big).\]
By choosing the cutoff $\sigma$ to be sufficiently smooth, we may ensure that
\[ \sum_{r_1,\dots,r_{s+1}} |c_{r_1,\dots, r_{s+1}}|  = O_{s}(1).\]
Thus to show that \eqref{eq8.13} is uniformly bounded it suffices to show the same for 
\begin{equation}\label{eq8.14a}
\E_{n \in \Z_{N'}; {h} \in \Z_{N'}^{s+1}} f(n) e \big( (r_1 h_1 + \dots + r_{s+1} h_{s+1})/N\big)\prod_{\omega \in \{0,1\}^{s+1}_*} \b_{\omega}(n + \omega \cdot {h})
\end{equation}
for all $r_1,\dots,r_{s+1} \in \Z_{N'}$. It is easy to see that the exponential may be split up and incorporated into the $\b_{\omega}( \; )$ terms, and therefore we have reduced the matter to placing a bound on
\begin{equation}\label{eq8.14}
\E_{n \in \Z_{N'}; {h} \in \Z_{N'}^{s+1}} f(n)\prod_{\omega \in \{0,1\}^{s+1}_*} \b_{\omega}(n + \omega \cdot {h}).
\end{equation}
Now we are assuming that $\Vert f \Vert_{U^{s+1}[N]} \leq 1$. By Lemma \ref{normlemma} this implies that $\Vert f \Vert_{U^{s+1}(\Z_{N'})}$ $= O_s(1)$. The boundedness now follows from the Gowers-Cauchy-Schwarz inequality \eqref{gcz-u}. Tracing backwards, we see in turn that \eqref{eq8.14}, \eqref{eq8.13}, \eqref{eq8.12} and \eqref{eq8.11} are all $O_s(1)$, thereby concluding the proof.
\end{proof}

Although we will not need this fact here, it is interesting to note that Proposition \ref{decomp-prop} allows one to
extend Proposition \ref{nil-gow} from bounded $f$ to integrable $f$:

\begin{corollary}[Nilsequences obstruct uniformity, II]\label{nil-gow-cor} Let $s \geq 0$ and $\delta \in (0,1)$. Let $G/\Gamma = (G/\Gamma, d_{G/\Gamma})$ be a nilmanifold with some fixed smooth metric $d_{G/\Gamma}$, and let $(F(g^n x))_{n \in \N}$ be a bounded $s$-step nilsequence with Lipschitz constant at most $M$.  Let $f : [N] \rightarrow \R$ be a function for which
$$ \E_{n \in [N]} |f(n)| \leq 1 $$
and
\[ |\E_{n \in [N]} f(n) F(g^n x)| \geq \delta.\]
Then we have
\[ \Vert f \Vert_{U^{s+1}[N]} \gg_{s,\delta,M,G/\Gamma} 1.\]
\end{corollary}

\begin{proof} We apply Proposition \ref{decomp-prop} with $\eps$ equal to a small multiple of $\delta$, and conclude from the triangle inequality that
\[ |\E_{n \in [N]} f(n) F_1(n)| \geq \delta/2.\]
Since $F_1$ has a $U^{s+1}[N]^*$ norm of $O_{s,\delta,M,G/\Gamma}(1)$, the claim follows.
\end{proof}

\section{A splitting of the von Mangoldt function}\label{wtrick-sec}

To summarise so far, we have reduced the task of proving that the $\GI(s)$ and $\MN(s)$ conjectures imply the Main Theorem to 
the much easier task of establishing Proposition \ref{manortho-avg}.  This is a correlation estimate involving $\Lambda'_{b,W}$.  It is convenient to return at this point to the original von Mangoldt function $\Lambda$.  The contribution from the prime powers which are introduced when $\Lambda'_{b,W}$ is replaced by $\Lambda_{b,W}$ is easily seen to be negligible, and so it suffices to establish the estimate
\[ \E_{n \in [N]} (\Lambda_{b,W}(n)-1) F_1(n) = o_{M,M',G/\Gamma,s}(1).\]
Recalling the definition \eqref{lbn} of $\Lambda_{b,W}(n)$, we are thus trying to establish the bound
\begin{equation}\label{lb}
\E_{n \in [N]} (\frac{\phi(W)}{W} \Lambda(Wn+b)-1) F_1(n) = o_{M,M',G/\Gamma,s}(1).
\end{equation}

At this point we perform a standard decomposition of $\Lambda$ into a ``smooth'' piece $\Lambda^\sharp$ corresponding to small divisors and a ``rough'' piece $\Lambda^\flat$ corresponding to large divisors.  We take a small exponent $\gamma = \gamma_s > 0$, whose exact value will be specified later, and set $R := N^\gamma$.  Observe from \eqref{lambda-mu} that
$$ \Lambda(n) = - \log R \sum_{d|n} \mu(d) \chi(\frac{\log d}{\log R})$$
where $\chi: \R^+ \to \R^+$ is the identity function $\chi(x) := x$.  We now perform a smooth splitting 
$\chi = \chi^\sharp + \chi^\flat$, where $\chi^\sharp(x)$ vanishes for $|x| \geq 1$ and $\chi^\flat(x)$ vanishes
for $|x| \leq 1/2$, the precise form of this splitting being unimportant.  
This induces a splitting $\Lambda = \Lambda^\sharp + \Lambda^\flat$, where
\begin{equation}\label{loof}
\Lambda^\sharp(n) := - \log R \sum_{d|n} \mu(d) \chi^\sharp(\frac{\log d}{\log R}) \quad
 \mbox{and} \quad \Lambda^\flat(n) := - \log R \sum_{d|n} \mu(d) \chi^\flat(\frac{\log d}{\log R}).
\end{equation}
Thus to prove \eqref{lb} it will suffice to show the estimates
\begin{equation}\label{lb-sharp}
\E_{n \in [N]} (\frac{\phi(W)}{W} \Lambda^\sharp(Wn+b)-1) F_1(n)=  o_{s,M'}(1)
\end{equation}
and
\begin{equation}\label{lb-flat}
\E_{n \in [N]} \frac{\phi(W)}{W} \Lambda^\flat(Wn+b) F_1(n) = o_{M,G/\Gamma,s}(1).
\end{equation}
We begin by establishing the bound \eqref{lb-sharp}.  It is here that we need the dual norm bound \eqref{f1un}.  Indeed, from that bound we have
\begin{align*}
\big|\E_{n \in [N]} (\frac{\phi(W)}{W} \Lambda^\sharp(Wn+b)-1) F_1(n)\big| &\leq \big\| \frac{\phi(W)}{W} \Lambda^\sharp(Wn+b)-1 \big\|_{U^{s+1}[N]} \Vert F_1 \Vert_{U^{s+1}[N]*} \\ &
\leq M' \big\|\frac{\phi(W)}{W} \Lambda^\sharp(Wn+b)-1 \big\|_{U^{s+1}[N]}.
 \end{align*}
 It suffices, then to show that
 \begin{equation}\label{to-prove-appD} \big\|\frac{\phi(W)}{W} \Lambda^\sharp(Wn+b)-1 \big\|_{U^{s+1}[N]} = o_s(1).\end{equation}
This is a multilinear correlation estimate for a truncated divisor sum, and can be treated by standard sieve theory methods related to the correlation estimates of Goldston and Y{\i}ld{\i}r{\i}m \cite{gy-1,gy-2,gy-3} provided that the exponent $\gamma$ is sufficiently small (an appropriate choice would be, for example, $\gamma_s := \frac{1}{10}2^{-s}$). We provide the details of this computation in Appendix \ref{gy-sec}. This establishes
\eqref{lb-sharp}.

It remains to establish the bound \eqref{lb-flat}.  Recall that $F_1$ is an averaged nilsequence.  From the triangle inequality, it will thus suffice to prove the bound
\begin{equation}\label{eq8.177}
\E_{n \in [N]} \frac{\phi(W)}{W} \Lambda^\flat(Wn+b) F(g^n x) = o_{M,G/\Gamma,s}(1)
\end{equation}
for all $1$-bounded $s$-step nilsequences $F(g^n x)$ of Lipschitz constant $M$. We emphasise that the $o$-term is required to depend only on $M,G/\Gamma$ and $s$, and should be otherwise be independent of $F, g$ and $x$.  

We will eventually apply the $\MN(s)$ conjecture, which comes with the safety net of an error term which decays like $\log^{-A} N$ for any $A$. With this in mind, we begin by removing the $W$-dependence in \eqref{eq8.177} in a rather crude fashion. Since $\phi(W)/W \leq 1$, we ignore this factor completely. 

Now by a simple substitution we have
\begin{equation}\label{eq8.775}
\E_{n \in [N]} \Lambda^{\flat}(Wn + b) F_1(g^n x) = W \E_{b < n \leq WN + b} 1_{n \equiv b \mdsub{W}} \Lambda^{\flat}(n) F_1(g^{(n-b)/W} x).
\end{equation}
Now any Lie group $G$ over $\R$ for which the exponential map $\exp : \mathfrak{g} \rightarrow G$ from the associated Lie algebra is surjective is \emph{divisible}, meaning that given any $g \in G$ and any positive integer $m$ there is an element $g^{1/m} \in G$ with $(g^{1/m})^m = g$. When $G$ is simply-connected and nilpotent, $\exp$ is a homeomorphism (see \cite{bourbaki} for details). In our setting, write $g' := g^{1/W}$ and $x' := g^{-b/W}x$. Then for all $n \equiv b \md{W}$ we have 
\begin{equation}\label{eq8.776} F_1(g^{\prime n} x') = F_1(g^{(n-b)/W}x).\end{equation} Note that the left-hand side here makes perfect sense for \emph{any} $n$, not just for $n$ such that $n \equiv b \md{W}$.

The constraint $1_{n \equiv b \mdsub{W}}$ may be expanded as a Fourier series 
\[ 1_{n \equiv b \mdsub{W}} = \frac{1}{W}\sum_{r \in \Z_W} e(-rb/W) e(rn/W)\]
on $\Z_W$. We substitute this and \eqref{eq8.776} into \eqref{eq8.775}, noting that each function $n \mapsto e(rn/W)$ may be realised as a $1$-bounded, $O(1)$-Lipschitz nilsequence on the $1$-step nilmanifold $\R/\Z$. Replacing $G/\Gamma$ with $G/\Gamma \times \R/\Z$, we see that in order to prove \eqref{eq8.177} it suffices to show that
\begin{equation}\label{eq8.188} W\E_{b < n \leq WN+ b} \Lambda^{\flat}(n) F(g^n x) = o_{M,G/\Gamma,s}(1)\end{equation} for all $M$-Lipschitz $1$-bounded nilsequences $(F(g^n x))_{n \in N}$ on an $s$-step nilmanifold $G/\Gamma$.  

In fact we will establish the stronger estimate
\begin{equation}\label{eq8.189}
|\sum_{n \in [N]} \Lambda^\flat(n) F(g^n x)| \ll_{M,G/\Gamma,s,A} N \log^{-A} N\end{equation}for any $A > 0$. Note that $w$ was chosen to be so slowly growing that $W = O(\log N)$, so this estimate really is stronger than \eqref{eq8.188}. We expand the left-hand side of \eqref{eq8.189} using \eqref{loof} and reduce to showing that
\begin{equation}\label{eq8.193} |\sum_{n \in [N]} \sum_{d|n} \mu(d) \chi^\flat(\frac{\log d}{\log R}) F(g^n x)| \ll_{M,G/\Gamma,s,A} N \log^{-A} N.\end{equation}
The left-hand side may be rearranged as
$$ |\sum_{m \in [N]} \sum_{d \in [N/m]}  \mu(d) \chi^\flat(\frac{\log d}{\log R}) F((g^m)^d  x)|.$$
Observe that $\chi^\flat$ is supported on $|x| \geq 1/2$, and so the summand vanishes unless $d \geq R^{1/2}$, in which case $m \leq N/R^{1/2}$.  We now apply the \mobname conjecture $\MN(s)$. Together with a straightforward summation by parts to remove the smooth cutoff $\chi^{\flat}$ this shows that
$$ |\sum_{d \in [N/m]}  \mu(d) \chi^\flat(\frac{\log d}{\log R}) F((g^m)^d  x)| \ll_{M,G/\Gamma,s,A} \frac{N}{m} \log^{-A} \frac{N}{m}.$$
Note that we are making critical use here of the fact that the bounds in the $\MN(s)$ conjecture are uniform in the $g$ parameter in order to deal with the fact that we have dilated $g$ to $g^m$.  Since $m \leq N/R^{1/2}$, we see that $\log^{-A}(N/m) \ll_A \log^{-A} N$.  Summing in $m$ and absorbing the logarithmically divergent sum $\sum_{m \in [N]} \frac{1}{m}$ into the $\log^{-A} N$ factor we obtain \eqref{eq8.193} as desired. This in turn implies \eqref{eq8.189} and hence, by our earlier series of reductions, \eqref{lb-flat}. Together with \eqref{lb-sharp}, which we have already established, this concludes the proof of Proposition \ref{manortho-avg}. By our long series of earlier reductions, this (finally!) completes the proof of the Main Theorem.\endproof

\section{Variations on the main argument and other remarks}\label{remarks-sec}

It is conceivable that our methods here extend to certain ``finite complexity'' multilinear averages involving systems of \emph{polynomials} $\psi_j( n)$ rather than affine-linear forms.  Indeed, the machinery of ``PET induction'' (see e.g. \cite{pet}) allows us in principle to use repeated applications of Cauchy-Schwarz to control certain of these averages by Gowers uniformity norms. A model problem would be to count the number of $p,n$ for which the numbers $p, p+n, p+n^2, \ldots, p+n^k$ are all prime.  A na\"{\i}ve attempt to do this meets with what seems to be an insurmountable obstacle. Namely, in order to restrict the range of the primes concerned to an interval such as $[N]$, certain other parameters (for example the ``shifts'' ${h}$ in the definition of the Gowers norms) have to be restricted to a much smaller range, say of size $O(N^{1/100})$.  This makes it impossible to pass back and forth between $[N]$ and $\Z_{N'}$ as we have done above, and the evaluation of exponential sums with $\mu$ or $\Lambda$ on such a range seems to be beyond hope, even assuming the GRH. It may be that the PET induction scheme can be ``globalised'' to avoid these issues, but we do not know how to address this at present.

For the benefit of readers who are only interested in the unconditional ``quadratic'' ($s=2$) applications of this paper such as Corollary \ref{maincor} or Examples \ref{example1}-\ref{example3} we outline
a shorter path to the Main Theorem in that case. This approach avoids Lie theory completely, and probably represents the best approach to obtaining bounds for error terms. Note, however, that with either approach our error terms are completely ineffective unless the GRH is assumed. The introduction of Lie theory, though strictly speaking unnecessary, seems to make our work easier to understand from the conceptual point of view. This is especially the case when $s \geq 3$, where it is not even clear how Lie theory-free analogues of the $\GI(s)$ and $\MN(s)$ conjecture might be formulated.

In the quadratic case it is possible to replace the concept of a $2$-step nilsequence by more concrete objects. In a sense these are more basic than $2$-step nilsequences, if only because in \cite{green-tao-u3inverse} we introduce these objects first and then build nilsequences from them. Note, however, that this may be an artifact of our approach.

These more basic objects can then be manipulated by hand without resorting to machinery such as the Host-Kra theory 
in Appendix \ref{nil-app}. Let us consider, by way of illustration, the following more concrete version of the inverse Gowers-norm conjecture $\GI(2)$ which was proven in   \cite{green-tao-u3inverse}.

\begin{theorem}[$U^3$ inverse theorem with bracket polynomials]\label{bracketpoly} Let $f: [N] \to [-1,1]$ be such that 
$\|f\|_{U^3[N]} > \delta$ for some $0 < \delta \leq 1$ and $N \geq 1$.
Then there exists a positive 
integer $J = O_\delta(1)$ and real numbers $a_j, b_j, \xi_{j,1}, \xi_{j,2}, \xi_{j,3}$ for $j \in [J]$ such that
\begin{equation}\label{piece}
| \E_{n \in [N]} f(n) e( \phi(n) ) | \gg_\delta 1
\end{equation}
where $\phi$ is the function
$$ \phi(n) := - \sum_{j \in [J]}\left( a_j \{ \xi_{j,1} n \} \{ \xi_{j,2} n \} + b_j \{ \xi_{j,3} n \}\right).$$
\end{theorem}
 
\begin{remark} As before, $\{x\}$ denotes the fractional part of $x$, which we take to lie in $(-\frac{1}{2},\frac{1}{2}]$.\end{remark}

This result follows quickly from \cite[Theorem 10.9]{green-tao-u3inverse} using Lemma \ref{normlemma} to work in a cyclic group of prime order. We refer to the phase $\phi(n)$ \eqref{piece} as a ``bracket polynomial''. 
By modifying the arguments in \S \ref{sec5}, one can transfer this theorem to the case when $f$ is bounded by a pseudorandom measure $\nu$ rather than by $1$, thereby reducing Theorem \ref{gowers-norm} to the establishment of the exponential sum estimate
$$  \E_{n \in [N]} (\Lambda_{b,W}(n)-1) e\big( - \sum_{j \in [J]}\big( a_j \{ \xi_{j,1} n \} \{ \xi_{j,2} n \} + b_j \{ \xi_{j,3} n \}\big)\big)  = o_{J}(1)$$
uniformly over all $b \in [W]$ with $\gcd(b,W)=1$.  This could in principle\footnote{Indeed, this exponential sum is a more complicated variant of the more traditional exponential sum $\sum_{n \in [N]} \Lambda(n) e( \alpha n^2 )$, which was considered for instance in \cite{ghosh,hua} .} be established directly by Vinogradov's method, following
the machinery in \cite{green-tao-u3mobius}, though the argument would be rather lengthy.  Alternatively one can deduce this result from the corresponding results for the M\"obius function established in \cite{green-tao-u3mobius} using a variant of the arguments in this paper.

A key difference is that the Host-Kra machinery and the machinery of averaged nilsequences are no longer required. Instead, 
the above function $e(\phi(n))$ can be replaced by a smoother variant, constructed for instance using a variant of the dual function machinery in
\cite{green-tao-longprimeaps}, in order to obtain a function which is bounded in $(U^3)^*$. This provides an analogue of Proposition \ref{manortho-avg}, and from that point onwards one may proceed similarly.

One could also use a still more ``basic'' type of obstruction for the $U^3$-norm, namely phases which are locally quadratic on Bohr sets (cf. \cite[\S 2]{green-tao-u3inverse}). These require even less unpacking than the bracket quadratics above, and indeed it was found to be rather convenient to work with these functions in \cite{green-tao-u3mobius}. It takes a while to even define these functions properly, however, and they suffer from a few technical deficiencies which affect various other steps of the argument. Perhaps the most serious is that if $n \mapsto f(n)$ is such a function then $n \mapsto f(dn)$ need not quite be, a phenomenon which causes trouble in \S \ref{wtrick-sec}.\endproof

\section{A brief discussion of bounds}\label{bounds-sec}

We have shied away from giving any explicit bounds on our $o(1)$ error terms. There are at least two reasons for this. Firstly, it is notationally easier to avoid doing so. Secondly, and much more importantly, unless one assumes the GRH we do not have any explicit bounds!

By way of illustration, let us consider the statement
\begin{equation}\label{mobius-4-aps} \E_{x,d\leq N} \mu(x)\mu(x+d)\mu(x+2d)\mu(x+3d) = o(1),\end{equation}
which follows from the case $s=2$ of Proposition \ref{mu-lam-cor}. A discussion of correlations involving $\Lambda$ would go along similar lines, but there is the distraction of the singular product $\beta_{\infty} \prod_p \beta_p$.

As we remarked, the error term here is completely ineffective without assuming GRH. Indeed to show that the left-hand side in \eqref{mobius-4-aps} is at most $\delta$, we would ultimately (deep inside the paper \cite{green-tao-u3mobius}) need estimates for the sum of the M\"obius function over arithmetic progressions with common difference $q \sim \log^{A(\delta)}N$. Although such estimates exist, the error terms involve an ineffective constant $C(A(\delta))$ due to the possible presence of Landau-Siegel zeros.

Assuming the GRH one could prove using our methods that 
\[ |\E_{x,d \leq N} \mu(x)\mu(x+d)\mu(x+2d)\mu(x+3d)| \leq C\log^{-c} N\] for some explicit $C$ and some explicit (but small) $c > 0$. To obtain such a result it would be best to avoid the use of Lie theory as outlined in \S \ref{remarks-sec}, since the many approximation arguments involved in that theory are quite costly from the quantitative point of view.

Improved results in additive combinatorics (particularly a solution to the so-called Polynomial Freiman-Ruzsa conjecture, which could be used as an input in \cite{green-tao-u3inverse}) could lead to a bound of the shape $\exp(-\log^c N)$. However it seems that obtaining a bound $N^{-c}$ is very difficult. 

Unconditionally, a bound in \eqref{mobius-4-aps} of the form $O(f(n))$ for some explicit function $f(n)$ tending to zero as $n \to \infty$ and some ineffective implied constant $O( )$ would be very interesting.

To set the above discussion in context, we mention the best available results for three-term progressions, which follow from estimates for $\sup_{\alpha \in \R/\Z}|\E_{n \leq N} \mu(n) e(\alpha n)|$. These seem to be as follows.

\[ \E_{x,d \leq N} \mu(x) \mu(x+d) \mu(x + 2d) \ll \left\{\begin{array}{lll} C_A\log^{-A}N & \mbox{any $A > 0$} & \mbox{Davenport \cite{davenport-old}} \\ C_\epsilon N^{-1/4 + \epsilon} & \mbox{on GRH} & \mbox{Baker--Harman \cite{baker-harman}.}  \end{array} \right. \]

Bounds of a similar type could be obtained for any instance of Proposition \ref{mu-lam-cor} with $s = 1$.

\appendix

\section{Elementary convex geometry}\label{convexgeom}

In this appendix we recall some profoundly classical facts concerning convex bodies which will allow us to manipulate
cutoffs such as $1_K$ readily, beginning with an ancient observation of Archimedes.  

\begin{lemma}[Archimedes comparison principle]  Let $K_1 \subseteq K_2 \subseteq \R^d$ be bounded convex bodies.  Then the surface area of $K_1$ is less than or equal to the surface area of $K_2$.
\end{lemma}

\begin{proof} It is easy to see that the intersection of $K_2$ with a half-space has lesser or equal surface area than $K_2$.  Since $K_1$ can be approximated to arbitrary accuracy by the intersection of finitely many half-spaces, the claim follows.
\end{proof}

\begin{corollary}[Boundary region estimate]\label{boundary-volume}  Let $K \subseteq [-N,N]^d$ be a convex body.  If $\eps \in (0,1)$, then the $\eps N$-neighbourhood of the boundary $\partial K$ has volume $O_d(\eps N^d)$.
\end{corollary}

\begin{proof} Rescale so that $N=1$. By differentiating in $\eps$ we see that it suffices to show that any convex body in $[-2,2]^d$ has surface area $O_d(1)$.  But this follows from the Archimedes comparison principle. One could also derive this fact using the theory of mixed volumes; see \cite{steiner}.
\end{proof}

At this point we can now readily prove \eqref{beta-gauss} using the Gauss volume-packing argument.  By intersecting $K$ with the half-spaces $\{ x \in \R^d: \psi_j(x) > 0 \}$ it suffices to show that
$$ 
|K \cap \Z^d| = \vol_d(K) + O_d(N^{d-1})$$
for all convex bodies $K \subseteq [-N,N]$.  However, given that $|K \cap \Z^d|$ is equal to the volume of the set
$(K \cap \Z^d) + [-1/2,1/2]^d$, which differs from $K$ only on the $O_d(1)$-neighbourhood of $\partial K$, the claim then follows from Corollary \ref{boundary-volume}.

Now we give an analytic consequence of Corollary \ref{boundary-volume}.

\begin{corollary}[Lipschitz approximation of convex indicators]\label{lip-approx}  Let $K \subseteq [-N,N]^d$ be a convex body and let $\eps \in (0,1)$.  Then we can write $1_K = F_\eps + O(G_\eps)$, where $F_\eps, G_\eps$ are non-negative Lipschitz functions on $[-2N,2N]^d$ with Lipschitz constants $O( \frac{1}{\eps N} )$ and bounded in magnitude by $1$, and where $\int_{\R^d} G_\eps(x)\ dx = O_d( \eps N^d )$.
\end{corollary}

\begin{proof}  We take
$$ F_\eps(x) := \max( 1 - \frac{\dist_{\R^d}( x, K )}{\eps N}, 0 ) \quad \mbox{and} \quad G_\eps(x) := 
\max( 1 - \frac{\dist_{\R^d}( x, \partial K )}{2\eps N}, 0 ).$$
The claim follows easily from Corollary \ref{boundary-volume}.
\end{proof}

In practice, Corollary \ref{lip-approx} allows us to replace a rough cutoff such as $1_K$ with the smoother operation of Lipschitz
cutoffs.  This can then be combined with Fourier analysis to replace the Lipschitz cutoffs in turn with modulations by linear phases, which turn out to be utterly harmless in our analysis.  This might remind readers of the P\'olya-Vinogradov completion-of-sums method, or the Erd\H{o}s-Tur\'an inequality.

\section{Gowers norm theory}\label{gowersnorm-sec}

In this appendix we develop the general ``elementary'' theory of Gowers uniformity norms, which were introduced in
\cite{gowers-long-aps} and subsequently, in the rather different context of ergodic theory, in \cite{host-kra}.  By elementary in this context, we basically mean
that we only pursue here those results which can be obtained as an easy consequence of the Cauchy-Schwarz inequality. This is
in contrast to the more advanced inverse theory involving nilsequences, Fourier analysis, and suchlike.  The theory here is an amalgam of parts of \cite[\S 3]{gowers-long-aps}, \cite{green-icm}, \cite[\S 5]{green-tao-longprimeaps},  \cite[\S 1]{green-tao-u3inverse},  \cite{host-kra}, \cite[\S 3]{tao-elescorial}, \cite{tao-hyper,tao-multiprime}, or \cite[Ch. 11]{tao-vu}.

It is convenient to work rather abstractly at first, dealing with complex-valued functions of many variables. This level of abstraction will be useful for us when we prove the generalised von Neumann theorem, Proposition \ref{gvn}, in \S \ref{gvn-app}. The argument is essentially that of \cite[\S 5]{green-tao-longprimeaps}, generalised to handle arbitrary systems of linear forms rather than merely $k$-term APs, but the introduction of extra notation somewhat eases the process of actually carrying this out.

\begin{definition}[Gowers box norms]  Let $(X_\alpha)_{\alpha \in A}$ be a finite non-empty collection of finite non-empty sets, and for any $B \subseteq A$ write $X_B := \prod_{\alpha \in B} X_\alpha$ for the Cartesian product.  If $f: X_A \to \C$ is a complex-valued function, we define the \emph{Gowers box norm} $\|f\|_{\Box(X_A)} \in \R^+$ to be the unique non-negative real number such that
\begin{equation}\label{fa}
 \|f\|_{\Box(X_A)}^{2^{|A|}} := \E_{x^{(0)}_A, x^{(1)}_A \in X_A} \prod_{\omega_A \in \{0,1\}^A} {\mathcal C}^{|\omega_A|} f( x^{(\omega_A)}_A ) 
 \end{equation}
where ${\mathcal C}: z \mapsto \overline{z}$ is complex conjugation, and for any $x^{(0)}_A = (x^{(0)}_\alpha)_{\alpha \in A}$ and $x^{(1)}_A = (x^{(1)}_\alpha)_{\alpha \in A}$ in $X_A$ and
$\omega_A = (\omega_\alpha)_{\alpha \in A}$ in $\{0,1\}^A$, we write $x^{(\omega)}_A := (x^{(\omega_\alpha)}_\alpha)_{\alpha \in A}$
and $|\omega_A| := \sum_{\alpha \in A} \omega_\alpha$.  We adopt the convention that if $A$ is empty (so that $f$ is a constant), then $\|f\|_{\Box(X_A)} := f$.
\end{definition}

It is not immediately obvious that the right-hand side of \eqref{fa} is non-negative, or that the term ``norm'' is appropriate. We will establish both of these facts below.

\begin{examples}  If $A = \{1\}$, then
$$ \| f \|_{\Box(X_1)} = \big(\E_{x^{(0)}_1,x^{(1)}_1 \in X_1} f(x^{(0)}_1) f(x^{(1)}_1)\big)^{1/2} = |\E_{x_1 \in X_1} f(x_1)|$$
while if $A = \{1,2\}$, then $\| f \|_{\Box(X_{1,2})} =$
\[ \bigg(\E_{x^{(0)}_1, x^{(1)}_1 \in X_1; x^{(0)}_2, x^{(1)}_2 \in X_2}
f(x^{(0)}_1, x^{(0)}_2) \overline{f(x^{(0)}_1, x^{(1)}_2)} \overline{f(x^{(1)}_1, x^{(0)}_2)} 
f(x^{(1)}_1, x^{(1)}_2)\bigg)^{1/4}.
\]
In general, the $2^{|A|}$th power of the $\Box(X_A)$ norm on $f_A$ is a multilinear average of $f_A$ over $|A|$-dimensional boxes (hence the name).
\end{examples}

It is easy to verify the recursive relationship
\begin{equation}\label{recursive}
 \| f \|_{\Box(X_A)}^{2^{|A|}} = \E_{x^{(0)}_\alpha, x^{(1)}_\alpha \in X_\alpha} 
\|f(\cdot, x^{(0)}_\alpha) \overline{f(\cdot, x^{(1)}_\alpha)} \|_{\Box(X_{A \backslash \{\alpha\}})}^{2^{|A|-1}}
\end{equation}
whenever $\alpha \in A$, which can be used as an alternate definition of the box norms.  In particular we 
see that the box norms $\| f\|_{\Box(X_A)}$ are non-negative for $A$ non-empty.  These norms are also conjugation-invariant,
homogeneous, and enjoy the positivity property
\begin{equation}\label{pointbound}
\| f \|_{\Box(X_A)} \leq \| \nu \|_{\Box(X_A)}
\end{equation}
whenever $f: X_A \to \C$ and $\nu: X_A \to \R^+$ obey the pointwise bound $|f(x_A)| \leq \nu(x_A)$ for all $x_A \in X_A$.

The box norms are also invariant under a large class of phase modulations.  Indeed one easily verifies from \eqref{recursive} and
induction that
\begin{equation}\label{phase-modulation}
\| f e( \sum_{B \subsetneq A} \phi_B ) \|_{\Box(X_A)} = \| f\|_{\Box(X_A)} 
\end{equation}
where $e: \R/\Z \to \C$ is the standard character $e(x) := e^{2\pi i x}$ and for each proper subset $B \subseteq A$, the phase function $\phi_B: X_B \to \R/\Z$ is arbitrary.  Thus the $\Box(X_A)$ norm is insensitive to ``lower order'' modulations which involve only a proper subset of the variables in $X_A$.

A fundamental inequality\footnote{In our treatment here, this inequality plays a more central role than in earlier papers; we are using it as a kind of ``universal Cauchy-Schwarz inequality'', in the sense that any other inequality that we need, which would in earlier papers be proven by multiple applications of the ordinary Cauchy-Schwarz inequality, is instead proven here by a single application of the Gowers-Cauchy-Schwarz inequality.  This seems to fit with the philosophy that the Gowers norms are somehow ``universal'' or ``characteristic'' for all averages of a certain complexity.}
 concerning these norms is

\begin{lemma}[Gowers-Cauchy-Schwarz inequality]\label{gczow}  Let $(X_\alpha)_{\alpha \in A}$ be a finite collection of finite non-empty sets.  For every $\omega_A \in \{0,1\}^A$ let $f_{\omega_A}: X_A \to \C$ be a function.
Then
\begin{equation}\label{gcz-box}
\big|\E_{x_A^{(0)}, x_A^{(1)} \in X_A}
\prod_{\omega_A \in \{0,1\}^A} {\mathcal C}^{|\omega_A|} f_{\omega_A}(x_A^{(\omega_A)})\big|
\leq \prod_{\omega_A \in \{0,1\}^A} \|f_{\omega_A} \|_{\Box(X_A)}.
\end{equation}
\end{lemma}

\begin{proof} We induct on $|A|$.  When $|A|=0$ the claim trivially holds, and in fact there is equality.  Now suppose that $|A| \geq 1$ and the claim has already been proven for smaller sets $A$. 

Partition $A$ as $A' \cup \{\alpha\}$ for some $\alpha \in A$.  We can rewrite the left-hand side of \eqref{gcz-box} as
$$ |\E_{x^{(0)}_{A'}, x^{(1)}_{A'} \in X_{A'}} \prod_{\omega_\alpha \in \{0,1\}} {\mathcal C}^{\omega_\alpha} F_{\omega_\alpha}( x^{(0)}_{A'}, x^{(1)}_{A'} )|$$
where
$$ F_{\omega_\alpha}( x^{(0)}_{A'}, x^{(1)}_{A'} )
:=  \E_{x^{(\omega_\alpha)}_\alpha \in X_\alpha}
\prod_{\omega_{A'} \in \{0,1\}^{A'}} {\mathcal C}^{|\omega_{A'}|} f_{(\omega_{A'}, \omega_\alpha)}(x^{(\omega_{A'})}, x_\alpha).
$$
By Cauchy-Schwarz it thus suffices to show that
$$ \E_{x^{(0)}_{A'}, x^{(1)}_{A'} \in X_{A'}} |F_{\omega_\alpha}( x^{(0)}_{A'}, x^{(1)}_{A'} )|^2
\leq \prod_{\omega_{A'} \in \{0,1\}^{A'}} \|f_{(\omega_{A'}, \omega_\alpha)} \|_{\Box(X_A)}^2
$$ 
for each $\omega_\alpha \in \{0,1\}$.  We can expand the left-hand side as 
$$
\E_{x^{(0)}_\alpha, x^{(1)}_\alpha \in X_\alpha}
\E_{x^{(0)}_{A'}, x^{(1)}_{A'} \in X_{A'}}
\prod_{\omega_{A'} \in \{0,1\}^{A'}} {\mathcal C}^{|\omega_{A'}|} 
\bigg(f_{(\omega_{A'},\omega_\alpha)}(x^{(\omega_{A'})}_{A'},x^{(0)}_\alpha)
\overline{f_{(\omega_{A'},\omega_\alpha)}(x^{(\omega_{A'})}_{A'},x^{(1)}_\alpha)}\bigg).$$
Applying the induction hypothesis, we can bound this by
$$ 
\E_{x^{(0)}_\alpha, x^{(1)}_\alpha \in X_\alpha}
\prod_{\omega_{A'} \in \{0,1\}^{A'}} 
\| f_{(\omega_{A'},\omega_\alpha)}(\cdot,x^{(0)}_\alpha)
\overline{f_{(\omega_{A'},\omega_\alpha)}(\cdot,x^{(1)}_\alpha)} \|_{\Box(X_{A'})}$$
and the claim now follows from H\"older's inequality and \eqref{recursive}.
\end{proof}

From \eqref{gcz-box} we easily deduce the \emph{Gowers triangle inequality}
$$ \| f + g \|_{\Box(X_A)} \leq \| f \|_{\Box(X_A)} + \| g \|_{\Box(X_A)} $$
as can be seen by raising both sides to the power $2^{|A|}$. Let us also observe, setting all but one of the functions
in \eqref{gcz-box} to be Kronecker delta functions, that if $\|f\|_{\Box(X_A)} = 0$ and $|A| \geq 2$ then $f$ vanishes identically.
Thus we see that the $\Box(X_A)$-norm is indeed a norm for $|A| \geq 2$, whilst for $|A|=1$ it is merely a semi-norm.

As a consequence of the Gowers-Cauchy-Schwarz inequality we obtain

\begin{corollary}[Second Gowers-Cauchy-Schwarz inequality]\label{gczow-2}  Let $(X_\alpha)_{\alpha \in A}$ be a collection of finite non-empty sets.  For every $B \subseteq A$ let $f_B: X^{B} \to \C$ be a function.  Then
\begin{equation}\label{gczow2}
|\E_{x_A \in X_A} \prod_{B \subseteq A} f_{B}(x_{B})|
\leq \prod_{B \subseteq A} \| f_{B}^{\overline{2}^{|A|-|B|}} \|_{\Box(X_B)}^{1/2^{|A|-|B|}}
\end{equation}
where $x_{B} \in X_{B}$ is the restriction of $x_A$ to the indices $B$, and for any complex number $z$ we define
$z^{\overline{2}^n} := z$ when $n = 0$ and $z^{\overline{2}^n} := |z|^{2^n}$ for $n > 0$.
\end{corollary}

\begin{proof}  
For each $\omega_A \in \{0,1\}^A$ we let $f_{\omega_A}: X_A \to \C$ be the function
$$ f_{\omega_A}(x_A) := {\mathcal C}^{|\omega_A|} f_B( x_B )$$
where $B := \{ \alpha \in A: \omega_\alpha = 1 \}$.  Then we can rewrite the above left-hand side as
$$ |\E_{x^{(0)}_A, x^{(1)}_A \in X_A} \prod_{\omega_A \in \{0,1\}^A} {\mathcal C}^{|\omega_A|} f_{\omega_A}(x^{(\omega_A)}_A)|$$
which by the Gowers-Cauchy-Schwarz inequality is bounded by
$$ \prod_{\omega_A \in \{0,1\}^A} \|f_{\omega_A}\|_{\Box(X_A)}.$$
However, direct calculation (using \eqref{recursive}, for instance) shows that
$$ \|f_{\omega_A}\|_{\Box(X_A)} = \| f_{B}^{\overline{2}^{|A|-|B|}} \|_{\Box(X_B)}^{1/2^{|A|-|B|}}$$
where $B := \{ \alpha \in A: \omega_\alpha = 1 \}$, and the claim follows.
\end{proof}

As a special case of Corollary \ref{gczow-2} (together with \eqref{pointbound}), we see that
\begin{equation}\label{baby-gvn}
|\E_{x_A \in X_A} f_A(x_A) \prod_{B \subsetneq A} f_{B}(x_{B})| \leq \|f_A \|_{\Box(X_A)}
\end{equation}
whenever the functions $f_B$ are bounded in magnitude by $1$ for $B \subsetneq A$; compare this with \eqref{phase-modulation}.  
The inequality \eqref{baby-gvn} asserts that the $\Box$ norm is stable with respect to lower order functions
and can be viewed as a type of generalised von Neumann theorem.  

\begin{remark} If $f_A$ is also bounded by $1$, then there is a converse to \eqref{baby-gvn}, namely that there exist bounded functions $f_B$ for which
$$ |\E_{x_A \in X_A} f_A(x_A) \prod_{B \subsetneq A} f_{B}(x_{B})| \geq \|f_A \|_{\Box(X_A)}^{2^{|A|}}.$$
Indeed this follows easily from raising \eqref{fa} to the power $2^{|A|}$ and using the pigeonhole principle to freeze the $x^{(1)}_A$ variables.  Thus we see that the lower order functions $\prod_{B \subsetneq A} f_B(x_B)$ are ``characteristic'' for
the $\Box(X_A)$ norm: if $\Vert f_A \Vert_{\Box(X_A)}$ is large then $f_A$ correlates with a function of the form $\prod_{B \subsetneq A} f_B(x_B)$. One can pursue this idea to eventually obtain the hypergraph version of the Szemer\'edi 
regularity lemma, a task which was carried out fully in \cite{tao-hyper}.  
\end{remark}

In our applications we will need to generalise \eqref{baby-gvn} to the case where the $f_B$ are bounded by some other functions $\nu_B$. Fortunately this is also an easy consequence of Corollary \ref{gczow-2}:

\begin{corollary}[Weighted generalised von Neumann theorem]\label{gczow-gvn}  Let $(X_\alpha)_{\alpha \in A}$ be a finite collection of finite non-empty sets.  For every $B \subseteq A$ let $f_B: X_B \to \C$ and $\nu_B: X_B \to \R^+$ be functions such that $|f_B(x_B)| \leq \nu_B(x_B)$ for all $x_B \in X_B$.  Then
\begin{equation}\label{gczow-smash}
|\E_{x_A \in X_A} \prod_{B \subseteq A} f_B(x_B)|
\leq \| f_A \|_{\Box^A(\nu;X_A)} \prod_{B \subsetneq A} \| \nu_B \|_{\Box^B(\nu;X_B)}^{1/2^{|A|-|B|}}
\end{equation}
where for any $B \subseteq  A$ and $g_B: X_B \to \C$ we define $\Vert g_B \Vert_{\Box^B(\nu;X_B)}$ to be the unique nonnegative real number satisfying
$$ \|g_B\|_{\Box^B(\nu;X_B)}^{2^{|B|}} := \E_{x^{(0)}_B, x^{(1)}_B \in X_B}
\big(\prod_{\omega_B \in \{0,1\}^B} {\mathcal C}^{|\omega_B|} g_B( x^{(\omega_B)}_B) \big)
\prod_{C \subsetneq B} \prod_{\omega_C \in \{0,1\}^C} \nu_{C}( x^{(\omega_C)}_C ).$$
\end{corollary}

\begin{remark} It follows from \eqref{star-88} below that the right-hand side of the last equation is non-negative, and so $\Vert g_B \Vert_{\Box(\nu;X_B)}$ is well-defined.
Note for instance that
\begin{equation}\label{nub}
 \| \nu_B \|_{\Box^B(\nu;X_B)} := \big(\E_{x^{(0)}_B, x^{(1)}_B \in X_B}
\prod_{C \subseteq  B} \prod_{\omega_C \in \{0,1\}^C} \nu_{C}( x^{(\omega_C)}_C ) 
\big)^{1/2^{|B|}}.
\end{equation}
and
\[ \Vert f_B \Vert_{\Box(1;X_B)} = \Vert f_B \Vert_{\Box(X_B)}.\]
\end{remark}

\begin{proof} By a limiting argument we may assume that the $\nu_B$ are strictly positive throughout $X_B$.  We refactorise
$$ \prod_{B \subseteq A} f_B(x_B) = \prod_{B \subseteq A} \tilde f_B(x_B)$$
where
$$ \tilde f_B(x_B) := \frac{f_B(x_B)}{\nu_B(x_B)} \prod_{C \subseteq  B} \nu_C(x_B)^{1/2^{|A|-|C|}}.$$
Applying Corollary \ref{gczow-2} we can thus bound the left-hand side of \eqref{gczow-smash} by
$$ \| \tilde f_A \|_{\Box(X_A)} \prod_{B \subsetneq A} \| \tilde f_B^{\overline{2}^{|A|-|B|}} \|_{\Box(X_B)}^{1/2^{|A|-|B|}}.$$
However, direct calculation shows that
\begin{equation}\label{star-88} \| \tilde f_A \|_{\Box(X_A)} = \| f_A \|_{\Box( \nu;X_A)},\end{equation}
whilst the pointwise bound
$$ |\tilde f_B(x_B)| \leq \prod_{C \subseteq  B} \nu_C(x_B)^{1/2^{|A|-|C|}}$$
together with \eqref{pointbound} gives
\begin{align*}
|\| \tilde f_B^{\overline{2}^{|A|-|B|}} \|_{\Box(X_B)}^{1/2^{|A|-|B|}}| 
&\leq \big\| \prod_{C \subseteq  B} \nu_C(x_B)^{1/2^{|B|-|C|}}\big\|_{\Box(X_B)}^{1/2^{|A|-|B|}} \\
&= \| \nu_B \|_{\Box( \nu;X_B )}^{1/2^{|A|-|B|}}
\end{align*}
and the claim follows.
\end{proof}

\begin{remark} In order for this inequality to be useful, one needs to compare the weighted $\Box$ norm
$\| f \|_{\Box(\nu;X_A)}$ with the unweighted norm $\| f \|_{\Box(X_A)}$.  For any fixed set of weights $\nu$, this is not possible when the $\nu$ are unbounded; however, if the $\nu$ also depend on an additional parameter $y$, then we will be able to establish comparability estimates of this type after averaging in $y$, assuming that $\nu$ obeys suitable ``linear forms conditions''. See Appendix \ref{gvn-app}; similar ideas appear in \cite{green-tao-longprimeaps,tao-multiprime}.
\end{remark}

Now we pass from this abstract setting to a more ``additive'' setting.  Given any $s \geq 0$, any finite additive group $Z$ and any function $f: Z \to \C$, we define the \emph{Gowers uniformity norm} $\|f\|_{U^{s+1}(Z)}$ by the formula
$$ \|f\|_{U^{s+1}(Z)} := \| f( x_1 + \ldots + x_{s+1} ) \|_{\Box^{s+1}( Z^{s+1} )}.$$
Equivalently, we have
\begin{align*}
\|f\|_{U^{s+1}(Z)}^{2^{s+1}} &=  \E_{x^{(0)}, x^{(1)} \in Z^{s+1}} \prod_{\omega \in \{0,1\}^{s+1}} {\mathcal C}^{|\omega|} 
f( \sum_{j=1}^{s+1} x^{(\omega_j)}_j )  \\
&=  \E_{x \in Z; {h} \in Z^{s+1}} \prod_{\omega \in \{0,1\}^{s+1}} {\mathcal C}^{|\omega|} 
f( x + \sum_{j=1}^{s+1} \omega_j h_j ) .
\end{align*}
Because the $U^{s+1}(Z)$ norm is derived from the box norm of dimension $s+1$, many properties of the latter norm automatically descend to the former norm.  For instance, the $U^{s+1}(Z)$ norm is indeed a norm for $s \geq 1$, and from \eqref{phase-modulation} we have the
invariance
\begin{equation}\label{character}
\| e(\phi) f \|_{U^{s+1}(Z)} = \|f\|_{U^{s+1}(Z)}
\end{equation}
whenever $s \geq 1$ and $\phi: Z \to \R/\Z$ is an affine-linear phase or more generally a polynomial phase of degree at most $s$.  In our applications we shall take $Z$ to be a cyclic group $\Z_{N'}$, and our functions $f$ shall usually be real-valued.  Also, from Lemma \ref{gczow} we have the Gowers-Cauchy-Schwarz inequality for $Z$, which was first observed in \cite{gowers-long-aps} and reads as follows:
\begin{equation}\label{gcz-u}
|\E_{x \in Z; {h} \in Z^{s+1}} \prod_{\omega \in \{0,1\}^{s+1}} {\mathcal C}^{|\omega|} 
f_\omega( x + \sum_{j=1}^{s+1} \omega_j h_j )| \leq \prod_{\omega \in \{0,1\}^{s+1}} \| f_\omega\|_{U^{s+1}(Z)}.
\end{equation}

For technical reasons we shall need to localise the Gowers norms slightly.  Let $A$ be any finite non-empty subset of an additive group $Z$, which may or may not be finite.  Then for any $f: A \to \C$, we define the Gowers uniformity norm $\|f\|_{U^{s+1}(A)}$ by the formula
\begin{equation}\label{ukdef}
\begin{split}
\|f\|_{U^{s+1}(A)}^{2^{s+1}} &= \E_{x^{(0)}, x^{(1)}: \sum_{j=1}^{s+1} x^{(\omega_j)}_j \in A \;\;\forall \omega \in \{0,1\}^{s+1} } \prod_{\omega \in \{0,1\}^{s+1}} {\mathcal C}^{|\omega|} f( \sum_{j=1}^{s+1} x^{(\omega_j)}_j )  \\
&= \E_{x,{h}: x + \sum_{j=1}^{s+1} \omega_j h_j \in A \;\;\forall \omega \in \{0,1\}^{s+1}} \prod_{\omega \in \{0,1\}^{s+1}} {\mathcal C}^{|\omega|} f( x + \sum_{j=1}^{s+1} \omega_j h_j ) .
\end{split}
\end{equation}
In the particular case $A = [N]$, which is used several times in the paper, we shall adopt the abbreviation
\[ \Vert f \Vert_{U^{s+1}[N]} := \Vert f \Vert_{U^{s+1}([N])}.\]

If $A$ is contained in a \emph{finite} additive group $Z$, then these local Gowers norms are related to their global counterparts by the identity
\begin{equation}\label{gowers-compat}
\| f \|_{U^{s+1}(A)} = \| f 1_A \|_{U^{s+1}(Z)} / \| 1_A \|_{U^{s+1}(Z)}
\end{equation}
for any $f: A \to \C$, where $f 1_A: Z \to \C$ is the extension by zero of $f$ from $A$ to $Z$.  The local norm $U^{s+1}(A)$ is also intrinsic in the following sense: if $A \subseteq Z$, $A' \subseteq Z'$, and $\phi: A \to A'$ is a \emph{Freiman isomorphism} in the sense that it is 1-1 onto its image and for any $a_1,a_2,a_3,a_4 \in A$, we have $a_1+a_2=a_3+a_4$ if and only if $\phi(a_1)+\phi(a_2)=\phi(a_3)+\phi(a_4)$, then we have $\| f \circ \phi \|_{U^{s+1}(A)} = \| f \|_{U^{s+1}(A')}$ for all $f: A' \to \C$.  A particular consequence of this is the following lemma.

\begin{lemma}[Comparability of $U^{s+1}(I)$ and $U^{s+1}(\Z_{N'})$]\label{normlemma}  Let $N' \geq 1$ be an integer, let $\alpha > 0$, and let $I = \{ a, a+1, \ldots, b\}$ be an interval of integers whose length satisfies $\alpha N' \leq |I| \leq N'/2$.  Let $f: I \to \C$ be a function on $I$, and let $\tilde f: \Z_{N'} \to \C$ be the function formed from $f$ by identifying $I$ with a subset of $\Z_{N'}$ and setting $\tilde f(x) = 0$ for $x \notin I$.  Then we have
\begin{equation}\label{norm-compare}
\| \tilde f \|_{U^{s+1}(\Z_{N'})} = c \| f \|_{U^{s+1}(I)}
\end{equation}
where $c = c_{I,N',s} > 0$ is a constant which is independent of $f$, and which is bounded above and below by quantities depending only on $\alpha$ and $s$.
\end{lemma}

\begin{proof} As $|I| \leq N'/2$, the interval $I \subseteq \Z$ is Freiman isomorphic to its counterpart in $\Z_{N'}$.  The claim then follows from \eqref{gowers-compat} together the easily confirmed observation that $\| 1_I \|_{U^{s+1}(\Z_{N'})}$ is bounded above and below by quantities depending only on $\alpha$ and $s$.
\end{proof}

\begin{remark} We will typically apply this lemma with $I = [N]$ and with $N'$ comparable to a moderately large multiple of $N$. See, for example, the proof of Proposition \ref{gow-ps}.
\end{remark}

\section{Proof of the generalised von Neumann theorem}\label{gvn-app}

The purpose of this appendix is to prove Proposition \ref{gvn}. 

\begin{gvn-restate}[Generalised von Neumann theorem] Let $s,t,d,L$ be positive integer parameters as usual. Then there are constants $C_1$ and $D$, depending on $s,t,d$ and $L$, such that the following is true.  Let $C$, $C_1 \leq C \leq O_{s,t,d,L}(1)$, be arbitrary and suppose that $N' \in [CN,2CN]$ is a prime. Let $\nu:\Z_{N'} \to \R^+$ be a $D$-pseudorandom measure, and suppose that $f_1,\dots, f_t : [N] \rightarrow \R$ are functions with $|f_i(x)| \leq \nu(x)$ for all $i \in [t]$ and $x \in [N]$. Suppose that $\Psi = (\psi_1,\ldots,\psi_t)$ is a system of affine-linear forms in $s$-normal form with $\|\Psi\|_N \leq L$.  Let $K \subseteq  [-N,N]^d$ be a convex body such that $\Psi(K) \subseteq  [N]^t$.  Suppose also that
$$ \min_{1 \leq j \leq t}  \Vert f_j \Vert_{U^{s+1}[N]} \leq \delta$$
for some $\delta > 0$. Then we have
\begin{equation}
\sum_{{n} \in K} \prod_{i \in [t]} f_i(\psi_i({n})) =  \; 
o_{\delta}(N^d) + \kappa(\delta)N^{d}.\end{equation}
\end{gvn-restate}

Recall that this is a variant of \cite[Proposition 5.3]{green-tao-longprimeaps}, which was proven by a long series of applications of the Cauchy-Schwarz inequality.  We shall phrase our argument using Corollary \ref{gczow-2}, but the argument is essentially that of \cite[\S 5]{green-tao-longprimeaps}. It is also necessary to perform some regularisation to deal with the convex body $K$, a technical feature not present in \cite[Proposition 5.3]{green-tao-longprimeaps}.

\textsc{Moving to a cyclic group.} Let us first make some very minor reductions. We start by moving the whole problem to the group $\Z_{N'}$. We will always assume that $N' = O_{s,t,d,L}(N)$, but one may wish to take $N'$ to be quite a bit larger than $N$ in order that a pseudorandom measure $\nu$ can be constructed so as to make Proposition \ref{gvn} applicable. We embed $[N]$ inside $\Z_{N'}$ in the usual manner, and extend the functions $f_1,\ldots,f_t$ to all of $\Z_{N'}$ by defining them to be zero outside of $[N]$. From Lemma \ref{normlemma} we then have
$$
 \Vert f_j \Vert_{U^{s+1}(\Z_{N'})} \ll_C \delta
$$
for some $j \in \{1,\dots,t\}$.
Similarly, we may identify the set $K \cap \Z^d$ with a subset $K'$ of 
$\Z_{N'}^d$.  We can also view $\Psi$ as a map from $\Z_{N'}^d$ to $\Z_{N'}^t$.  Note that $\Psi$ will then map $K'$ to
$[N]^d$. To summarise, we have reduced matters to establishing the following.

\begin{gvn-restate-again}[Transfer to $\Z_{N'}$] 
Let $s,t,d,L$ be positive integer parameters as usual. Then there is a constant $D$, depending on $s,t,d$ and $L$, such that the following is true.  Let $\nu:\Z_{N'} \to \R^+$ be a $D$-pseudorandom measure, and suppose that $f_1,\dots, f_t : \Z_{N'} \rightarrow \R$ are functions with $|f_i(x)| \leq \nu(x)$ for all $i \in [t]$ and $x \in \Z_{N'}$. Suppose that $\Psi = (\psi_1,\ldots,\psi_t)$ is a system of affine-linear forms in $s$-normal form with $\|\Psi\|_N \leq L$. Let $K' \subseteq \Z_{N'}^d$ be identified with $K \cap \Z^d$ for some convex $K \subseteq [-\frac{1}{4}N',\frac{1}{4}N']^d$.  Suppose also that
$$ \min_{1 \leq j \leq t}  \Vert f_j \Vert_{U^{s+1}(\Z_{N'})} \leq \delta$$
for some $\delta > 0$. Then we have
\begin{equation}\label{lin-w-prime} 
\E_{{n} \in \Z_{N'}^d} 1_{K'}({n}) \prod_{i \in [t]} f_i(\psi_i({n})) =  \; 
o_{\delta}(1) + \kappa(\delta).\end{equation}
\end{gvn-restate-again}

\begin{remark}
Note the disappearance of $C$. This was an artefact of the relationship between $N$ and $N'$, which has now been forgotten.
\end{remark}

From this point onwards we do our linear algebra over $\Z_{N'}$, rather than over $\Q$. Note that the notion of $s$-normal form coincides in the two settings provided that $N' \geq N_0(s,t,d,L)$ is sufficiently large. Furthermore no two of the homogeneous parts  $\dot{\psi}_i$ are parallel when considered $\md{N'}$. This fact (which is very easily checked) is a simple instance of a kind of ``Lefschetz principle''.

\textsc{Removing the convex cutoff.} The next step is to partially eliminate the cutoff $1_{K'}( n)$ by replacing it by a more analytically tractable Lipschitz cutoff.  We introduce a metric on $\Z_{N'}^d$ by declaring the distance between $(n_1,\ldots,n_d)$ and $(m_1,\ldots,m_d)$ to be $(\sum_{j=1}^d \| \frac{n_i-m_i}{N'} \|_{\R/\Z}^2)^{1/2}$, where $\|x\|_{\R/\Z}$ denotes the distance to the nearest integer. This is the metric induced from the standard embedding of $\Z_{N'}^d$ into the torus $(\R/\Z)^d$.  To establish Proposition $\mbox{\ref{gvn}}^{\prime}$, we claim that it suffices
to establish the bound
\begin{equation}\label{lipcon}
\E_{{n} \in \Z_{N'}^d} F( n) \prod_{i \in [t]} f_i(\psi_i({n})) = o_{\delta,M}(1) + \kappa_M(\delta)
\end{equation}
whenever $M > 0$, $F: \Z_{N'}^d \to [-1,1]$ has Lipschitz constant $M$ and the functions $f_i$ are bounded pointwise by $\nu$ and satisfy $\min_{1 \leq i \leq t}\Vert f_i \Vert_{U^{s+1}(\Z_{N'})} \leq \delta$.  To see why, let $\eps > 0$ be a small quantity to be chosen later.  It will suffice to prove that
$$
\E_{{n} \in \Z_{N'}^d} 1_{K'}( n) \prod_{i \in [t]} f_i(\psi_i({n})) =  \; 
o_{\eps}(1) +  \kappa_{\eps}(\delta) + \kappa(\eps),$$
as the claim then follows by setting $\eps$ to be a sufficiently slowly decaying function of $\delta$.

To establish this bound, we apply Corollary \ref{lip-approx} to effect the decomposition
$$ 1_{K'}( n) = F_\eps(  n ) + O( G_\eps(  n ) )$$
for all $ n \in \Z_{N'}^d$, where $F_\eps, G_\eps: \Z_{N'}^d \to [0,1]$ are Lipschitz in the above metric with constant $O(1/\eps)$.  Furthermore, from the Lipschitz and integral bounds in Corollary \ref{lip-approx} we easily obtain the estimate
\begin{equation}\label{lipest}
 \E_{ n \in \Z_{N'}} G_\eps( n) = o_{\eps}(1) + \kappa(\eps).
\end{equation}
Here we are basically using nothing more than the standard fact that Lipschitz functions are uniformly Riemann integrable.
From \eqref{lipcon} we have
$$ \E_{{n} \in \Z_{N'}^d} F_\eps( n) \prod_{i \in [t]} f_i(\psi_i({n})) = 
o_{\eps}(1)  + \kappa_{\eps}(\delta)$$
and so by the triangle inequality and the fact that $|f_i(x)| \leq \nu(x)$ it is enough to show that
\begin{equation}\label{gbomb}
\E_{{n} \in \Z_{N'}^d} G_\eps( n) \prod_{i \in [t]} \nu(\psi_i({n})) = 
o_{\eps}(1)  + \kappa(\eps).
\end{equation}
Now a standard application of the linear forms condition (see \cite[Lemma 5.2]{green-tao-longprimeaps}) gives
$$ \| \nu - 1 \|_{U^{s+1}(\Z_{N'})} = o(1).$$
Now the function $\frac{1}{2}(\nu - 1)$ satisfies $\frac{1}{2}|\nu(x) - 1| \leq \frac{1}{2}(\nu(x) + 1)$, and this latter function is easily seen to be a pseudorandom measure (see \cite[Lemma 3.4]{green-tao-longprimeaps}). Thus from \eqref{lipcon} we have
\[  
\E_{{n} \in \Z_{N'}^d} G_{\eps}( n) \prod_{i \in [t]} g_i(\psi_i({n})) = o_{\eps}(1)
\]
whenever all the functions $g_i$ are either $1$ or $\nu - 1$, and not all of them are $1$. When $g_i = 1$ for all $i$ we have the bound $o_{\eps}(1) + \kappa(\epsilon)$, from \eqref{lipest}.
The bound \eqref{gbomb} now follows immediately upon writing $\nu = 1 + (\nu - 1)$ and expanding as a sum of $2^t$ terms.

It remains to prove \eqref{lipcon}.  We now claim that we may dispense with the Lipschitz cutoff $F$ entirely, and reduce to proving the estimate
\begin{equation}\label{nolipcon}
\E_{{n} \in \Z_{N'}^d} \prod_{i \in [t]} f_i(\psi_i({n})) =  \; 
o_{\delta}(1) + \kappa(\delta),
\end{equation}
which involves no cutoff function at all.
To see this, first observe that \eqref{nolipcon} implies the extension
\begin{equation}\label{nolipcon-fourier}
\E_{{n} \in \Z_{N'}^d} e( m \cdot  n / N) \prod_{i \in [t]} f_i(\psi_i({n})) =  \; 
o_{\delta}(1) + \kappa(\delta).
\end{equation}
for any frequency $m \in \Z_{N'}^d$.  Indeed, if $m$ lies in the span of $\dot \psi_1, \ldots, \dot \psi_t$ then we may simply 
factor $e( m \cdot  n / N)$ into terms that can be absorbed into the $f_1,\ldots,f_t$ factors, noting that we can trivially extend \eqref{nolipcon} to cover the case when $f_1,\ldots,f_t$ are complex-valued instead of real-valued.  If $m$ does not lie in this span, then it is easy to see that the left-hand side of \eqref{nolipcon-fourier} in fact vanishes.

Now we return to \eqref{lipcon}. Let $X > 0$ be arbitrary.
By a standard Fourier-analytic argument, given in detail in \cite[Lemma A.9]{green-tao-u3mobius}, we may decompose
$$ F( n) = \sum_{j=1}^J c_j e( m_j \cdot  n / N ) + O_d(M\log X/X)$$
where $J = O_d(X^{d})$, $c_j = O(1)$ are coefficients, and $m_j \in \Z_{N'}^d$ are frequencies.  Inserting this into \eqref{lipcon}, we have
\begin{align*}
&\E_{{n} \in \Z_{N'}^d} F({n}) \prod_{i \in [t]} f_i (\psi_i({n})) \\
&= \sum_{j =1}^Jc_j \E_{{n} \in \Z_{N'}^d}e( m_j \cdot  n / N )\prod_{i \in [t]} f_i (\psi_i({n})) + O_d\big( \frac{M\log X}{X} \big)\E_{{n} \in \Z_{N'}^d} \prod_{i \in [t]} \nu(\psi_i({n})).
\end{align*}
Using \eqref{nolipcon-fourier} to control the first term and the linear forms condition to estimate the second, we see that this is bounded by
\[ O_d(X^d) (o_{\delta}(1) + \kappa(\delta)) + O_d\big( \frac{M\log X}{X} \big)(1 + o(1)).\]
Taking $X$ to be a sufficiently slowly growing function of $1/\delta$ we obtain \eqref{lipcon} as desired.

\textsc{Main argument.} It remains to prove \eqref{nolipcon}.  By symmetry we may assume that $f_1$ is the function with minimal $U^{s+1}$ norm, thus
$$
\Vert f_1 \Vert_{U^{s+1}(\Z_{N'})} \leq \delta.
$$
Recall that the system $\Psi: \Z^d \to \Z^t$ is in $s$-normal form.  By permuting the basis vectors $e_1,\ldots,e_d$ if necessary, we may then assume that $\prod_{j=1}^{s+1} \dot \psi_i(e_j)$ vanishes for $i \neq 1$ and is non-zero for $i=1$.  

In summary, we are reduced to proving

\begin{gvn-restate-yetagain}[Reduced generalised von Neumann theorem] 
Let $s,t,d,L$ be positive integer parameters as usual. Then there is a constant $D$, depending on $s,t,d$ and $L$, such that the following is true.  Let $\nu:\Z_{N'} \to \R^+$ be a $D$-pseudorandom measure, and suppose that $f_1,\dots, f_t : \Z_{N'} \rightarrow \R$ are functions with $|f_i(x)| \leq \nu(x)$ for all $i \in [t]$ and $x \in \Z_{N'}$. Suppose that $\Psi = (\psi_1,\ldots,\psi_t)$ is a system of affine-linear forms such that
$\prod_{j=1}^{s+1} \dot \psi_i(e_j)$ vanishes for $i \neq 1$ and is non-zero for $i=1$.
Then we have
\begin{equation}\label{rgvn-est}
\big|\E_{{n} \in \Z_{N'}^d} \prod_{i \in [t]} f_i(\psi_i({n}))\big| \leq 
\| f_1 \|_{U^{s+1}(\Z_{N'})} + o(1).
\end{equation}
\end{gvn-restate-yetagain}

To prove the estimate \eqref{rgvn-est}, note first that the coefficients $\dot \psi_1(e_j)$, $j \in [s+1]$, are non-zero and bounded by $O_{s,t,d,L}(1)$, and hence are invertible in $\Z_{N'}$ provided that $N \geq N_0(s,t,d,L)$.  Thus we may dilate\footnote{This dilation converts the coefficients from bounded integers, to rationals with bounded numerator and denominator.  However, when the time comes to apply the linear forms condition, one can clear denominators and reduce back to estimates involving only bounded integers again.} the first $s+1$ variables and assume that $\dot \psi_1(e_j) = 1$ for $j \in [s+1]$, a manoeuvre which affords a little notational simplicity if nothing more. With this normalisation we have, writing $n = (x_1,\dots, x_d)$ and $y = (x_{s+2},\dots, x_d)$, that
$$ \psi_1( x_1,\ldots,x_{s+1},y) = x_1 + \ldots + x_{s+1} + \psi_1(0,y).$$
The other forms $\psi_i$, $i = 2,\dots,t$ do not involve all of the variables $x_1,\dots, x_{s+1}$, since the system $\Psi$ is in normal form. This will be a crucial fact for us and to handle it we look, for each $\psi_i$, at the set $\Omega(i)$ of indices $j \in [s+1]$ for which $\dot \psi_i (e_j) \neq 0$, and then group the forms according to their associated set $\Omega(i)$. Thus $\Omega(1) = [s+1]$ and $\Omega(i) \subsetneq [s+1]$ for $i = 2,\dots,t$. Observe that the indices $j = s+2,\dots,d$ and the associated variable $y = (x_{s+2},\dots,x_d)$ will be largely irrelevant in the sequel. With this nomenclature we may write the left-hand side of \eqref{rgvn-est} as
\begin{equation}\label{rgvn-2}
| \E_{y \in \Z_{N'}^{d-s-1}} \E_{x_{[s+1]} \in \Z_{N'}^{s+1}} \prod_{B \subseteq  [s+1]} F_{B,y}( x_B ) |
\end{equation}
where $x_{[s+1]} = (x_j)_{j \in [s+1]}$, $x_B$ is the restriction of $x_{[s+1]}$ to $B$, and
$$ F_{B,y}(x_B) := \prod_{i \in [t]: \Omega(i) = B} f_i(\psi_i(x_B,y)).$$
We have abused notation ever so slightly by regarding $f_i$ as a function on $\Z_{N'}^B \times \Z_{N'}^{d-s-1}$ rather than on $\Z_{N'}^{s+1} \times \Z_{N'}^d$, supressing mention of the irrelevant variables $x_j$, $j \in [s+1]\setminus \Omega(i)$. Observe that
$$ F_{[s+1],y}(x_{[s+1]}) = f_1( \psi_1( x_{[s+1]}, y ) ) = f_1( x_1 + \ldots + x_{s+1} + \psi_1(0,y) ).$$
Now we have the pointwise bounds $|F_{B,y}(x_B)| \leq \nu_{B,y}(x_B)$, where
$$ \nu_{B,y}(x_B) := \prod_{i \in [t]: \Omega(i) = B } \nu(\psi_i(x_B,y)).$$
Invoking Corollary \ref{gczow-gvn}, we may bound \eqref{rgvn-2} by
$$ \E_{y \in \Z_{N'}^{d-s-1}} \| F_{[s+1],y} \|_{\Box( \nu_{[s+1],y}; \Z_{N'}^{[s+1]})} 
\prod_{B \subsetneq [s+1]}
\| \nu_{B,y} \|_{\Box( \nu_{B,y}; \Z_{N'}^B )}^{1/2^{s+1-|B|}}.$$
The reader may wish to recall the definition of the quantities appearing here, which are provided in the statement of Corollary \ref{gczow-gvn}. 

Applying H\"older's inequality\footnote{This is really an application of the Cauchy-Schwarz inequality several times, since the exponent is a power of two.}, we see that to show \eqref{rgvn-est} it suffices to show that
\begin{equation}\label{rgv-1}
\E_{y \in \Z_{N'}^{d-s-1}} \| F_{[s+1],y} \|_{\Box^{s+1}( \nu_{[s+1],y};\Z_{N'}^{[s+1]} )}^{2^{s+1}}
\leq
\| f_1 \|_{U^{s+1}(\Z_{N'})} + o(1)
\end{equation}
and that
\begin{equation}\label{rgv-2}
\E_{y \in \Z_{N'}^{d-s-1}} \| \nu_{B,y} \|_{\Box( \nu_{\subsetneq B,y} )}^{2^{|B|}}
= 1 + o(1)
\end{equation}
for all non-empty $B \subsetneq [s+1]$. Note that except for $f_1$, the unknown functions $f_2,\ldots,f_t$ have all been eliminated. This procedure will be familiar to readers who have looked at (for example) \cite[Ch. 5]{green-tao-longprimeaps}.

We begin with \eqref{rgv-2}.  We expand the left-hand side, obtaining
\begin{align*}
\E_{y \in \Z_{N'}^{d-s-1}} \| \nu_{B,y} \|_{\Box^B( \nu_{\subsetneq B,y} )}^{2^{|B|}} &= \E_{x_B^{(0)}, x_B^{(1)} \in X_B} \prod_{C \subseteq B} \prod_{\omega_C \in \{0,1\}^C} \nu_{C,y}(x_C^{(\omega_C)}) \\ &= \E_{x_B^{(0)}, x_B^{(1)} \in X_B} \prod_{C \subseteq B} \prod_{\omega_C \in \{0,1\}^C} \prod_{i \in [t] : \Omega(i) = C}\nu(\psi_i(x_C^{(\omega_C)},y)).
\end{align*}
Because of the definition of $\Omega(i)$, and the hypothesis that no two of the $\psi_i$ were affine-linear combinations of each other, we see that the affine-linear forms
$$ (x^{(0)}_B, x^{(1)}_B, y) \mapsto \psi_i( x^{(\omega_C)}_C, y ),$$
as $C$ varies over subsets of $B$ and $i$ varies over those $i \in [t]$ such that $\Omega(i) = C$, also have the property that no two forms are affine-linear combinations of each other. In other words, this system has finite complexity.  Thus \eqref{rgv-2} will follow from the linear forms condition \eqref{lfc} provided that the degree $D$ of pseudorandomness is sufficiently large.

Now we turn to \eqref{rgv-1}.  The left-hand side expands as
\begin{align*} \E_{x^{(0)}_{[s+1]}, x^{(1)}_{[s+1]} \in \Z_{N'}^{s+1}; y \in \Z_{N'}^{d-s-1}}
\!\!\!\!\prod_{\omega \in \{0,1\}^{s+1}} f_1 & ( \sum_{j=1}^{s+1} x^{(\omega_j)}_i + \psi_1(0,y) ) \times \\ & \times
\prod_{C \subsetneq [s+1]} \prod_{\omega_C \in \{0,1\}^C} \prod_{i \in [t] : \Omega(i) = C} \nu( \psi_i( x^{(\omega_C)}_C, y ) ).\end{align*} 
Substituting ${h} := x^{(1)}_{[s+1]} - x^{(0)}_{[s+1]}$ and $z := x^{(0)}_1 + \ldots + x^{(0)}_{s+1} + \psi_1(0,y)$, we may 
rewrite this as
\begin{align*} \E_{x^{(0)}_{[s+1]}, {h} \in \Z_{N'}^{s+1}; y \in \Z_{N'}^{d-s-1}}
\prod_{\omega \in \{0,1\}^{s+1}} f_1( z + \sum_{j=1}^{s+1} \omega_j h_j )
\prod_{C \subsetneq [s+1]} \prod_{\omega_C \in \{0,1\}^C} &\prod_{i \in [t]: \Omega(i) = C} \nu\big( \psi_i( x^{(0)}_C, y ) + \\ & + \sum_{j \in C} \omega_j \dot \psi_i( e_j ) h_j \big).\end{align*}
Observe that for fixed ${h}$, the map $(x^{(0)}_{[s+1]},y) \mapsto z$ is uniform, in the sense that each $z$ is mapped to by
exactly $(N')^{d-1}$ preimages.  Thus we may rewrite the preceding expression as
$$ \E_{z \in \Z_{N'}; {h} \in \Z_{N'}^{s+1}} W(z,{h}) \prod_{\omega \in \{0,1\}^{s+1}} f_1( z + \sum_{j=1}^{s+1} \omega_j h_j ) $$
where
\[
 W(z,{h}) := \E_{\substack{x^{(0)}_{[s+1]} \in \Z_{N'}^{s+1}; y \in \Z_{N'}^{d-s-1} \hfill\\ z = x^{(0)}_1 + \ldots + x^{(0)}_{s+1} + \psi_1(0,y)}}
\prod_{C \subsetneq [s+1]} \prod_{\omega_C \in \{0,1\}^C} \prod_{i \in [t]: \Omega(i) = C} \nu( \psi_i( x^{(0)}_C, y ) + \sum_{j \in C} \omega_j \dot \psi_i( e_j ) h_j ).
\]
Comparing this with \eqref{ukdef}, we see that to prove \eqref{rgv-1} it suffices to show that
$$ \E_{z \in \Z_{N'}; {h} \in \Z_{N'}^{s+1}} (W(z,{h})-1) \prod_{\omega \in \{0,1\}^{s+1}} f_1( z + \sum_{j=1}^{s+1} \omega_j h_j )  = o(1).
$$
By Cauchy-Schwarz and the hypothesis $|f_1(x)| \leq \nu(x)$, it suffices to establish the estimates
$$ \E_{z \in \Z_{N'}; {h} \in \Z_{N'}^{s+1}}  |W(z,{h})-1|^n \prod_{\omega \in \{0,1\}^{s+1}} \nu( z + \sum_{j=1}^{s+1} \omega_j h_j ) = 0^n + o(1)$$
for $n=0$ and $n=2$. Expanding, we reduce to showing that
$$ \E_{z \in \Z_{N'}; {h} \in \Z_{N'}^{s+1}} W(z,{h})^n \prod_{\omega \in \{0,1\}^{s+1}} \nu( z + \sum_{j=1}^{s+1} \omega_j h_j )  = 1 + o(1)$$
for $n=0,1,2$.

This will follow from the linear forms condition.  We shall just verify the case $n=2$, as the cases $n=0,1$ follow from that case (they utilise a subset of the linear forms that are used in the $n=2$ case).  When $n=2$, we can expand out the left-hand side as
\begin{equation}\label{enu}
\begin{split}
& \E_* \big(\prod_{\omega \in \{0,1\}^{s+1}} \nu( z + \sum_{j=1}^{s+1} \omega_j h_j )\big) \times\\
&\prod_{C \subsetneq [s+1]} \prod_{\omega_C \in \{0,1\}^C} \prod_{i \in [t] : \Omega(i) = C} \nu\big( \psi_i( x^{(0)}_C, y \big) + \sum_{j \in C} \omega_j \dot \psi_i( e_j ) h_j ) \nu\big( \psi_i( \tilde x^{(0)}_C, \tilde y ) + \sum_{j \in C} \omega_j \dot \psi_i( e_j ) h_j\big)
\end{split}
\end{equation}
where the average $\E_*$ is over all sextuples
$$ (z, {h}, x^{(0)}_{[s+1]}, \tilde x^{(0)}_{[s+1]}, y, \tilde y ) \in \Z_{N'} \times \Z_{N'}^{s+1} \times \Z_{N'}^{s+1}
\times \Z_{N'}^{s+1} \times \Z_{N'}^{d-s-1} \times \Z_{N'}^{d-s-1}$$
subject to the affine constraints
\begin{equation}\label{constraints}
 z = x^{(0)}_1 + \ldots + x^{(0)}_{s+1} + \psi_1(0,y) = \tilde x^{(0)}_1 + \ldots + \tilde x^{(0)}_{s+1} + \psi_1(0,\tilde y).
\end{equation}

Naturally, we wish to apply the linear forms condition, on the assumption that $\nu$ is $D$-pseudorandom for sufficiently large $D$. To do this we must first eliminate the constraints \eqref{constraints}. To do this, we substitute for $x_{s+1}^{(0)}$ and $\tilde x_{s+1}^{(0)}$ in terms of the other variables, that is to say we write
\[ x_{s+1}^{(0)} = z - x_1^{(0)} - \dots - x_{s}^{(0)} - \psi_1(0,y)\]
and
\[ \tilde x_{s+1}^{(0)} = z - \tilde x_1^{(0)} - \dots - \tilde x_{s}^{(0)} - \psi_1(0,\tilde y).\]
In this way we may rewrite \eqref{enu} as an unconstrained average over the $2d + s - 1$ variables ${h}, x^{(0)}_{[s]}, \tilde x^{(0)}_{[s]}, y, \tilde y$.

When written in terms of this set of variables, it is clear that all the linear forms in \eqref{enu} have integer coefficients which are bounded in terms of $s,t,d$ and $L$. To apply the linear forms condition, all we must do is satisfy ourselves that no two of these forms are affinely dependent, that is to say no two of them have parallel homogeneous parts.

To see this, first observe that the $2^{s+1}$ homogeneous forms $z + \sum_{j=1}^{s+1} \omega_j h_j$ are pairwise distinct, and that they are also different from any other form appearing in \eqref{enu} even after performing the above substitutions, because the latter forms all involve at least one of the variables from $x^{(0)}_{[s]}$, $\tilde x^{(0)}_{[s]}$ (here we are using the fact that $C$ is a \emph{proper} subset of $[s+1]$).
	
Now consider an affine form $\psi_i( x^{(0)}_C, y ) + \sum_{j \in C} \omega_j \dot \psi_i( e_j ) h_j$ appearing in \eqref{enu}. Recalling that $C = \Omega(i)$, the set of all $j$ for which $\dot \psi_i(e_j) \neq 0$, we see that in our new system of variables this form may be written as the slightly alarming expression
\begin{equation}\label{eq4.776} \sum_{j = 1}^s \dot \psi_i(e_j) (x_j^{(0)} + \omega_j h_j) + \dot \psi_i(e_{s+1})  (z - x_1^{(0)} - \dots - x_{s}^{(0)} - \psi_1(0,y) + \omega_{s+1}h_{s+1}) + \psi_i(0,y).\end{equation}
There is a similar expression involving tildes. We claim first of all that at least one of the variables $x_1^{(0)},\dots,x_s^{(0)}$ must appear with non-zero coefficient. If this were not the case then we would have $\dot \psi_i(e_j) = \dot \psi_i(e_{s+1}) $ for $j \leq s$ and hence, since $C \subsetneq [s+1]$, $C$ is empty. Hence so is the product over $\omega_C \in \{0,1\}^C$ in \eqref{enu}. Thus no form \eqref{eq4.776} with this property appears in \eqref{enu}, thereby confirming the claim.

The claim just proved immediately implies that no form \eqref{eq4.776} has homogeneous part parallel to that of a form with tildes. It remains to prove that the forms in \eqref{eq4.776} have pairwise non-parallel homogeneous parts.

Suppose that we are given a form \eqref{eq4.776} written as
\[ q\big(r_1 x_1^{(0)} + \dots + r_s x_s^{(0)} + r'z + (\mbox{terms involving $h$ and $y$})\big)\] where $q \neq 0$. We claim that
the set $C$ from which the form came may be identified knowing only $r_1,\dots,r_s,r'$. Indeed we must have $qr' = \dot \psi_i(e_{s+1})$, whence $\dot \psi_i(e_j) = \lambda(r_j + r')$ for $j \leq s$. The set $C$ may be found simply by looking at which of these quantities do not vanish. It is immediately clear that $\omega_C \in \{0,1\}^C$ may also be recovered.

The only way in which two forms \eqref{eq4.776} could have parallel homogeneous parts, then, is if there is some fixed choice of $\omega$, some $i \neq i'$ and some rational $q,q' \neq 0$ such that
\begin{align*} 
& q \big(\sum_{j = 1}^s \dot \psi_i(e_j) (x_j^{(0)} + \omega_j h_j) + \dot \psi_i(e_{s+1})  (z - \sum_{j=1}^s x_j^{(0)}  - \psi_1(0,y) + \omega_{s+1}h_{s+1}) + \psi_i(0,y)\big) \\ &=q'\big(\sum_{j = 1}^s \dot \psi_{i'}(e_j) (x_j^{(0)} + \omega_j h_j) + \dot \psi_{i'}(e_{s+1})  (z - \sum_{j=1}^s x_j^{(0)} - \psi_1(0,y) + \omega_{s+1}h_{s+1}) + \psi_{i'}(0,y)\big)
\end{align*}
for all choices of the variables. After some simple manipulations one confirms that $q \dot \psi_i(e_j) = q' \dot \psi_{i'}(e_j)$ for $j \leq s+1$ and that $q \psi_i(0,y) = q' \psi_{i'}(0,y)$. Thus $\dot \psi_i$ is parallel to $\dot \psi_{i'}$, contrary to the assumption that the system $\Psi = (\psi_i)_{i=1}^t$ has finite complexity.

We have verified that it was valid to invoke the linear forms condition, provided that $D$ is large enough. This completes the proof of Proposition \ref{gvn}' and hence that of Proposition \ref{gvn}.
\endproof

\section{Goldston-Y{\i}ld{\i}r{\i}m correlation estimates}\label{gy-sec}

One aim of this section is to construct a pseudorandom measure $\nu$ such that a suitable multiple of $\nu$ majorises the modified von Mangoldt function $\Lambda'_{b,W}$. Specifically, we will prove Proposition \ref{pseudodom}. This was essentially carried out in \cite[Chs. 9, 10]{green-tao-longprimeaps}, building on work of Goldston and Y{\i}ld{\i}r{\i}m \cite{gy-1,gy-2,gy-3,pintz}, but the argument there only led to a majorant for \emph{one} function $\Lambda'_{b,W}$, whereas in the present work we need to simultaneously majorise $\Lambda'_{b_1,W},\dots,\Lambda'_{b_t,W}$. A few small modifications to the argument in \cite{green-tao-longprimeaps} would, however, achieve this. Another aim of this section is to prove \eqref{to-prove-appD}, a crucial estimate on the Gowers norm of a certain truncated von Mangoldt function $\Lambda^{\sharp}$. This does not follow immediately from the results in \cite{green-tao-longprimeaps}, though can be proved using similar ideas. We take the opportunity to give a brief but more-or-less self-contained account of these ideas here, while also providing some simplifications.

The heart of the matter is the establishment of correlation estimates for truncated divisor sums $\Lambda_{\chi,R,a}: \Z \to \R$
of the form
$$ \Lambda_{\chi,R,a}(n) := \log R (\sum_{d|n} \mu(d) \chi(\frac{\log d}{\log R}))^a.$$
In this expression $R$ is a moderately large number, which in practice will be a small power of $N$, $\chi: \R \to \R$ is a smooth, compactly supported function, and $a \in \N$. In our applications we only ever take $a=1$ or $a=2$.  We extend $\Lambda_{\chi,R,a}$ to the negative numbers in the obvious manner. Indeed, the compact support of $\chi$ ensures that $\Lambda_{\chi,R,a}$ is periodic.

\begin{remark} Observe that $\Lambda_{\chi,R,a} = \chi(0)^a \log R$ on ``almost primes'' - numbers coprime to $\prod_{p \leq R} p$.  For the purposes of gaining intuition about these functions one might think of them heuristically as being weights on the almost primes, though they do also have some weight on other numbers.  The reason we need to deal with $\Lambda_{\chi,R,2}(n)$ is to correct for the rather unfortunate fact that $\Lambda_{\chi,R,1}(n)$ can be negative. This trick is of course closely related to the $\Lambda^2$ sieve of Selberg.  
\end{remark}

Associated to these truncated divisor sums are certain numbers which we call \emph{sieve factors}.

\begin{definition}[Sieve factors]\label{sieve-factor}  Let $\chi: \R \to \R$ be compactly supported and suppose that $a \geq 1$.  Then we define the \emph{sieve factor} $c_{\chi,a}$ by the formula
\begin{equation}\label{chi-def}
c_{\chi,a} := 
\int_\R \ldots \int_\R \prod_{B \subseteq [a]} \big( \sum_{j \in B} (1+i \xi_j) \big)^{(-1)^{|B|-1}} 
\prod_{j = 1}^a \varphi(\xi_j)\ d\xi_j,
\end{equation}
where $\varphi$ is the modified Fourier transform of $\chi$, defined by the formula
\begin{equation}\label{varphi-def}
 e^x \chi(x) = \int_{-\infty}^\infty \varphi(\xi) e^{-ix\xi}\ d\xi.
\end{equation}
\end{definition}

The sieve factor $c_{\chi,a}$ looks very complicated (though explicitly computable), but in the special cases $a=1,2$ it has a particularly simple form:

\begin{lemma}\label{sieve-compute}  We have $c_{\chi,1} = - \chi'(0)$ and $c_{\chi,2} = \int_0^\infty |\chi'(x)|^2\ dx$.  More generally, $c_{\chi,a}$ is a real number.
\end{lemma}

\begin{proof} We deal first with the case $a = 1$.  From \eqref{chi-def} and \eqref{varphi-def} we have
$$
c_{\chi,1} = \int_\R (1 + i \xi) \varphi(\xi)\ d\xi = \chi(0) - \frac{d}{dx}( e^x \chi(x) )|_{x=0}
$$
and the claim follows.  Now we handle the case $a = 2$.  We have
$$ c_{\chi,2} = \int_\R \int_\R \frac{(1 + i \xi) (1 + i \xi')}{2 + i (\xi + \xi')} \varphi(\xi) \varphi(\xi')\ d\xi d\xi'$$
Using the identity
$$ \frac{1}{2 + i (\xi + \xi')} = \int_0^\infty e^{-(1 + i \xi)x} e^{-(1 + i \xi')x}\ dx$$
we can rewrite $c_{\chi,2}$ as
$$ \int_0^\infty \left(\int_\R \varphi(\xi) (1 + i \xi) e^{-(1 + i \xi)x}\ d\xi\right)^2 dx.$$
But from differentiating \eqref{varphi-def} we see that the expression in parentheses is $-\chi'(x)$, and the claim follows.

Finally, for general $a$, we observe that since $\chi$ is real, we have $\varphi(-\xi) = \overline{\varphi(\xi)}$.  Taking complex conjugates of \eqref{chi-def} and substituting $\xi_j \mapsto -\xi_j$ we obtain the claim.
\end{proof}

Roughly speaking, we will be able to show the analogue of the generalised Hardy-Littlewood conjecture for these sums $\Lambda_{\chi,R,a}$ so long as
$\chi$ is suitably smooth and $R$ is a sufficiently small power of $N$. More precisely, we prove the following.

\begin{theorem}[Goldston-Y{\i}ld{\i}r{\i}m estimate]\label{fundlemma}  Let $t,d, L$ be positive integers, let $N$ be a large positive integer as usual, and let $\Psi = (\psi_1,\ldots,\psi_t)$ be a system of affine-linear forms with $\|\Psi\|_N \leq L$. Assume that no two of the forms $\psi_i$ are rational multiples of one another. Let $ a = (a_1,\ldots,a_t) \in \N^t$ be a $t$-tuple of integers.
Let $K \subseteq [-N,N]^d$ be a convex body, and let $\chi_1,\ldots,\chi_t: \R \to \R$ be smooth, compactly supported functions.  Let $R = N^\gamma$, where $\gamma > 0$ is sufficiently small depending on $t,d,L,\chi$ and $a$.  Call a prime $p$ \emph{exceptional} if there exist two forms $\psi_i, \psi_j$ which are linearly dependent modulo $p$, and let $P_\Psi$ denote the set of all exceptional primes. Write $X := \sum_{p \in P_\Psi} p^{-1/2}$.  Then we have
\begin{equation}\label{dickson-eq-fund}
\sum_{ n \in K \cap \Z^d} \prod_{i \in [t]} \Lambda_{\chi_i,R,a_i}( \psi_i( n) ) = \prod_{i \in [t]} c_{\chi_i,a_i}\cdot\vol_d(K) \cdot \prod_p \beta_p + O(\frac{N^d}{\log^{1/20} R} e^{O(X)} ), 
\end{equation}
where the local factors $\beta_p$ for each prime $p$ were defined in \eqref{beta-inf}, \eqref{beta-p}, and the sieve factors $c_{\chi,a}$ were defined in Definition \ref{sieve-factor}.
The implied constants here can depend on $t,d,L,\chi_1,\ldots,\chi_t$ and ${a}$.
\end{theorem}

\begin{remarks} Note that we are \emph{not} assuming that the system $\Psi$ has finite complexity but, as stated, we do assume that no two of the forms $\psi_i$ are rational multiples of one another.
This means that $P_{\Psi}$ is finite but not necessarily bounded in terms of $t,d,L$. If, for example, we have $d = 1$, $t = 2$ and $\psi_1(n) = n$, $\psi_2(n) = n+ M$, then $P_{\Psi}$ can be somewhat large if $M$ has many prime factors. If $\Psi$ does have finite complexity then $X$ is bounded in terms of $t,d,L$ and the error term becomes $o(N^d)$.  In other situations this term can be more substantial. We have not attempted to find an error term which is best possible, being happy to settle for one that suffices for our application, and in particular for the correlation condition (Definition \ref{correlation-condition}).  
\end{remarks}

Theorem \ref{fundlemma} should be compared with Conjecture \ref{dickson}.  The space $\Psi^{-1}((\R^+)^t)$, which appears in \eqref{beta-inf}, is not present here because the truncated divisor sums $\Lambda_{\chi_i,R,a_i}$ extend periodically to the negative numbers, in contrast to the von Mangoldt function $\Lambda$.

\begin{remark} In the works of Goldston, Pintz and Y{\i}ld{\i}r{\i}m \cite{gy-1,gy-2,gy-3,pintz} the choice of cutoff $\chi$ was critically important. In our analysis it is not, ultimately because the inverse Gowers-norm conjecture $\GI(s)$ applies even for arbitrarily small $\delta > 0$. This allows us to use simpler and smoother enveloping sieves in which the sieve factors are large. We do, of course, require these factors to be independent of $N$. In taking $\chi$ to be very smooth, a number of simplifications are possible. Following notes of the second author \cite{tao-coates,tao-gy-notes} (see also \cite{host-survey}), we avoid the use of any deep facts from analytic number theory such as the classical zero-free region for the Riemann zeta function. One  may instead make do with the elementary observation that the Riemann zeta function $\zeta(s)$ has the asymptotic $\zeta(s) = \frac{1}{s-1}+O(1)$ for $s$ near $1$ and $\Re(s) > 1$. We note that these simplifications could also be applied (retrospectively) to Chapters 9 and 10 of \cite{green-tao-longprimeaps}.
\end{remark}

\begin{remark} Observe that if $R = N^\gamma$, then $0 \leq \Lambda'(n) \leq \frac{1}{\gamma \chi(0)^2} \Lambda_{\chi,R,2}(n)$ for all $n$, $R < n \leq N$.  Thus we can use Theorem \ref{fundlemma} to obtain upper bounds for the expression \eqref{dickson-eq} which
lose a multiplicative factor of $\left(\frac{ c_{\chi,2} }{ \gamma \chi(0)^2 }\right)^t$, which is independent of $N$. This observation, coupled with a good choice of $\chi$ and $\gamma$, is rather close to the Selberg $\Lambda^2$ sieving technique. As is well-known there are significant barriers (the ``parity problem'') to reducing this multiplicative loss to something approaching 1.
\end{remark}

\textsc{Proof of Theorem \ref{fundlemma}.}  To simplify the notation we allow all implicit constants to depend on $t,d,L,\chi_1,\ldots,\chi_t$ and $a$. We may assume that $N$ (and hence $R$) are large with respect to these parameters, as the claim is trivial otherwise.

It is convenient to introduce the index set
$$ \Omega := \{ (i,j): i \in [t]; j \in [a_i] \} \subseteq \N^2.$$
With this notation, it is a simple matter to expand the left-hand side of \eqref{dickson-eq-fund} as
$$ \log^t R
\sum_{(m_{i,j})_{(i,j) \in \Omega} \in \N^\Omega} \bigg(\prod_{(i,j) \in \Omega} \mu(m_{i,j}) \chi_i\big(\frac{\log m_{i,j}}{\log R}\big)\bigg)
\sum_{ n \in K \cap \Z^d} \prod_{(i,j) \in \Omega} 1_{m_{i,j} | \psi_i(  n ) }.$$
The $\mu$ factors allow us to restrict $m_{i,j}$ to $\N_*$, the set of square-free natural numbers.
If, for each $i \in [t]$, we set $m_i := \lcm( m_{i,1},\ldots,m_{i,a_i} )$, then we can rewrite the above expression as
$$ \log^t R
\sum_{(m_{i,j})_{(i,j) \in \Omega} \in \N_*^\Omega} \bigg(\prod_{(i,j) \in \Omega} \mu(m_{i,j}) \chi_i\big(\frac{\log m_{i,j}}{\log R}\big)\bigg)
\sum_{ n \in K \cap \Z^d} \prod_{i \in [t]} 1_{m_i | \psi_i(  n ) }.$$
Since $\chi$ is compactly supported we may restrict $m_i$ to be at most $R^{O(1)}$ for all $i$.
In particular if we set $m := \prod_{i \in [t]} m_i$ then $m \leq R^{O(1)}$ also.
From the Chinese remainder theorem we see that as a function of $ n$, the expression $\prod_{i \in [t]} 1_{m_i | \psi_i(  n ) }$ is periodic with respect to the lattice $m \cdot \Z^d$. By a volume packing argument similar to that used
to prove \eqref{beta-gauss} in Appendix \ref{convexgeom}, we have
$$
\sum_{ n \in K \cap \Z^d} \prod_{i \in [t]} 1_{m_i | \psi_i(  n ) }
= \vol_d(K) \alpha_{m_1,\ldots,m_t}
+ O( m N^{d-1} )$$
where $\alpha_{m_1,\ldots,m_t}$ is the local factor
$$ \alpha_{m_1,\ldots,m_t} := \E_{ n \in \Z_m^d} \prod_{i \in [t]} 1_{m_i | \psi_i(  n ) }.$$
The total contribution of the error term $O(m N^{d-1})$ to \eqref{dickson-eq-fund} can be estimated crudely by
$O( R^{O(1)} N^{d-1} \log^t R )$, which will be $o(N^d)$ if the exponent $\gamma$ that defines $R$ is sufficiently small.  Thus we can discard this term and reduce our task to that of showing that
\begin{align}\nonumber
\log^t R \!\!\!\!
\sum_{(m_{i,j})_{(i,j) \in \Omega} \in \N_*^\Omega} \bigg(\prod_{(i,j) \in \Omega} \mu(m_{i,j}) \chi_i\big(\frac{\log m_{i,j}}{\log R}\big)\bigg)&
\alpha_{m_1,\ldots,m_t}
= \prod_{i \in [t]} c_{\chi_i,a_i}\cdot \prod_p \beta_p \\ &+ O( e^{O(X)} \log^{-1/20} R ).\label{logger}
\end{align}
Note that we have eliminated the convex body $K$ and the scale parameter\footnote{Observe that, although the $o$-notation concerns the situation when $N \rightarrow \infty$, this is exactly the same as letting $R \rightarrow \infty$.} $N$.   From the Chinese remainder theorem we make the key observation that $\alpha_{m_1,\ldots,m_t}$ is multiplicative in $m_1,\ldots,m_t$, so that if we decompose $m_i = \prod_p p^{r_{p,i}}$
then
\begin{equation}\label{mult}
\alpha_{m_1,\ldots,m_t} = \prod_p \alpha_{p^{r_{p,1}}, \ldots, p^{r_{p,t}}}.
\end{equation}
Note that as the $m_{i,j}$ are square-free, the $r_{p,i}$ are either $0$ or $1$.  

The next step is to use Fourier expansion to replace the weights $\chi_i$ by more multiplicative functions.  Indeed, as $\chi_i$ is smooth and compactly supported we have the Fourier expansion \eqref{varphi-def}
for some smooth $\varphi_i$ which is rapidly decreasing in the sense that $|\varphi_i(\xi)| \ll_A (1+\xi)^{-A}$ for all $A > 0$.  Thus we have
$$ \chi_i\big(\frac{\log m_{i,j}}{\log R}\big) = \int_{-\infty}^\infty m_{i,j}^{-\frac{1+i\xi}{\log R}} \varphi_i(\xi)\ d\xi.$$
We could insert this Fourier expansion into \eqref{logger} directly, but it will be easier if we first take advantage of the rapid decrease of $\varphi_i$ to truncate the Fourier integral to the interval $I := \{ \xi \in \R: |\xi| \leq \log^{1/2} R \}$ (say), thereby obtaining
$$ \chi_i\big(\frac{\log m_{i,j}}{\log R}\big) = \int_I m_{i,j}^{-\frac{1+i\xi}{\log R}} \varphi_i(\xi)\ d\xi + O_A( m_{i,j}^{-1/\log R} \log^{-A} R )$$
for any $A$.  Since $\chi_i(\log m_{i,j}/\log R)$ is itself bounded by $O(m_{i,j}^{-1/\log R})$, we conclude that
\begin{equation}\label{D-star}\prod_{(i,j) \in \Omega} \chi_i\big(\frac{\log m_{i,j}}{\log R}\big) 
= 
\int_I \ldots \int_I \prod_{(i,j) \in \Omega} m_{i,j}^{-z_{i,j}} \varphi_i(\xi_{i,j})\ d\xi_{i,j} 
+ O_A\big( \log^{-A} R \prod_{(i,j) \in \Omega} m_{i,j}^{-1/\log R}\big),\end{equation}
where we have written $z_{i,j} := (1+i\xi_{i,j})/\log R$.
Let us first deal with the contribution of the error term $O_A( \log^{-A} R \prod_{(i,j) \in \Omega} m_{i,j}^{-1/\log R})$
to \eqref{logger}.  Taking absolute values everywhere, we can bound these contributions by
$$ \ll_A (\log R)^{O(1)-A} \sum_{(m_{i,j})_{(i,j) \in \Omega} \in \N^\Omega} \alpha_{m_1,\ldots,m_t} \prod_{(i,j) \in \Omega} m_{i,j}^{-1/\log R}.$$
Using the multiplicativity, we can factorise this expression as an Euler product
$$ (\log R)^{O(1)-A}  \prod_{p} \sum_{(r_{i,j})_{(i,j) \in \Omega} \in \N^\Omega}
\alpha_{p^{r_1},\ldots,p^{r_t}} p^{-(\sum_{(i,j) \in \Omega} r_{i,j}) / \log R}$$
where $r_i := \max(r_{i,1},\ldots,r_{i,a_i})$.  Crude computations
then show that $\alpha_{p^{r_1},\ldots,p^{r_t}}$ is equal to $1$ when $r_1=\ldots=r_t=0$ and $O(1/p)$ otherwise (cf. Lemma \ref{sing}) and hence we can bound the above expression by
$$ (\log R)^{O(1)-A} \prod_{p} (1 - \frac{1}{p^{1 + 1/\log R}} )^{-O(1)}.$$
Since the Riemann zeta function $\zeta(s) = \prod_p (1 - \frac{1}{p^s})^{-1}$ has a simple pole at $s=1$ with residue $1$, we see that
\begin{equation}\label{prod}
\prod_p (1 - \frac{1}{p^s})^{-1} = \frac{1}{s-1} + O(1)
\end{equation}
whenever $\Re(s) > 1$ and $s-1$ is sufficiently close to $1$. This allows us to bound the above expression by
$$ O_A( (\log R)^{O(1)-A} )$$
which will be acceptable if $A$ is large enough.  Thus we only need to deal with the contribution of the main term of \eqref{D-star} to \eqref{logger}.  After swapping sums and integrals\footnote{One can justify the exchange of integrals and summations because $I$ is compact, and the summation can be shown to be absolutely convergent, either by using the crude bounds above, or by using bounds such as \eqref{qual} below.}, we
write this term as
$$
\log^t R
\int_I \ldots \int_I \sum_{(m_{i,j})_{(i,j) \in \Omega} \in \N_*^\Omega} \prod_{(i,j) \in \Omega} \mu(m_{i,j})
m_{i,j}^{-z_{i,j}} \alpha_{m_1,\ldots,m_t} \varphi_i(\xi_{i,j})\ d\xi_{i,j}.$$
Using the multiplicativity of $\alpha$ once more, we can write this expression as
\begin{equation}\label{euler}
\log^t R
\int_I \ldots \int_I  \prod_p E_{p,\xi}\cdot \prod_{(i,j) \in \Omega} \varphi_i(\xi_{i,j})\ d\xi_{i,j},
\end{equation}
where $\xi = (\xi_{i,j})_{(i,j) \in \Omega} \in I^\Omega$ and $E_{p,\xi}$ is the Euler factor
$$
E_{p,\xi} := \sum_{(m_{i,j})_{(i,j) \in \Omega} \in \{1,p\}^\Omega}  \big( \prod_{(i,j) \in \Omega} \mu(m_{i,j})
m_{i,j}^{-z_{i,j}} \alpha_{m_1,\ldots,m_t}\big).$$
Our task is to show that \eqref{euler} is $\big(\prod_{i \in [t]} c_{\chi_i,a_i}\big) \prod_p \beta_p + O( e^{O(X)} \log^{-1/20} R )$.
To tackle this we must understand the Euler factors $E_{p,\xi}$.  We may rewrite this expression as
\begin{equation}\label{ep-def}
E_{p,\xi} = \sum_{B \subseteq \Omega} (-1)^{|B|} \frac{\alpha(p,B)}{p^{\sum_{(i,j) \in B} z_{i,j}}}.\end{equation}
In this expression $\alpha(p,B) := \alpha_{p^{r_1},\ldots,p^{r_t}}$, where $r_i := 1$ whenever $(i,j) \in B$ for at least one $j$, and $r_i := 0$ otherwise. Note that $\alpha(p,\emptyset) = 1$.

Call a set $B \subseteq \Omega$ \emph{vertical} if it is non-empty and contained inside a vertical fibre $\{i\} \times [a_i]$ for some $i \in [t]$.  If $B$ is vertical then $\alpha(p,B) = \E_{{n} \in \Z_p^d}1_{p | \psi_i({n})}$, which is equal to $1/p$ if $p \geq p_0(t,d,L)$ is sufficiently large. To say something about $\alpha(p,B)$ when $B$ is neither empty nor vertical, recall that we described a prime $p$ as exceptional, and wrote $p \in P_{\Psi}$, if there exist $i,i'$ such that $\psi_i$ is a multiple of $\psi_{i'}$ in $\Z_p$.
For $p \notin P_{\Psi}$, we see from Lemma \ref{sing} that $\alpha(p,B) = O(1/p^2)$ whenever $B$ is not vertical or empty. If $p \in P_{\psi}$ then the best we can say in general is that $\alpha(p,B) = O(1/p)$.

From the above discussion we have 
\begin{equation}\label{qual}
 E_{p,\xi} = (1 + O(1/p^2)) E'_{p,\xi}\quad \hbox{ for } p \notin P_{\Psi},
\end{equation}
where $E'_{p,\xi}$ is the Euler factor\footnote{To provide a link to the discussion of \cite{green-tao-longprimeaps}, we observe that 
\[ \prod_p E'_{p,\xi} = \prod_{B \subseteq \Omega, B \vertical} \zeta\big(1 + \sum_{(i,j) \in B} z_{i,j} \big)^{(-1)^{|B|}}.\]}
\begin{equation}\label{ep-dash-def} E'_{p,\xi} := \prod_{B \subseteq \Omega, B \vertical} \left(1 - \frac{1}{p^{1 + \sum_{(i,j) \in B} z_{i,j}}}\right)^{(-1)^{|B|-1}}.\end{equation}
For $p \in P_{\Psi}$ we must rely instead on the far weaker bound
\begin{equation}\label{p-excep}
 E_{p,\xi} = (1 + O(1/p)) E'_{p,\xi}.
\end{equation}
From the estimate \eqref{prod} and the fact that $|z_{i,j}| = O(\log^{-1/2} R)$ when $\xi_{i,j} \in I$ we have 
\begin{align}\nonumber
\prod_p E'_{p,\xi} & = \prod_{B \subseteq \Omega, B \vertical} \bigg( \frac{1}{\sum_{(i,j) \in B} z_{i,j}} + O(1)  \bigg)^{(-1)^{|B|}}\\
 &= 
(1 + O(\log^{-1/2} R)) \prod_{B \subseteq \Omega, B \vertical} \big( \sum_{(i,j) \in B} z_{i,j} \big)^{(-1)^{|B|-1}}.\label{euler-conv}
\end{align}
Our aim now is to establish a corresponding estimate for $\prod_p E_{p,\xi}$. Note that we cannot afford the loss of a multiplicative constant which would result from a na\"{\i}ve application of \eqref{qual}. 

\begin{proposition}[Euler product estimate]\label{eulerprod}  We have
$$ \prod_p E_{p,\xi} = \left(\prod_p \beta_p + O( e^{O(X)} \log^{-1/20} R )\right)  \prod_p E'_{p,\xi}$$ 
for any $\xi \in I^\Omega$.
\end{proposition}

\begin{proof}   
From Lemma \ref{sing} we have $\beta_p = 1 + O(1/p)$ for $p \in P_\Psi$ and $\beta_p = 1 + O(1/p^2)$ otherwise. For starters this implies the very crude bound
\begin{equation}\label{v-crude} \prod_p \beta_p \leq e^{O(X)},\end{equation} which we will use later on.
Our first main task is to dispose of the contribution of the large primes $p$, when (say) $p > \log^{1/10} R$.
Using the estimates for $\beta_p$ just mentioned, we have
\begin{equation}\label{v-crude-ii} \prod_{p \leq \log^{1/10} R} \beta_p \leq e^{O(X)}.\end{equation}
We also have
\begin{align*}
\prod_{p > \log^{1/10} R} \beta_p &\leq \exp \big( O (\sum_{p > \log^{1/10} R : p \in P_{\Psi}} p^{-1})\big)\\
& \leq \exp\big( O(X\log^{-1/20} R) \big) \\
&= 1 + O( e^{O(X)} \log^{-1/20} R ),
\end{align*}
where the last bound follows from the elementary inequality $e^{\lambda X} \leq 1 + \lambda e^{X}$, valid for $\lambda \leq 1$ and $X \in \R_{\geq 0}$. Similarly, using the inequality $e^{-\lambda X} \geq 1 - \lambda e^X$, we obtain the corresponding lower bound, and thus
\begin{equation}\label{beta-crude} \prod_{p > \log^{1/10} R} \beta_p = 1 +  O( e^{O(X)} \log^{-1/20} R ). \end{equation}
From this and \eqref{v-crude} we see that it will suffice to show that
$$ \prod_p E_{p,\xi} = \bigg(\prod_{p \leq \log^{1/10} R} \beta_p + O( e^{O(X)} \log^{-1/20} R )\bigg) \prod_p E'_{p,\xi} .$$ 
Now from \eqref{qual}, \eqref{p-excep} we have
$$ \prod_{p > \log^{1/10} R} E_{p,\xi} = \exp\bigg( \sum_{p>\log^{1/10} R} p^{-2} +  \sum_{p \in P_\Psi: p>\log^{1/10} R} p^{-1} \bigg) \prod_{p > \log^{1/10} R} E'_{p,\xi}.$$
Since $\sum_{p>\log^{1/10} R} p^{-2} = O( \log^{-1/10} R)$ and 
$\sum_{p \in P_\Psi: p>\log^{1/10} R} p^{-1} = O( X \log^{-1/20} R )$, we conclude that
\begin{align*} \prod_{p > \log^{1/10} R} E_{p,\xi} &= \exp \big(O(1+X) \log^{-1/20} R\big) \prod_{p > \log^{1/10} R} E'_{p,\xi}\\ &= \big(1 + O(e^{O(X)} \log^{-1/20} R)\big) \prod_{p > \log^{1/10} R}E'_{p, \xi},\end{align*}the last step following as in the proof of \eqref{beta-crude}. From this and \eqref{v-crude} we see that it suffices to show that
\begin{equation}\label{to-prove-11}
\prod_{p \leq \log^{1/10} R} 
E_{p,\xi} = \bigg(\prod_{p \leq \log^{1/10} R} \beta_p + O( e^{O(X)} \log^{-1/20} R )\bigg) \prod_{p \leq \log^{1/10} R} E'_{p,\xi}.
\end{equation}
To do this, we will prove the following lemma.
\begin{lemma}\label{taylor-exp}
We have 
\[ E_{p,\xi} = \big(\beta_p + O( \frac{\log p}{\log^{1/2} R} ) \big) E'_{p,\xi}.\]
for all $p \leq \log^{1/10} R$.
\end{lemma}
\emph{Proof that Lemma \ref{taylor-exp} implies \eqref{to-prove-11}.} Suppose first that there is $p_0 \leq \log^{1/10} R$ such that $\beta_{p_0} = 0$. Then, using the fact that $\beta_p = 1 + O(1/p)$, we have
\[ \prod_{p \leq \log^{1/10} R} \big(\beta_p + O( \frac{\log p}{\log^{1/2} R} ) \big) = O(\frac{\log p}{\log^{1/2} R}) e^{O(X)},\] which is acceptable. If no $\beta_p$ is vanishes then, since $\beta_p = 1 + O(1/p)$ and $\beta_p$ is a rational with denominator dividing $p^d$, we have a bound $\beta_p \gg 1$ with the implied constant depending only on the global parameters $t,d,L$. Thus, using \eqref{v-crude-ii}, we have
\begin{align*} \prod_{p \leq \log^{1/10} R} \big(\beta_p + O( \frac{\log p}{\log^{1/2} R} ) \big) & = \prod_{p \leq \log^{1/10} R} \beta_p \cdot \prod_{p \leq \log^{1/10} R} \big(1 + O(\frac{\log p}{\log^{1/2} R})  \big)\\
& = \big(\prod_{p \leq \log^{1/10} R} \beta_p\big) \cdot \big( 1 + O(\log^{-1/3} R))\\
& = \prod_{p \leq \log^{1/10} R} \beta_p + O(e^{O(X)} \log^{-1/3} R)).\end{align*}
Thus \eqref{to-prove-11} holds in this case also.\endproof

\emph{Proof of Lemma \ref{taylor-exp}.} Observe that since $\xi \in I^\Omega$, we have \[ p^{\sum_{(i,j) \in B} z_{i,j}} = 1 + O( \log p / \log^{1/2} R )\] for 
all $B$ and all $p \leq \log^{1/10} R$.  Dividing \eqref{ep-def} by \eqref{ep-dash-def} (noting that the latter has magnitude comparable to $1$) and performing Taylor expansion in $w = p^{\sum z_{i,j}}$ about $w = 1$ it is not hard to check that
$$\frac{E_{p,\xi}}{E'_{p,\xi}} = \frac{\widetilde E_p}{\widetilde E'_p} + O( \frac{\log p}{\log^{1/2} R} ),
$$
where $\widetilde E_p, \widetilde E'_p$ are defined setting all the $z_{i,j}$ equal to zero in \eqref{ep-def} and \eqref{ep-dash-def} respectively. Thus
\begin{equation}\label{ep-tilde-def}
\widetilde E_p := \sum_{B \subseteq \Omega} (-1)^{|B|} \alpha(p,B)
\end{equation}
and
\begin{equation}\label{ep-tilde-dash-def}
\widetilde E'_p := \sum_{B \subseteq \Omega, B \vertical} \big(1 - \frac{1}{p}\big)^{(-1)^{|B| - 1}}.
\end{equation}
To prove the lemma, then, it suffices
to prove the identity
\begin{equation}\label{betaform}
\beta_p = \frac{\widetilde{E}_p}{\widetilde E'_p}.
\end{equation}
Recalling \eqref{ep-tilde-def} and \eqref{ep-tilde-dash-def}, it will suffice to show that
\begin{equation}\label{beta-ident}
\sum_{B \subseteq \Omega} (-1)^{|B|} \alpha(p,B) = \beta_p \prod_{B \subseteq \Omega, B \vertical} (1 - \frac{1}{p})^{(-1)^{|B|-1}}.
\end{equation}
Using the binomial theorem, the right-hand side of \eqref{beta-ident} simplifies to $\beta_p (1 - \frac{1}{p})^{t}$, which by \eqref{beta-p}
is equal to
$$ \E_{ n \in \Z_p^d} 1_{p \not | \psi_1( n)} \ldots 1_{p \not | \psi_t( n)}.$$
By the inclusion-exclusion principle this can be written as
$$ \sum_{r_1,\ldots,r_t \in \{0,1\}} (-1)^{r_1+\ldots+r_t} \E_{ n \in \Z_p^d} \prod_{r_i = 0} 1_{p | \psi_i( n)},$$
which in turn is just
$$ \sum_{r_1,\ldots,r_t \in \{0,1\}} (-1)^{r_1+\ldots+r_t} \alpha_{p^{r_1},\ldots,p^{r_t}}.$$
We are to show that this is equal to the left-hand side of \eqref{beta-ident}, namely
\[ \sum_{B \subseteq \Omega} (-1)^{|B|} \alpha(p,B).\]
To do this, we compare coefficients of $\alpha_{p^{r_1},\dots,p^{r_t}}$ on both sides. To evaluate the coefficient on the left-hand side, let $I$ be the set of indices for which $r_i \neq 0$. Then this coefficient is easily seen to be
\[ \prod_{i \in I} \sum_{B_i \subseteq [a_i]} (-1)^{|B_i|} \]
which, by the binomial theorem, is simply $(-1)^{|I|}$.  This gives \eqref{betaform}, and the claim follows.
\end{proof}

We return to the proof of \eqref{logger}. Recall that we had reduced this to the task of finding an approproate asymptotic for \eqref{euler}. Substituting the result of Proposition \ref{eulerprod} into \eqref{euler} and applying \eqref{v-crude}, it is easy to reduce this in turn to showing the following two facts.
Firstly, that
\begin{equation}\label{euler-1}
\log^t R
\int_I \ldots \int_I \big(\prod_p E'_{p,\xi}\big) \prod_{(i,j) \in \Omega} \varphi_i(\xi_{i,j})\ d\xi_{i,j}
= \prod_{i \in [t]} c_{\chi_i,a_i} + O(\log^{-1/20}R);
\end{equation}
and secondly that
\begin{equation}\label{euler-2}
\log^t R
\int_I \ldots \int_I 
\prod_p |E'_{p,\xi}| \prod_{(i,j) \in \Omega} |\varphi_i(\xi_{i,j})|\ d\xi_{i,j}
= O(1).
\end{equation}

Let us begin with the second task, that of proving \eqref{euler-2}. We simply substitute $z_{i,j} = (1 + i\xi_{i,j})/\log R$ into \eqref{euler-conv}. The contribution from the terms $\log R$ is precisely $\log^{-t} R$, by a simple application of the binomial theorem  $\sum_{B \subseteq C: B \neq \emptyset} (-1)^{|B|} = 0^{|C|} - 1$. For the terms involving the $\xi_{i,j}$ we have the crude estimate
\[ \prod_{(i,j) \in \Omega} (1 + |\xi_{i,j}|)^{O(1)},\] and so
$$ \prod_p |E'_{p,\xi}| \ll \frac{1}{\log^t R} \prod_{(i,j) \in \Omega} (1 + |\xi_{i,j}|)^{O(1)}.$$
However since $\chi$ is smooth its modified Fourier transform satisfies $|\varphi_i(\xi_{i,j})| \ll_A (1 + |\xi_{i,j}|)^{-A}$ for any $A > 0$, as we have already remarked.  The claim then follows by taking $A$ large enough.

Now we prove \eqref{euler-1}. Using the rapid decay of the functions $\varphi$ once more together with \eqref{euler-conv} we see that it suffices to show that
\begin{align*}
\log^t R
\int_I \ldots \int_I \prod_{B \subseteq \Omega, B \vertical} \big( \sum_{(i,j) \in B} z_{i,j} \big)^{(-1)^{|B|-1}} 
& \prod_{(i,j) \in \Omega} \varphi_i(\xi_{i,j})\ d\xi_{i,j}
 \\ & = \prod_{i \in [t]} c_{\chi_i,a_i} + O(\log^{-1/20}R).
\end{align*}
The first move is to reinstate the integrals over all of $\R$, rather than just over $I$. Doing this introduces an error which is $\ll_A \log^{-A} R$ for any $A > 0$, on account of the rapid decrease of $\varphi$.  Once this is done the multiple integral is easily seen to factor, there being one integral for each index $i$. After scaling out the factors of $\log R$, the claim follows from the definition \eqref{chi-def} of the sieve weights $c_{\chi,a}$.
The result follows, and we have concluded the proof of Theorem \ref{fundlemma}.
\endproof

\textsc{Construction of the enveloping sieve.} Now we are ready to prove Proposition \ref{pseudodom}, the statement of which was as follows.
\begin{pseudodom-repeat}[Domination by a pseudorandom measure] 
Let $D > 1$ be arbitrary. Then there is a constant $C_0 := C_0(D)$ such that the following is true. Let $C \geq C_0$, and suppose that $N' \in [CN,2CN]$.  Let
$b_1,\ldots,b_t \in \{0,1,\ldots,W-1\}$ be coprime to $W := \prod_{p \leq w} p$.  Then there exists a $D$-pseudorandom 
measure $\nu: \Z_{N'} \to \R^+$ which obeys the pointwise bounds
$$ 1 + \Lambda'_{b_1,W}(n) + \ldots + \Lambda'_{b_t,W}(n) \ll_{D,C} \nu(n) $$
for all $n \in [N^{3/5}, N]$, where we identify $n$ with an element of $\Z_{N'}$ in the obvious manner.
\end{pseudodom-repeat}
The definition of $D$-pseudorandom was given, and discussed, in \S \ref{envelope-sec}. See in particular Definitions \ref{linear-forms-condition} and \ref{correlation-condition} and the paragraphs following the latter.
Let $\gamma = \gamma(C,D) > 0$ be a parameter to be chosen later and set $R := N^\gamma$.  Fix an arbitrary smooth even function $\chi: \R \to \R$ which is supported on $[-1,1]$ and satisfies  $\chi(0)=1$ and $\int_0^1 |\chi'(x)|^2\ dx = 1$.  For such a function we have $c_{\chi,2} = 1$, thanks to Lemma \ref{sieve-compute}.

We define the preliminary weight $\tilde \nu: [N] \to \R^+$ by setting
$$ \tilde \nu(n) := \E_{i \in [t]} \frac{\phi(W)}{W} \Lambda_{\chi,R,2}(Wn+b_i)$$
and then transfer this to $\Z_{N'}$ by setting $\nu(n) := \frac{1}{2} + \frac{1}{2} \tilde \nu(n)$ when $n \in [N]$ 
and $\nu(n) := 1$ otherwise.  

By construction, $\tilde \nu$ is certainly non-negative.  To verify the pointwise bounds, it suffices to show that
$$ \Lambda'_{b_i,W}(n) \ll_{C,D} \frac{\phi(W)}{W} \Lambda_{\chi,R,2}(Wn+b_i)$$
for all $i \in [t]$ and $n \in [N^{3/5}, N]$.  The left-hand side is only non-zero when $Wn+b_i$ is a prime
which is greater than $N^{3/5}$. Supposing that $\gamma < 3/5$, we see that in this case the left-hand side is equal to $\frac{\phi(W)}{W} \log N$, while the right-hand side is $\frac{\phi(W)}{W} \log R$.  Since $R = N^\gamma$ and $\gamma$ depends only on $C,D$, the claim follows.

It remains to show that $\nu$ is a $D$-pseudorandom measure.  Our argument here shall follow that in \cite{green-tao-longprimeaps} rather closely, but will use Theorem \ref{fundlemma} as a substitute for \cite[Propositions 9.5,9.6]{green-tao-longprimeaps}. For that reason we shall skip some of the details which are more or less exact repetition of those in \cite{green-tao-longprimeaps}.

Let us first verify the $(D,D,D)$-linear forms condition.  By decomposing $\nu$ up into its various components as in \cite{green-tao-longprimeaps}, it certainly suffices to establish the somewhat general bound
$$ \sum_{ n \in K \cap \Z^d} \prod_{j \in [m]} \tilde \nu( \psi_j(  n ) ) 
= \vol_d(K) + o(N^d)$$
where $\Psi = (\psi_1,\ldots,\psi_m)$ is a system of affine-linear forms, no two of which are affinely related, $m,d, \|\Psi\|_N$ are all $O_D(1)$, and $K \subseteq [-N,N]^d$ is a convex body with $\Psi(K) \subseteq [-N,N]^m$.  Splitting $\tilde \nu$ up further,
we thus reduce to showing that
\begin{equation}\label{phew}
 \left(\frac{\phi(W)}{W}\right)^{m} \sum_{ n \in K \cap \Z^d} \prod_{j \in [m]} \Lambda_{\chi,R,2}( \psi_j( Wn+b_{i_j} ) )
= \vol_d(K) + o(N^d)
\end{equation}
for all $i_1,\ldots,i_m \in [t]$.  

Now we apply Theorem \ref{fundlemma}.  As we are assuming that no two of the forms $\psi_j(n)$ are affinely related, the same is true for the forms $\psi_j(Wn + b_{i_j})$.  In particular we see that the exceptional primes, if they exist, are bounded in size by $O(w) = O(\log\log N)$. In particular we have $X = O(\log \log^{1/2} N)$ and so $e^{O(X)} \log^{-1/20} R = o(1)$.
We can thus write the left-hand side of \eqref{phew} as
$$ \left(\frac{\phi(W)}{W}\right)^{m}  c_{\chi, 2}^{m} \vol_d(K) \prod_p \beta_p + o(N^d)$$
where we suppress the dependence of constants on $t,m,d,L,D$.  Because all the $b_{i_j}$ are coprime to $W$, we see that
$\beta_p = (\frac{p}{p-1})^t$ for all $p \leq w$, and in particular $\prod_{p \leq w} \beta_p = \left(\frac{W}{\phi(W)}\right)^t$.
Also, for $p > w$ we see from Lemma \ref{sing} that $\beta_p = 1 + O(1/p^2)$, and so $\prod_{p > w} \beta_p = 1 + o(1)$.  Since $c_{\chi,2}=1$, the claim follows.

Now we verify the $D$-correlation condition for $\nu$.  As before we can pass from $\nu$ to $\tilde \nu$, and reduce to showing that
$$ \sum_{n \in I} \prod_{j \in [m]} \tilde \nu(n+h_j) \ll N \sum_{1 \leq j < j' \leq m} \tau(h_j - h_{j'})$$
for all $m = O_D(1)$, all $h_1,\ldots,h_m \in [N]$, and all intervals $I \subseteq [N]$, and where $\tau: [-N,N] \to \R^+$ obeys the moment bounds $\E_{n \in [-N,N]} \tau(n)^q \ll_q 1$ for all $q > 0$.  We may assume that no two of the $h_i$ are equal as in this case one can use crude divisor estimates, setting $\tau(0)$ to be moderately large 
(see \cite{green-tao-longprimeaps} for details).  Again, we split up $\tilde \nu$ and reduce to showing that
$$ \left(\frac{\phi(W)}{W}\right)^{m}
\bigg(\sum_{n \in I} \prod_{j \in [m]} \Lambda_{\chi,R,2}(W(n+h_j)+b_{i_j})\bigg) \ll N \sum_{1 \leq j < j' \leq m} \tau(h_j - h_{j'})$$
whenever $i_1,\ldots,i_m \in [t]$.  We can apply Theorem \ref{fundlemma} with the system of forms $\Psi = (W(n + h_j) + b_{i_j})_{j = 1}^m$ and write the left-hand side as
$$ \left(\frac{\phi(W)}{W}\right)^{m} \bigg(c_{\chi,2}^m |I| \prod_p \beta_p + O(N e^{O(X)} \log^{-1/20} R) \bigg) .$$
As before we can discard the sieve factor $c_{\chi,2}=1$, and we have $\prod_{p \leq w} \beta_p = \left(\frac{W}{\phi(W)}\right)^m$. 

It thus suffices to show that
$$ \prod_{p > w} \beta_p +  e^{O(X)} \log^{-1/20} R \ll \sum_{1 \leq j < j' \leq m} 
\tau(h_j - h_{j'}).$$

From Lemma \ref{sing} we see that for $p > w$ we have $\beta_p = 1+O(1/p)$, with the improvement $\beta_p = 1 + O(1/p^2)$
as long as $p \notin P_{\Psi}$, that is as long as $p$ does not divide $W(h_j-h_{j'})+b_{i_j} - b_{i_{j'}}$ for any $1 \leq j < j' \leq m$.  Thus
\begin{equation}
\prod_{p > w} \beta_p  \ll \prod_{\substack{p > w \\ p \in P_{\Psi}}} \big(1 + O(\frac{1}{p})\big) \ll \exp \bigg( O\big( \sum_{\substack{p > w \\ p \in P_{\Psi}}}\frac{1}{p} \big) \bigg) .\label{use-soon-7}
\end{equation}
On the other hand, since $w = O(\log \log N)$ is so small we have
\begin{align*}
 e^{O(X)} \log^{-1/20} R 
&\leq \exp( O( \sum_{p \in P_\Psi} \frac{1}{p^{1/2}} ) ) \log^{-1/20} R \\
&\ll  \exp \bigg( O\big( \sum_{\substack{p > w \\ p \in P_{\Psi}}}\frac{1}{p^{1/2}} \big) \bigg).
\end{align*}
It follows from this analysis that if we set
$$ \tau(n) := \sum_{1 \leq j < j' \leq m} \exp\bigg( O\big( \sum_{\substack{p > w \\ p | Wn+b_{i_j} - b_{i_{j'}}}}  \frac{1}{p^{1/2}}\big) \bigg)$$
then we obtain the desired correlation estimate.  To show the moment bounds on $\tau$ it suffices to show that
$$ \E_{n \in [N]} \exp( q \sum_{\substack{p > w\\ p | Wn + h}} \frac{1}{p^{1/2}} ) \ll_q 1$$
for all $h = O(W)$. By repeating the proof of
\cite[Lemma 9.9]{green-tao-longprimeaps} we can deduce this bound from
$$ \sum_{n \in [N]} \prod_{\substack{p > w\\ p | Wn + h}} (1 + p^{-1/4}) \ll_q N.$$
Using the bound
$$ \prod_{\substack{p > w \\ p | Wn+h}} (1 + p^{-1/4}) \leq \sum_{\substack{(d,W) = 1 \\ d | Wn+h }} d^{-1/4}$$
we reduce to showing that
$$ \sum_{\substack{(d, W) = 1 \\ d = O(WN)}} d^{-1/4} \sum_{\substack{n \in [N] \\ d | Wn+h}} 1 \ll_q N.$$
But we have \[ \sum_{\substack{n \in [N] \\ d | Wn+h}} 1 = O(1 + N/d)\] by the Chinese remainder theorem, and the claim then follows easily.
This concludes the proof of Proposition \ref{pseudodom}.
\endproof
  
\textsc{The correlation estimate for $\Lambda^\sharp$.} 
The final task of this appendix is to establish the correlation estimate \eqref{to-prove-appD}, which was the crucial fact that $\Lambda^{\sharp}_{b,W} - 1$ has small Gowers norm. We allow all constants to depend on $s$.
Expanding out the $U^{s+1}[N]$ norm, it suffices to show the slightly more general bound
$$ \sum_{(n,  h) \in K}
\prod_{\omega \in \{0,1\}^{s+1}} (\frac{\phi(W)}{W} \Lambda^\sharp(W(n+\omega \cdot  h)+b) - 1)
= o(N^{s+2})
$$
whenever $K$ is a convex body in $[-N,N]^{s+2}$.  Expanding out the product, it suffices to show that
$$ \sum_{(n,  h) \in K}
\prod_{\omega \in B} \frac{\phi(W)}{W} \Lambda^\sharp(W(n+\omega \cdot  h)+b)
= \vol_{s+2}(K) + o(N^{s+2})$$
for all $B \subseteq \{0,1\}^{s+1}$.  Now observe that $\Lambda^\sharp = - \Lambda_{\chi^\sharp,R,1}$, and so we may invoke
Theorem \ref{fundlemma} with the system of forms $\Psi = (W(n + \omega \cdot {h}) + b)_{\omega \in B}$ to write the left-hand side as
$$
\big(- \frac{\phi(W)}{W}\big)^{|B|}
c_{\chi^\sharp,1}^{|B|} \vol_{s+2}(K) \prod_p \beta_p + O(N^{s+2} e^{O(X)} / \log^{1/20} R ).$$
As in the preceding section, we compute that $\prod_{p \leq w} \beta_p = \left(\frac{W}{\phi(W)}\right)^{|B|}$, while
$\beta_p = 1 + O(1/p^2)$ for $p > w$. Furthermore all exceptional primes $p$ have $p \leq w$, and 
thus since $w$ is so small
\[  e^{O(X)} / \log^{1/20} R = O(\log^{-1/20} R) \exp\big(O ( \sum_{p \leq w} \frac{1}{p^{1/2}}) \big) = o(1).\]
Finally, from Lemma \ref{sieve-compute}, we have
$c_{\chi^\sharp,1} = -1$.  The claim follows.
\endproof

\section{Nilmanifold constraints; Host-Kra cube groups}\label{nil-app}

Our aim in this appendix is prove Proposition \ref{cube}, which asserts a constraint concerning parallelepiped in nilmanifolds.  It 
turns out to be convenient to generalise the notion of a parallelepiped to a more general object, namely a \emph{Host-Kra cube}.  Thus much of this appendix will be devoted to the algebraic theory of these cubes.  We will first introduce such parallelepipeds in the Lie group $G$, establish the constraint there, and then descend to the quotient space $G/\Gamma$ and show that the constraint persists down to the quotient. In preparing the material that follows we benefitted much from conversations with Sasha Leibman, and also from remarks made by one of the anonymous referees.

\textsc{Host-Kra cube groups in $G$.} Let $G$ be a connected Lie group with identity $\id_G$, with the associated lower central series $G_{\bullet}$ given by
$$ G = G_0 = G_1 \supseteq G_2 \supseteq \ldots,$$ where $G_0 = G_1 = G$ and $G_{i+1} = [G, G_i]$.  We recall the standard facts that
$[G_i,G_j] \subseteq  G_{i+j}$, and that each $G_i$ is a closed connected normal Lie subgroup of $G$; see
for instance \cite[Ch. 3, \S 9, Corollary to Prop. 4]{bourbaki}.  In particular the quotient
groups $G_i \backslash G$ are also Lie groups.

To define the Host-Kra cube group $\HK^{s+1}(G_{\bullet})$ we first need some combinatorial notation.

\begin{definition}[Simple combinatorics of $\{0,1\}^{s+1}$]
We refer to $\{0,1\}^{s+1}$ as \emph{the} cube. Its elements $\omega$ may be partially ordered by decreeing that $\omega \leq \omega'$ if $\omega_j \leq \omega'_{j}$ for $j = 1,\dots,s+1$. A \emph{hyperplane} is any set of the form $H_{j,a} := \{\omega : \omega_j = a\}$. If $0 \leq d \leq s+1$ then we say that a \emph{face of codimension $d$} is any non-empty intersection $F$ of $d$ distinct hyperplanes, and we write $d = \codim(F)$. Thus any vertex in $\{0,1\}^{s+1}$ is a face of codimension $s+1$, whilst the whole cube $\{0,1\}^{s+1}$ is a face of codimension $0$. We say that two faces are \emph{parallel} if they have the same fixed coordinates, and hence the same codimension. Every face $F$ has a minimal element $\min(F)$ and a maximal element $\max(F)$.  We say that a face\footnote{With respect to the partial ordering $\leq$, a lower face is exactly the same concept as a principal filter.} is \emph{lower} if $\min(F) = 0^{s+1}$. Note that every face is parallel to exactly one lower face, and that lower faces $F$ are in one-to-one correspondence with their maximal elements $\max(F)$, which can be arbitrary. Finally, we say that two parallel faces are \emph{adjacent} if their union is a face of one lower codimension.
\end{definition}
 
\begin{definition}[Face groups]
Let $F \subseteq  \{0,1\}^{s+1}$ be an face of codimension $d$. For any element $g \in G$, we write $g^F$ for the element of $G^{\{0,1\}^{s+1}}$ such that $(g^F)_{\omega} = g$ when $\omega \in F$, and $(g^F)_{\omega} = \id_G$ otherwise.
The \emph{face group} $\Gamma_F$ is the group generated by all elements $g^F$ with $g \in G_{\codim(F)}$, thus
$\Gamma_F \cong G_{\codim(F)}$.
\end{definition}

\begin{definition}[Host-Kra cube group]
The Host-Kra cube group $\HK^{s+1}(G_{\bullet})$ is the subgroup of $G^{\{0,1\}^{s+1}}$ generated by all the face groups $\Gamma_F$, as $F$ ranges over faces of $\{0,1\}^{s+1}$.
\end{definition}

The Host-Kra cube group could be defined with a more general \emph{filtration} in place of the lower central series $G_{\bullet}$, that is to say a sequence of subgroups in which the condition that $G_{i+1} = [G, G_{i}]$ is relaxed to an inclusion $[G_i, G_j] \subseteq G_{i+j}$. We will not need this here.

The significance of the group $\HK^{s+1}(G_{\bullet})$ for us is that it contains the parallelepipeds:
 
\begin{lemma}[Parallelepipeds are Host-Kra cubes]\label{cubes-are-parallels}
Given any $g,x \in G$ and $n,h_1,\ldots,$ $h_{s+1}$ in $\Z$, the parallelepiped
$\g := (g^{n + \omega \cdot {h}} x)_{\omega \in \{0,1\}^{s+1}}$ lies in $\HK^{s+1}(G_{\bullet})$.
\end{lemma}

\begin{proof} We may write, in $G^{\{0,1\}^{s+1}}$,
\[ \g = (g^{h_{s+1}})^{F_{s+1}}(g^{h_s})^{F_s} \dots (g^{h_1})^{F_1}(g^nx)^{F_0},\]
where $F_0 := \{0,1\}^{s+1}$ and $F_i$ is the hyperplane $F_i := \{\omega : \omega_i = 1\}$ for $i = 1,\dots,s+1$. Thus $\g$ is the product of $s+2$ of the generators of $\HK^{s+1}(G_{\bullet})$.
\end{proof}

The face groups $G_F$ are related to each other in a pleasant way:

\begin{lemma}[Face relations]\label{face-relations-lem}
Let $F,F'$ be faces in $\{0,1\}^{s+1}$. 
\begin{enumerate}
\item If $F, F'$ are disjoint, then the elements in $\Gamma_F$ and $\Gamma_{F'}$ commute with one another.
\item If $F$ and $F'$ intersect then $[\Gamma_{F}, \Gamma_{F'}] \subseteq  \Gamma_{F \cap F'}$. 
\item If $F$ and $F'$ are adjacent and parallel, then $\Gamma_F \subseteq  \Gamma_{F'} \Gamma_{F \cup F'}$ and $\Gamma_{F'} \subseteq  \Gamma_F \Gamma_{F \cup F'}$.
\end{enumerate}
\end{lemma}

\begin{proof}
 (i) is immediate. To prove (ii), note that any element of $[\Gamma_F, \Gamma_F']$ has the form $x^{F \cap F'}$ for some $x \in [G_d,G_{d'}]$, where $d := \codim(F)$ and $d' := \codim(F')$. The result follows on noting that $\codim(F \cap F') \leq \codim(F) + \codim(F')$, and recalling from group theory that $[G_d, G_{d'}] \subseteq  G_{d + d'}$. (iii) is immediate; in this setting we have $x^F x^{F'} = x^{F \cup F'}$.
\end{proof}

From Lemma \ref{face-relations-lem} (iii) and an easy induction on the codimension we see that every face group $\Gamma_F$ lies in the group generated by the \emph{lower} face groups.  In particular this implies that the entire group $\HK^{s+1}(G_{\bullet})$ is generated by
the lower face groups. The same result holds for the \emph{upper faces}, but we will not have any further use of this here.

Now we seek a more explicit description of $\HK^{s+1}(G_{\bullet})$ by the lower face groups.  To achieve this, we need

\begin{definition}[Decreasing ordering of faces]
Let $F_1 > \dots > F_{2^{s+1}}$ be any ordering of the $2^{s+1}$ lower faces of $\{0,1\}^{s+1}$. We say that this ordering is \emph{decreasing} if whenever $F_i \supseteq F_j$ we have $i \geq j$. Thus $F_1 = \{0,1\}^{s+1}$ and $F_{2^{s+1}} = 0^{s+1}$. 
\end{definition}

Clearly, decreasing orders of faces exist; let us fix such an ordering.  Now, observe from Lemma \ref{face-relations-lem} (i),(ii) that if $i < j$, then we either have $\Gamma_{F_j} \cdot \Gamma_{F_i} \subseteq \Gamma_{F_i} \cdot \Gamma_{F_j}$ or $\Gamma_{F_j} \cdot \Gamma_{F_i} \subseteq \Gamma_{F_i} \cdot \Gamma_{F_j} \cdot \Gamma_{F_k}$ for some $k > j$.  From these inclusions we see that any product of elements from the lower face groups $\Gamma_{F_i}$ can eventually be contained in $\Gamma_{F_1} \cdot \Gamma_{F_2} \cdot \ldots \cdot \Gamma_{F_{2^{s+1}}}$, as one can use the above inclusions to move all occurrences of $\Gamma_{F_1}$ to the far left, use the closure property $\Gamma_{F_1} \cdot \Gamma_{F_1} = \Gamma_{F_1}$ to concatenate, then move all occurences of $\Gamma_{F_2}$ to be adjacent to $\Gamma_{F_1}$, and so forth.  Since the lower face groups generate $\HK^{s+1}(G_{\bullet})$, we have thus obtained the factorisation
$$
 \HK^{s+1}(G_{\bullet}) = \Gamma_{F_1} \cdot \Gamma_{F_2} \cdot \ldots \cdot \Gamma_{F_{2^{s+1}}}.
$$
Thus there exist functions $\tau_i: \HK^{s+1}(G_{\bullet}) \to \Gamma_{F_i}$ such that
\begin{equation}\label{hps}
\g = \tau_1(\g) \ldots \tau_{2^{s+1}}(\g)
\end{equation}
for all $\g \in \HK^{s+1}(G_{\bullet})$. 

\begin{remark} Since $\Gamma_{F_i} \cong G_{\codim(F_i)}$ is a closed connected Lie subgroup of $G^{\{0,1\}^{s+1}}$, we can conclude from the above factorisation that $\HK^{s+1}(G_{\bullet})$ is also a closed connected group Lie subgroup of $G^{\{0,1\}^{s+1}}$.  Furthermore, since the hyperplane face groups consist entirely of parallelepipeds, and the lower dimensional face groups can be expressed as commutators of the hyperplane face groups, we see that $\HK^{s+1}(G_{\bullet})$ is in fact the subgroup generated by the parallelepipeds.  Thus this is an extremely natural group for studying parallelepipeds.
\end{remark}

In the factorisation \eqref{hps}, the $\tau_i$ are unique: an inspection of the $\max(F_1)$ coefficients of both sides shows that $\tau_1(\g)$ is determined uniquely by $\g$, and then after factoring $\tau_1(\g)$ out, an inspection of the $\max(F_2)$ coefficients of both sides shows that $\tau_2(\g)$ is determined uniquely by $\g$, and so forth.  Indeed, this algorithm shows that if $\g = (g_\omega)_{\omega \in \{0,1\}^{s+1}}$, then $\tau_i(\g) \in \Gamma_{F_i}$ is a continuous function of the coordinates $g_{\max(F_1)}, \ldots, g_{\max(F_i)}$ only; indeed, equating $\Gamma_{F_i}$ with $G_{\codim(F_i)}$, the group element $\tau_i(\g)$ is an explicit word in these coordinates. Conversely, $g_{\max(F_i)}$ is a word in $\tau_1(\g),\dots,\tau_i(\g)$ only.

Recall that we are aiming to prove Proposition \ref{cube}, which establishes a constraint amongst the $2^{s+1}$ vertices of a parallelepiped in $G/\Gamma$, an $s$-step nilmanifold. Henceforth we assume that we are in this setting (the discussion up to now has been valid quite generally). The preceding observations allow one to prove a related fact, namely that if $G$ is $s$-step nilpotent and if $(g_{\omega})_{\omega \in \{0,1\}^{s+1}} \in \HK^{s+1}(G_{\bullet})$ then $g_{0^{s+1}}$ is a word in the $g_{\omega}$, $\omega \in \{0,1\}^{s+1}_*$. Indeed the nilpotence of $G$ implies that the final face group $\Gamma_{F_{2^{s+1}}}$ is trivial, and hence $\tau_{2^{s+1}}(\g) = \id$ for all $\g$. Thus $g_{0^{s+1}} = g_{\max(F_{2^{s+1}})}$ is a actually a word in $\tau_1(\g),\dots,\tau_{2^{s+1}-1}(\g)$, and hence in the $g_{\omega}$, $\omega \in \{0,1\}_*^{s+1}$.

To prove Proposition \ref{cube}, we must show how this constraint ``descends'' to $G/\Gamma$. A step in this direction is the following lemma, which follows immediately from the fact that $g_{0^{s+1}}$ is a word in the $g_{\omega}$, $\omega \in \{0,1\}^{s+1}_*$.

\begin{lemma}\label{gamm-lem}
Suppose that $g = (g_{\omega})_{\omega \in \{0,1\}^{s+1}} \in \HK^{s+1}(G_{\bullet})$ and that $g_{\omega} \in \Gamma$ for all $\omega \in \{0,1\}^{s+1}_*$. Then the remaining point $g_{0^{s+1}}$ lies in $\Gamma$ as well.
\end{lemma}

We have defined the Host-Kra cube \emph{group}; now we define the Host-Kra nilmanifold.

\begin{definition}[Host-Kra cube nilmanifold] We define the Host-Kra nilmanifold $\HK^{s+1}(G_{\bullet}/\Gamma)$ to be the $s$-step nilmanifold $\HK^{s+1}(G_{\bullet})/(\Gamma^{\{0,1\}^{s+1}} \cap \HK^{s+1}(G_{\bullet}))$. 
\end{definition}

\emph{A priori}, this definition does not make sense. The Lie group $\HK^{s+1}(G_{\bullet})$ is connected, simply-connected and $s$-step nilpotent Lie group (the nilpotence follows from the fact that it is a subgroup of $G^{\{0,1\}^{s+1}}$ and the simple-connectedness from the factorisation \eqref{hps} together with the simple-connectedness of the face groups $\Gamma_{F_i} \cong G_{\codim(F_i)}$). We have not, however, shown that $\Gamma^{\{0,1\}^{s+1}} \cap \HK^{s+1}(G_{\bullet})$ is cocompact inside it. This is the business of Lemma \ref{lemmaC} below. To prove it, we will need a basic topological property of nilmanifolds, first established in the foundational paper of Mal'cev \cite{malcev}.

\begin{lemma}\label{malcev-lemma}\cite{malcev}
Let $G$ be a connected, simply-connected nilpotent Lie group, and let $\Gamma$ be a discrete cocompact subgroup. Then for any $j \geq 1$ the group $\Gamma \cap G_j$ is discrete and cocompact in $G_j$.
\end{lemma}

\begin{remark} To obtain results such as the Main Theorem in the case $s = 2$, we need only consider nilmanifolds which are products of Heisenberg examples. This was observed in Proposition \ref{gi2-prop}. In this case, Lemma \ref{malcev-lemma} can easily be verified by hand using calculations along the lines of those in \cite[Appendix B]{green-tao-u3mobius}.
\end{remark}

\begin{lemma}\label{lemmaC} $\Gamma^{\{0,1\}^{s+1}} \cap \HK^{s+1}(G_{\bullet})$ is a discrete and cocompact subgroup of $\HK^{s+1}(G_{\bullet})$.
\end{lemma}
\proof The discreteness is obvious, since $\Gamma^{\{0,1\}^{s+1}}$ is discrete in $G^{\{0,1\}^{s+1}}$. Now by Lemma \ref{malcev-lemma} there is a a compact set $K_j \subseteq G_j$ such that $G_j = K_j \cap (\Gamma \cap G_j)$. For each $i$, consider the subgroup $H_i \leq \HK^{s+1}(G_{\bullet})$ consisting of those $\g$ such that $\tau_1(\g) = \dots = \tau_i(\g) = \id$. By our earlier observations this is the same as the subgroup $\{\g : g_{\max(F_1)} = \dots = \g_{\max(F_i)} = \id\}$, and hence in particular $H_i$ is normal in $\HK^{s+1}(G_{\bullet})$.

Suppose that $1 \leq i \leq 2^{s+1}$ and that $\g \in H_{i-1}$. Then $\g = \tau_{i}(\g) \mathbf{h}$, where $\mathbf{h} \in H_{i}$. Since $\tau_i(\g)$ lies in the fact group $\Gamma_{F_i}$, we may write it as $(k_i \gamma_i)^{F_i}$ where $k_i \in K_{\codim(F_i)}$ and $\gamma_i \in \Gamma \cap G_{\codim(F_i)}$. Since $H_i$ is normal, we may hence write
\[ \g = (k_i)^{F_i} \cdot \mathbf{h}' \cdot (\gamma_i)^{F_i},\] where $\mathbf{h'}$ is another element of $H_i$. 

Continuing inductively until $i = 2^{s+1}$, we eventually express an arbitrary element of $\HK^{s+1}(G_{\bullet})$ as a product of $(k_1)^{F_1} \dots (k_{2^{s+1}})^{F_{2^{s+1}}}$ times an element of $\Gamma^{\{0,1\}^{s+1}} \cap \HK^{s+1}(G_{\bullet})$. Since the set
\[ \{ (k_1)^{F_1} \dots (k_{2^{s+1}})^{F_{2^{s+1}}} : k_1 \in  K_1,\dots, k_{2^{s+1}} \in K_{2^{s+1}}\}\] is compact, this proves the lemma.\endproof

\emph{Proof of Proposition \ref{cube}.} The projection $G^{\{0,1\}^{s+1}} \rightarrow (G/\Gamma)^{\{0,1\}^{s+1}}$ induces a 1-1 continuous map from the compact set $\HK^{s+1}(G_{\bullet}/\Gamma)$ to $(G/\Gamma)^{\{0,1\}^{s+1}}$. Henceforth, we consider the former set as a a compact subset of the latter. Let $p$ be the restriction to $\HK^{s+1}(G_{\bullet}/\Gamma)$ of the obvious projection from $(G/\Gamma)^{\{0,1\}^{s+1}}$ to $(G/\Gamma)^{\{0,1\}_*^{s+1}}$, and let $\Sigma$ be the range of this map. It follows from Lemma \ref{gamm-lem} that this map is 1-1, and hence there is a unique map $P : \Sigma \rightarrow G/\Gamma$ such that $(P(x),x) \in \HK^{s+1}(G_{\bullet}/\Gamma)$ for every $x = (x_{\omega})_{\omega \in \{0,1\}_*^{s+1}} \in \Sigma$. The map $P$ is automatically continuous since its graph $\HK^{s+1}(G_{\bullet}/\Gamma)$ is compact and all spaces involved are Hausdorff.

Proposition \ref{cube} follows immediately from Lemma \ref{cubes-are-parallels}.\endproof

\end{document}